\documentclass[10pt]{article}
\usepackage[utf8]{inputenc}
\usepackage{fancyhdr}
\usepackage{extramarks}
\usepackage{amsmath}
\usepackage{amsthm}
\usepackage{amsfonts}
\usepackage{siunitx} 
\usepackage{tikz}
\usepackage[plain]{algorithm}
\usepackage{algpseudocode}
\usepackage{multirow}
\usepackage{booktabs}
\usepackage{palatino}
\usepackage{graphicx}
\usepackage{subcaption}
\usepackage[colorlinks,linkcolor=black,anchorcolor=black,citecolor=black,urlcolor=blue]{hyperref}
\usepackage{amsmath,bm}
\usepackage{booktabs}
\usepackage{mathtools}
\usepackage{amssymb}
\usepackage{caption}
\usepackage{capt-of}
\usepackage{mciteplus}
\usepackage{cite}
\usepackage{mathrsfs}
\usepackage[title,titletoc,toc]{appendix}
\usepackage{xr}
\usepackage{parskip}
\usepackage{textcomp}
\usepackage[colaction]{multicol}
\usepackage[switch]{lineno}
\usepackage{lipsum}
\usepackage{etoolbox}
\usepackage{longtable}
\usepackage{array}
\usepackage{rotating}
\usepackage{tablefootnote}
\graphicspath{ {./images/} }
\usepackage{epstopdf}
\usepackage{relsize}
\usepackage{float}
\newcolumntype{P}[1]{>{\centering\arraybackslash}p{#1}}
\newcolumntype{M}[1]{>{\centering\arraybackslash}m{#1}}
\newcolumntype{C}[1]{>{\centering\arraybackslash}p{#1}}
\usetikzlibrary{automata,positioning}

\usetikzlibrary{automata,positioning}

\begin{document}

\title{Geometric algebra generation of molecular surfaces}
\author{Azzam Alfarraj$^{1,2}$ and Guo-Wei Wei$^{1,3,4}$\footnote{
		Corresponding author.		E-mail: weig@msu.edu} \\
	$^1$ Department of Mathematics, \\
	Michigan State University, MI 48824, USA.\\
	$^2$ Mathematics Department,\\
	King Fahd University of Petroleum and Minerals, Dhahran 31261, KSA. \\
	$^3$ Department of Electrical and Computer Engineering,\\
	Michigan State University, MI 48824, USA. \\
	$^4$ Department of Biochemistry and Molecular Biology,\\
	Michigan State University, MI 48824, USA. \\
}

\date{} 

\maketitle

\begin{abstract}
Geometric algebra is a powerful framework that unifies mathematics and physics. Since its revival in the middle of the 1960s by David Hestenes, it attracts great attention and has been exploited in many fields such as physics, computer science, and engineering. This work introduces a geometric algebra method for the molecular surface generation that utilizes the Clifford-Fourier transform which is a generalization of the classical Fourier transform. Notably, the classical Fourier transform and Clifford-Fourier transform differ in the derivative property in $\mathbb{R}_k$ for $k$ even. This distinction is due to the noncommutativity of geometric product of pseudoscalars with multivectors and has significant consequences in applications. We use the Clifford-Fourier transform in $\mathbb{R}_3$ to benefit from the derivative property in solving partial differential equations (PDEs). The Clifford-Fourier transform is used to solve the mode decomposition process in PDE transform. Two different initial cases are proposed to make the initial shapes used in the present method. The proposed method is applied first to small molecules and proteins.   To validate the method, the molecular surfaces generated are compared to surfaces of other definitions. Applications are considered to protein electrostatic analysis. 
This work opens the door for further applications of geometric algebra and Clifford-Fourier transform in biological sciences. 
 
\end{abstract}

\pagenumbering{roman}
\begin{verbatim}
\end{verbatim}

{\setcounter{tocdepth}{4} \tableofcontents}
\newpage


\setcounter{page}{1}
\renewcommand{\thepage}{{\arabic{page}}}

\section{Introduction}\label{sec:introduction}

The structures of biomolecules, such as those of proteins,  DNAs,  molecular motor, subcellular organelles, and viruses are directly related to their interactions and functions \cite{giard2010molecular,zheng2012biomolecular,bates2008minimal}. Therefore, studying biomolecular structures is a major topic in molecular biology and is essential for understanding biological processes. For example, the visualization of biomolecular surfaces and their electrostatic potentials is vital in the analysis of biomolecular interactions like protein-nucleic acid and protein-protein interactions, ligand-receptor binding, macromolecular assembly, enzymatic mechanism, and drug discovery \cite{zheng2012biomolecular}. Also, geometric complementarity, which is essential in molecular docking can be predicted via molecular surface shapes \cite{giard2010molecular}. Hence, many theoretical methods and computational algorithms were proposed to characterize molecular shapes. In 1953, Corey and Pauling presented molecular models of atom and bond that are still pivotal in molecular sciences \cite{corey1953molecular}. The wide applications of surfaces raise the need for fast and reliable surface generation algorithms. This is due to the fact that in molecular simulations molecular surfaces are rendered millions of times repeatedly \cite{zheng2012biomolecular}.  For the case of large macromolecules that require excessive memory, a divide-and-conquer method was proposed to improve the efficiency of surface generation \cite{zhao2018divide}.

It is worth noting that when we talk about molecular surfaces we may not mean real physical surfaces. It is, indeed, a representation of molecular shapes and there are many representations available in the literature. The most popular representations are the van der Waals surface (VdWS), the solvent excluded surface (SES) and the solvent accessible surface (SAS) \cite{richards1977areas}. The VdWS is defined as a surface composed of overlapping rigid spheres with each sphere having a radius corresponding to the van der Waals radius of the corresponding atom. The SES is defined as the surface of the volume generated by moving a sphere representing a solvent molecule around the molecule and depicting the positions of the exterior surface of the sphere. Likewise, the SAS is depicting the positions of the center of the sphere \cite{giard2010molecular}. The importance of such representations is that they are often used as molecule-solvent interfaces. These interface models are crucial in illustrating how surfaces interact with surrounding molecules such as ions, counterions, and solvents \cite{zheng2012biomolecular}. These interactions may determine the stability and solubility of macromolecules in an aqueous environment. The essence of such investigations comes from the fact that the human cell mass has a great percentage of water in the range of $65-90\%$ and most biological processes occur in that aqueous part of the cell. With this being said, the above surface models have been exploited in studies of protein folding \cite{spolar1994coupling}, protein-protein interactions \cite{crowley2005cation}, drug classification \cite{bergstrom2003absorption}, solvation energies \cite{raschke2001quantification}, macromolecular docking \cite{dragan2004dna}, ion channel transport \cite{zheng2011poisson}, protein pocket detection \cite{zhao2018protein} and DNA binding  \cite{dragan2004dna}.

The aforementioned surfaces, i.e., VdWS, SES and SAS, admit geometric singularities which result in computational difficulties \cite{connolly1985depth, eisenhaber1993improved, gogonea1994implementation, sanner1996reduced,yu2008feature}. To overcome this issue, the energy minimization principle has been adopted for biomolecular surface construction. Partial different equation (PDE)-based biomolecular surfaces were proposed in 2005 \cite{wei2005molecular}. Inspired by geometric flows, the minimal molecular surface (MMS) was presented in 2006. In general,  minimal surfaces are widely seen in nature due to the energy minimization principle \cite{bates2008minimal}. These methods, unlike other popular methods, start with atomic coordinates and radii rather than some given surfaces. In addition,   Gaussian surfaces \cite{wang2021regularization, chen2011tmsmesh, decherchi2013general}, flexibility-rigidity index surfaces \cite{mu2017geometric, opron2014fast, xia2013multiscale}, level-set surfaces \cite{cheng2009coupling}, and skinning surfaces \cite{cheng2009quality} were also proposed to avoid singularity issues. 

In the past few decades, geometric flow algorithms were exploited in image analysis and surface processing. Witkin, in 1983, proposed an image denoising algorithm using diffusion equations that were presented to be formally equivalent to Gaussian low-pass filters \cite{witkin1984scale}.   Perona and Malik presented an anisotropic diffusion equation for image denoise without edges being smeared \cite{perona1990scale}. Generalized Perona-Malik equation with arbitrarily high order nonlinear PDEs was proposed for edge-preserving noisy image restoration  \cite{wei1999generalized}. Mode decomposition evolution equations were proposed to generalize nonlinear PDE-based high-pass filters. These equations perform a PDE transform, which splits the data, signals, and images into functional modes such as trend, edge, texture, noise, and so on, depending on frequencies \cite{wang2011partial}.  PDE transform was used to generate biomolecular surfaces \cite{zheng2012biomolecular}. The fast Fourier transform was incorporated in the PDE transform to avoid the stability constraints of solving high-order PDEs \cite{zheng2012biomolecular}. 

The Fourier transform is widely applied in science and engineering. Due to its great importance and impact on experiments and computational work, there have been many versions proposed to generalize or improve Fourier transform. In the field of geometric algebra, Clifford-Fourier transform was presented among other proposed transforms such as quaternion-Fourier transform \cite{brackx2005clifford,felsberg2001commutative, bayro2010geometric, hitzer2013quaternion}. The Clifford-Fourier transform exploits two main notions in geometric algebra: {\it geometric product} and {\it multivector}. It is noteworthy that the notion of multivectors is built on the geometric product that was proposed by Clifford in 1876 \cite{hestenes2012new}  to unify the work of Grassmann in the outer product, also known as wedge product, and the work of Hamilton in quaternions \cite{hestenes2012new}.  However, Clifford's work did not get much attention until the 1960s. In 1966, Davis Hestenes ignited the revival of geometric algebra and geometric calculus in his book \textit{Space-Time Algebra} \cite{bayro2010geometric}. In the beginning, the emphasis of Hestenes was mainly on physics before geometric algebra applications gotten wide recognition in other fields such as computer science \cite{dorst2010geometric} and image processing \cite{wang2019geometric,batard2010clifford}. Hestenes suggested geometric algebra to be the unifying language of mathematics and physics \cite{hestenes2012clifford}. Geometric algebra has also been applied to protein structure analysis  \cite{quine1999helix,hitzer2010interactive,chys2008application,lavor2019oriented,billinge2016assigned}.  Nowadays, one can say that geometric algebra offers a unified framework for diverse applications in mathematics, physics, computer science, engineering, and biology \cite{corrochano2001geometric}. Notably, other fields of mathematics that have a strong relationship with geometric algebra also have great potential in biophysical applications. Specifically, the evolutionary de Rham-Hodge method was proposed for molecular data representation and analysis \cite{chen2021evolutionary}. The wedge product and exterior calculus used in de Rham-Hodge theory are related to geometric calculus. Also, the $k$-forms of de Rham-Hodge theory play a very similar role as the $k$-vectors of geometric algebra. The evolutionary de Rham-Hodge method showed success in predicting the protein B-factors of some challenging cases and outperformed the present methods in protein flexibility analysis \cite{chen2021evolutionary}.   

The goal of this work is to develop a geometric algebra-based biomolecular surface generation algorithm. The Clifford-Fourier transform is used along with the PDE transform to define a new molecular representation. This work opens a new direction in geometric algebra-based biomolecular modeling and analysis. It may stimulate future applications of geometric algebra in biological sciences. This paper is organized into four sections. Section \ref{sec:introduction} is dedicated to a brief literature review on molecular surface generation methods, PDE transform, and geometric algebra, and calculus. Then, Section \ref{sec:Theories} presents the theoretical background of our method. It starts by giving a thorough introduction to geometric algebra stating the definitions and main properties of the outer product, geometric product, $k$-vectors, multivectors, and Clifford algebras. After that, Clifford-Fourier transform is presented where necessary definitions of multivector functions, derivatives, and integration are stated in the geometric calculus context. Then, specific cases of Clifford-Fourier transform in two-dimensional (2D) and 3D settings are discussed along with showing similarities and differences with the original Fourier transform. Next, we discuss the PDE transform. Afterward, Section \ref{sec:generation} is devoted to our biomolecular surface generation method. Two equations used in the construction of initial surfaces are given. Then, test cases are provided and investigated to explore the effects of changing the parameters, i.e., propagation time and isovalues. After investigating the parameters, surfaces of real proteins are generated and compared to those from well-known methods. Finally, Section \ref{sec:application} demonstrates some applications on the generation of the surface for the purposes of validation. First, the electrostatic surface potentials are calculated and mapped to surfaces generated using our Clifford-Fourier transform method and then to surfaces generated using the MSMS method \cite{sanner1996reduced}. The calculations are conducted using the APBS package in VMD \cite{jurrus2018improvements}. Second, the electrostatic solvation free energies of 21 proteins are calculated and compared to three other molecular surface generation methods presented in the literature. The energy calculations are carried out using MIBPB \cite{chen2011mibpb}.

\section{Theories and methods}\label{sec:Theories}

The Fourier transform has been used extensively in mathematics, science, and engineering. Many versions of the Fourier transform have been proposed in different fields of mathematics. In the field of geometric algebra, Clifford-Fourier transform and quaternion Fourier transform were presented. Since this work utilizes the Clifford-Fourier transform, a basic introduction to geometric algebra is given to define notions and establish notations. Then, a definition of Clifford-Fourier transform follows.

\subsection{Geometric algebra}
Geometric algebra presents a framework where operations are given as scalars, vectors, and multivectors irrespective of the grade of the vectors. This unification and generalization is due to two main concepts in geometric algebra:  \textbf{ geometric product} and \textbf{  multivector} \cite{ebling2005clifford,bhatti2020geometric}. To define the geometric product, we need first to introduce the outer product, also called the wedge product, operator of geometric algebra. 

\subsubsection{Outer product}
Given $x$ and $y$ in $\mathbb{R}^n$, their outer product is represented by $x \wedge y$. For any three vectors $x,y$ and $z$ in $\mathbb{R}^n$ and a scalar $\lambda$ in $\mathbb{R}$, the outer product has the following properties:

\begin{align} \label{p1}
x \wedge y &= - y \wedge x  \\
x \wedge x &= 0   \\
(\lambda x) \wedge y &= \lambda (x \wedge y)  \\
\lambda (x \wedge y) &=  (x \wedge y)\lambda \label{eq:comm property} \\
x \wedge (y+z) &= (x \wedge y) + (x \wedge z) \\
x \wedge (y \wedge z) &= (x \wedge y) \wedge z
\end{align}

\begin{figure}[ht]
	\centering
	\includegraphics[scale=0.7,trim={5cm 11cm 5cm 5cm},clip]{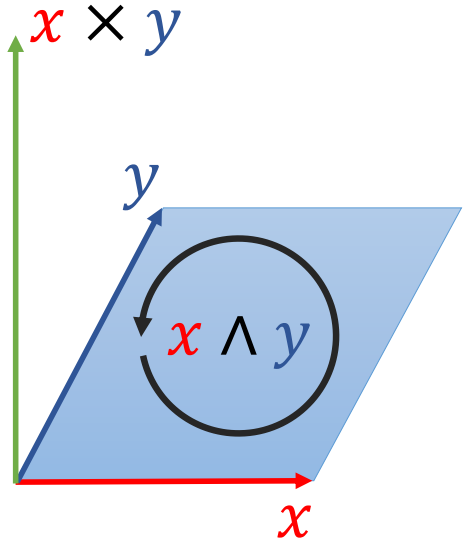}
	\caption{Visualization of \( x \wedge y \) compared to \( x \times y \). \( x \wedge y \) is an oriented parallelogram while \( x \times y\) is a vector perpendicular to both \( x \) and \( y \).}
	\label{fig:bivector}
\end{figure}

It is worth noting that the outer product is anticommutative as given in the property (\ref{p1}). The result of the outer product of two vectors $x$ and $y$ is called a \textbf{  bivector} and can be visualized as an oriented parallelogram with $x$ and $y$ as shown in Figure \ref{fig:bivector}.

Furthermore, the outer product of three vectors is called a  \textbf{  trivector} and can be visualized as an oriented parallelepiped, which has six oriented parallelogram as its faces. In general, the outer product of $k$ vectors is called a  \textbf{  $k$-vector} and obviously, it cannot be visualized in a 3D setting. Keep in mind that $k$-vectors are feasible if the vectors are in a Euclidean space of dimension $n$ where $k \leq n$. Otherwise, the outer product is zero if $k>n$. A $k$-vector is said to have a \textbf{  grade} $k$.
 

\subsubsection{Geometric product}
After this introduction of the outer product, the geometric product of $x$ and $y$ is represented as $xy$ and defined as:  
\begin{equation}
xy = x \cdot y + x \wedge y
\end{equation}
where $x \cdot y$ is the inner product.

Note that the product $xy$ is a combination of a scalar $x \cdot y$ and a bivector $x \wedge y$. This note leads us to the concept of \textbf{multivectors}. In an Euclidean space of dimension $n$, a \textbf{multivector} is a linear combination of different-grade $k$-vectors with $n$ being the highest grade, i.e. it is a linear combination of a scalar, vectors, bivectors, trivectors, ... , $(n-1)$-vectors and an $n$-vector. Let us deduce some properties of the geometric product. The geometric product of a vector $x$ with itself is the magnitude of $x$ squared as shown below:

\[
xx = x \cdot x + x \wedge x = ||x||^2 + 0 = ||x||^2,
\]
which leads to the fact that for any nonzero vector $x$, the vector $\dfrac{x}{||x||^2}$ is its inverse. Also, the inner product and the outer product of any two vectors $x$ and $y$ can be expressed in terms of their geometric products only as shown below:
\newline
Since \[xy = x \cdot y + x \wedge y \quad\text{and} \quad yx = x \cdot y - x \wedge y, \]
then \[xy + yx = 2 (x \cdot y) \quad \text{and} \quad xy - yx = 2 (x \wedge y), \]
which implied
\begin{align*}
x \cdot y &= \dfrac{xy + yx}{2},\\
x\wedge y &= \dfrac{xy - yx}{2}.
\end{align*}

\subsubsection{Clifford algebras $\mathbb{R}_n$}

First, we generalize the concept of basis to the geometric algebra setting. If 
$\lbrace e_i : i=1,...,n \rbrace$ 
is an orthonormal basis of $\mathbb{R}^n$, then 
\[ \lbrace e_i \wedge e_j : i,j=1,...,n \text{ and } i\neq j\rbrace \]
is a basis for bivectors and 
\[ \lbrace e_i \wedge e_j \wedge e_k : i,j,k=1,...,n \text{ and } i\neq j, j\neq k, i\neq k \rbrace \] is a basis for trivectors, and so on for the rest of $k$-vectors. This means that any $k$-vector can be written as a linear sum of the basis presented. This is useful when it comes to the summation of $k$-vectors and helps in finding out the resulting $k$-vector. 

Notably, the geometric product of orthonormal basis vectors has some special properties that simplify the calculations of geometric products of multivectors. These properties are:
 
\begin{align}
e_i^2 &= 1 , \qquad\qquad i = 1,2,\cdots,n \\
e_ie_j &= -e_je_i, \qquad\qquad i \neq j
\end{align}
since
\begin{align*}
e_i^2 &= e_i\cdot e_i + e_i\wedge e_i = e_i\cdot e_i = ||e_i||^2 = 1 , \qquad\qquad i = 1,2,\cdots,n \\
e_ie_j &= e_i\cdot e_j + e_i \wedge e_j = e_i \wedge e_j = -e_j \wedge e_i = -e_je_i, \qquad\qquad i \neq j.
\end{align*}

For the Euclidean space $\mathbb{R}^n$, we get a Clifford algebra $\mathbb{R}_n$ that has the dimension $2^n$. The basis of Clifford algebra $\mathbb{R}_n$ consists of the scalar $1$ and the basis of $\mathbb{R}^n$ and all different geometric products of the basis vectors. The Clifford algebra $\mathbb{R}_n$ contains all multivectors of grade $n$ or less where the multivector grade is the highest grade among its constituent $k$-vectors. For the sake of simplicity, we limit the discussion to $\mathbb{R}_2$ and $\mathbb{R}_3$ in the rest of this section. Therefore, the basis for $\mathbb{R}_2$ is
\[
\lbrace 1,e_1,e_2,e_3,e_1e_2 \rbrace
\]
and the basis for $\mathbb{R}_3$ is
\[
\lbrace 1,e_1,e_2,e_3,e_1e_2,e_2e_3,e_1e_3,e_1e_2e_3 \rbrace
\]
and to simplify the notations, $e_ie_j$ and $e_ie_je_k$ are denoted as $e_{ij}$ and $e_{ijk}$ from now on. The following are examples of multivectors in  $\mathbb{R}_2$ and  $\mathbb{R}_3$ respectively:
\begin{align*}
&M_2 = \alpha_0 + \alpha_{1}e_{1} + \alpha_{2}e_{2} + \alpha_{12}e_{12} \\
&M_3 = \alpha_0 + \alpha_{1}e_{1} + \alpha_{2}e_{2} + \alpha_{3}e_{3} + \alpha_{12}e_{12} + \alpha_{23}e_{23}+ \alpha_{13}e_{13} + \alpha_{123}e_{123}.
\end{align*}

It is noteworthy that the basis of the Clifford algebra $\mathbb{R}_2$ has only one element of grade 2 which is $e_{12}$ and therefore it is called \textbf{pseudoscalar} and denoted as $i_2$. In $\mathbb{R}_3$, similarly $e_{123}$ is called pseudoscalar and denoted as $i_3$.
The pseudoscalars have two important features. First, the square to -1 as follows
\begin{align}
i_2^2 &= e_1e_2e_1e_2 = -e_1e_1e_2e_2 = -(1)(1) = -1 \label{eq: comm property 2} \\
i_3^2 &= e_1e_2e_3e_1e_2e_3 = e_1e_2e_1e_2e_3e_3 = -e_1e_1e_2e_2e_3e_3 = -(1)(1)(1)=-1 .
\end{align}
Second, any multivector $M_2$ in $\mathbb{R}_2$ can be written as
\[
M_2 = \alpha + a +i_2(\beta +b)
\]
where $\alpha,\beta \in \mathbb{R}$ and $ a,b \in \mathbb{R}^2$, and any multivector $M_3$ in $\mathbb{R}_3$ can be written as
\[
M_3 = \alpha + a +i_3(\beta +b)
\]
where $\alpha,\beta \in \mathbb{R}$ and $ a,b \in \mathbb{R}^3$.

Moreover, $\lbrace \alpha +i_2 \beta | \alpha,\beta \in \mathbb{R}\rbrace \subset \mathbb{R}_2$ is isomorphic to $\mathbb{C}$ since $i_2^2=-1$. So, for any scalar $\gamma$
\[
e^{(i_2\gamma)} = \cos(\gamma) + i_2 \sin(\gamma).
\]
Likewise for $\mathbb{R}_3$ we have 
\[
e^{(i_3\gamma)} = \cos(\gamma) + i_3 \sin(\gamma).
\]
Note that $i_2$ is not commutative with multivectors in $\mathbb{R}_2$ which makes $e^{(i_2\gamma)}$ not commutative as well. Indeed, $i_2$ commutes with scalars and bivectors, as shown in Equation \ref{eq:comm property} and Equation \ref{eq: comm property 2},  and anticommutes with vectors as follows
\[
e_1i_2 = e_1e_1e_2 = -e_1e_2e_1 = -i_2e_1.
\]
The same can be said about $e_2$.

 On the other hand, $i_3$ is commutative with any multivector in $\mathbb{R}^3$ which means $e^{(i_3\gamma)}$ is commutative with any multivector in $\mathbb{R}^3$, and this can be proved in the way followed with $i_2$. This note makes a noticeable impact when discussing the properties of the Clifford-Fourier transform.

\subsection{Clifford-Fourier transform}
A multivector function $\mathbf{F}$ is a function whose range is a set of multivectors \cite{hestenes1968multivector}. Now, let $\mathbf{F}$ be a multivector function that is defined on $\mathbb{R}^n$

\[
\mathbf{F} : \mathbb{R}^n \longrightarrow \mathbb{R}_m,
\]
then its directional derivative in direction $r$ is defined as \cite{ebling2005clifford}
\[
\mathbf{F}_r =\lim_{\epsilon\rightarrow 0} \dfrac{\mathbf{F}(x+\epsilon r) - \mathbf{F}(x)}{\epsilon}
\],
where $\epsilon \in \mathbb{R}$. Also, its Riemannian integral is defined as \cite{ebling2005clifford}
\[
\int_{\mathbb{R}^n} \mathbf{F}(x) |dx| = \lim_{\substack{ n\rightarrow \infty \\ |\Delta x_i| \rightarrow 0}} \sum_{i=1}^{n}\mathbf{F}(x_i)\Delta x_i
\].

The Clifford Fourier transform  is presented for a 2D setting and then for a 3D setting. 

\subsubsection{Clifford-Fourier transform in 2D}

The Clifford-Fourier transform of a multivector function $\mathbf{F}: \mathbb{R}^n \longrightarrow \mathbb{R}_2$ is defined as 
\[
\mathcal{F}\lbrace \mathbf{F} \rbrace (u) = \int_{\mathbb{R}^n} \mathbf{F}(x)  e^{(-2 \pi i_2 \langle x,u \rangle)} |dx|
\]
and the inverse Clifford-Fourier transform is defined as
\[
\mathcal{F}^{-1}\lbrace \mathbf{F} \rbrace (x) = \int_{\mathbb{R}^n} \mathbf{F}(u)  e^{(2 \pi i_2 \langle x,u \rangle)} |du|
\],
where $x,u \in \mathbb{R}^n$, provided the integrals exist.
A multivector function $\mathbf{F}$ defined as 
\[
\mathbf{F}(x) = \mathbf{F}_0(x) + \mathbf{F}_1(x)e_1 +\mathbf{F}_2(x)e_2 +\mathbf{F}_{12}(x)e_{12}
\]
can be written as 
\begin{align*}
\mathbf{F}(x) &= \left[\mathbf{F}_0(x) +\mathbf{F}_{12}(x)i_2\right] \\
			&+ e_1\left[\mathbf{F}_1(x) +\mathbf{F}_2(x)i_2\right],
\end{align*}
which can be seen as two complex signals and interpreted as an element of $\mathbb{C}^2$.
The linearity of Clifford-Fourier transform would result in the following
\begin{align*}
\mathcal{F}\lbrace \mathbf{F} \rbrace (u) 
&= \left[\mathcal{F}\lbrace \mathbf{F}_0 + \mathbf{F}_{12}i_2 \rbrace (u)\right]1 \\
&+ e_1\left[\mathcal{F}\lbrace \mathbf{F}_1 + \mathbf{F}_2i_2 \rbrace (u)\right],
\end{align*}
which means that the Clifford-Fourier transform in 2D can be dealt with as a linear combination of two classical Fourier transforms.
One of the most powerful properties of the Fourier transform is the derivative property.  Fortunately, the Clifford-Fourier transform has derivative properties that might agree or disagree with the ones of the classical Fourier transform. To present the derivative property for a multivector function $\mathbf{F}$, one needs to decompose it into
\[
\mathbf{F} = \left[\mathbf{F}_0(x)+\mathbf{F}_{12}(x)e_{12}\right] 
+ \left[\mathbf{F}_1(x)e_1 +\mathbf{F}_2(x)e_2\right] = \mathbf{\textit{f}} + \mathbf{f},
\]
where $\mathbf{\textit{f }}$ is commutative with $i_2$ and $\mathbf{f}$ in anticommutative. With this being said, Ebling and Scheuermann \cite{ebling2005clifford} showed the following
\begin{align*}
\mathcal{F}\lbrace \nabla\mathbf{f} \rbrace (u) &= -2\pi i_2 u \mathcal{F}\lbrace \mathbf{f} \rbrace (u) \\
\mathcal{F}\lbrace \Delta\mathbf{f} \rbrace (u) &= 4\pi^2 u^2 \mathcal{F}\lbrace \mathbf{f} \rbrace (u),
\end{align*}
while,
\begin{align*}
\mathcal{F}\lbrace \nabla\mathbf{\textit{f}} \rbrace (u) &= 2\pi i_2 u \mathcal{F}\lbrace \mathbf{\textit{f}} \rbrace (u) \\
\mathcal{F}\lbrace \Delta\mathbf{\textit{f}} \rbrace (u) &= 4\pi^2 u^2 \mathcal{F}\lbrace \mathbf{\textit{f}} \rbrace (u).
\end{align*}
From the above, one can see that
\[
\mathcal{F}\lbrace \Delta\mathbf{F} \rbrace (u) = 4\pi^2 u^2 \mathcal{F}\lbrace \mathbf{F} \rbrace (u)
\]
while no rule can be written for $\mathcal{F}\lbrace \nabla\mathbf{F}\rbrace$.

\subsubsection{Clifford-Fourier transform in 3D}

The Clifford-Fourier transform of a multivector function $\mathbf{F}: \mathbb{R}^n \longrightarrow \mathbb{R}_3$ is defined as 
\[
\mathcal{F}\lbrace \mathbf{F} \rbrace (u) = \int_{\mathbb{R}^n} \mathbf{F}(x)  e^{(-2 \pi i_3 \langle x,u \rangle)} |dx|
\]
and the inverse Clifford-Fourier transform is defined as
\[
\mathcal{F}^{-1}\lbrace \mathbf{F} \rbrace (x) = \int_{\mathbb{R}^n} \mathbf{F}(u)  e^{(2 \pi i_3 \langle x,u \rangle)} |du|
\]
where $x,u \in \mathbb{R}^n$, provided the integrals exist.
A multivector function $\mathbf{F}$ defined as 
\begin{align*}
\mathbf{F}(x) &= \mathbf{F}_0(x) +\mathbf{F}_1(x)e_1 +\mathbf{F}_2(x)e_2+\mathbf{F}_3(x)e_3\\ 
&+\mathbf{F}_{12}(x)e_{12}+\mathbf{F}_{23}(x)e_{23}+\mathbf{F}_{31}(x)e_{31}+\mathbf{F}_{123}(x)e_{123}
\end{align*}
can be written as 
\begin{align*}
\mathbf{F}(x) &= \left[\mathbf{F}_0(x) +\mathbf{F}_{123}(x)i_3\right] \\
			&+ \left[\mathbf{F}_1(x) +\mathbf{F}_{23}(x)i_3\right]e_1\\
			&+ \left[\mathbf{F}_2(x) +\mathbf{F}_{31}(x)i_3\right]e_2\\
			&+ \left[\mathbf{F}_3(x) +\mathbf{F}_{12}(x)i_3\right]e_3, 
\end{align*}
which can be seen as four complex signals and interpreted as an element of $\mathbb{C}^4$.
The linearity of the Clifford-Fourier transform would result in the following
\begin{align*}
\mathcal{F}\lbrace \mathbf{F} \rbrace (u) 
&= \left[\mathcal{F}\lbrace \mathbf{F}_0 +\mathbf{F}_{123}i_3 \rbrace (u)\right]1 \\
&+ \left[\mathcal{F}\lbrace \mathbf{F}_1 +\mathbf{F}_{23}i_3 \rbrace (u)\right]e_1,\\
&+ \left[\mathcal{F}\lbrace \mathbf{F}_2 +\mathbf{F}_{31}i_3 \rbrace (u)\right]e_2,\\
&+ \left[\mathcal{F}\lbrace \mathbf{F}_3 +\mathbf{F}_{12}i_3 \rbrace (u)\right]e_3,
\end{align*}
and this makes it plausible to deal with Clifford-Fourier transform in 3D as a linear combination of four classical Fourier transform. The derivative property in 3D is similar to the one in the classical Fourier transform because $i_3$ is commutative with any multivector in $\mathbb{R}^3$. Therefore \cite{ebling2005clifford},
\begin{align*}
\mathcal{F}\lbrace \nabla\mathbf{F} \rbrace (u) &= 2\pi i_3 u \mathcal{F}\lbrace \mathbf{F} \rbrace (u) \\
\mathcal{F}\lbrace \Delta\mathbf{F} \rbrace (u) &= -4\pi^2 u^2 \mathcal{F}\lbrace \mathbf{F} \rbrace (u).
\end{align*}

 In our applications, we use the Clifford-Fourier transform in 3D, denoted as CFT3, since it acts like the classical Fourier transform in terms of derivative property. 

\subsection{PDE transform}

In this section, we present a brief review of the partial differential equation transform that was proposed in our earlier work \cite{wang2011partial}. This transform is used to generate the biomolecular surfaces by applying it to a specific initial data-driven by the coordinates of the atoms in the molecule and their van der Waals radii. The next section offers a detailed explanation of the methods used in getting the initial data as well as the surface construction procedure.

Motivated by many physical phenomena in biological systems and pattern formation in nature, a family of high order PDEs for image processing was introduced in 1999 
\begin{equation}\label{eq:wei1999}
\begin{split}
\dfrac{\partial u(\mathbf{r},t)}{\partial t} =& 
\mathlarger{\sum}_q \nabla \cdot \left[d_q (u(\mathbf{r},t),| \nabla u(\mathbf{r},t)|,t)\nabla\nabla^{2q}u(\mathbf{r},t)\right]\\
+&e(u(\mathbf{r},t),| \nabla u(\mathbf{r},t)|,t), \qquad\qquad q=0,1,2,\cdots
\end{split}
\end{equation}
where $u(\mathbf{r},t)$ is the image function, $\nabla = \frac{\partial u(\mathbf{r},t)}{\partial r}$, $d_q (u(\mathbf{r},t),| \nabla u(\mathbf{r},t)|,t)$ is the edge sensitive diffusion coefficient and $e(u(\mathbf{r},t),| \nabla u(\mathbf{r},t)|,t)$ is the enhancement operator. Equation \ref{eq:wei1999} is a generalization of the Perona-Malik equation \cite{perona1990scale} that can be recovered if the enhancement operator is set to zero and $q=0$. The diffusion coefficients $d_q (u(\mathbf{r},t),| \nabla u(\mathbf{r},t)|,t)$ were defined as 
\begin{equation}
d_q (u(\mathbf{r},t),| \nabla u(\mathbf{r},t)|,t) = d_{q0}\exp\left[-\dfrac{| \nabla u |^2}{2\sigma^2_q}\right],
\end{equation}
where the values of $d_{q0}$ depend on the noise level, and $\sigma_q$ for $q=0,1$ were defined in terms of local statistical variance of $u$ and $\nabla u$ as
\begin{equation}
\sigma^2_q(\mathbf{r}) = \overline{|\nabla^qu - \overline{|\nabla^qu} |^2} \qquad\qquad\qquad  q=0,1.
\end{equation}
 The notation $\overline{Y(r)}$ represents the local average of $Y(r)$ centered at $\mathbf{r}$. The importance of the statistical measure based on the local statistical variance comes from its role in discriminating image features from noise. This advantage gives the ability to bypass the preprocessing done to noisy images where they get convolved with a test function or smooth mask \cite{wang2011partial}.

\begin{figure}[h!]
\center
\begin{tabular}{ c c c }
\subcaptionbox{}
{\includegraphics[scale=0.7,trim={2cm 0 2cm 0},clip]{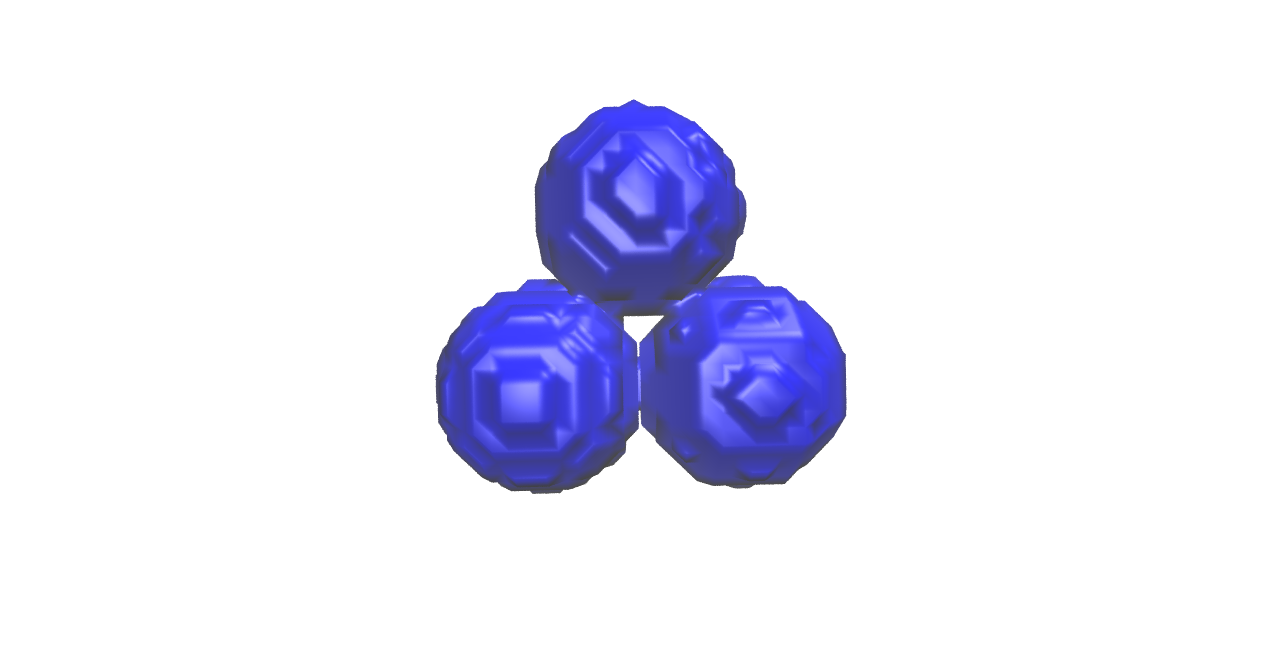}}
\subcaptionbox{}
{\includegraphics[scale=0.7,trim={2cm 0 2cm 0},clip]{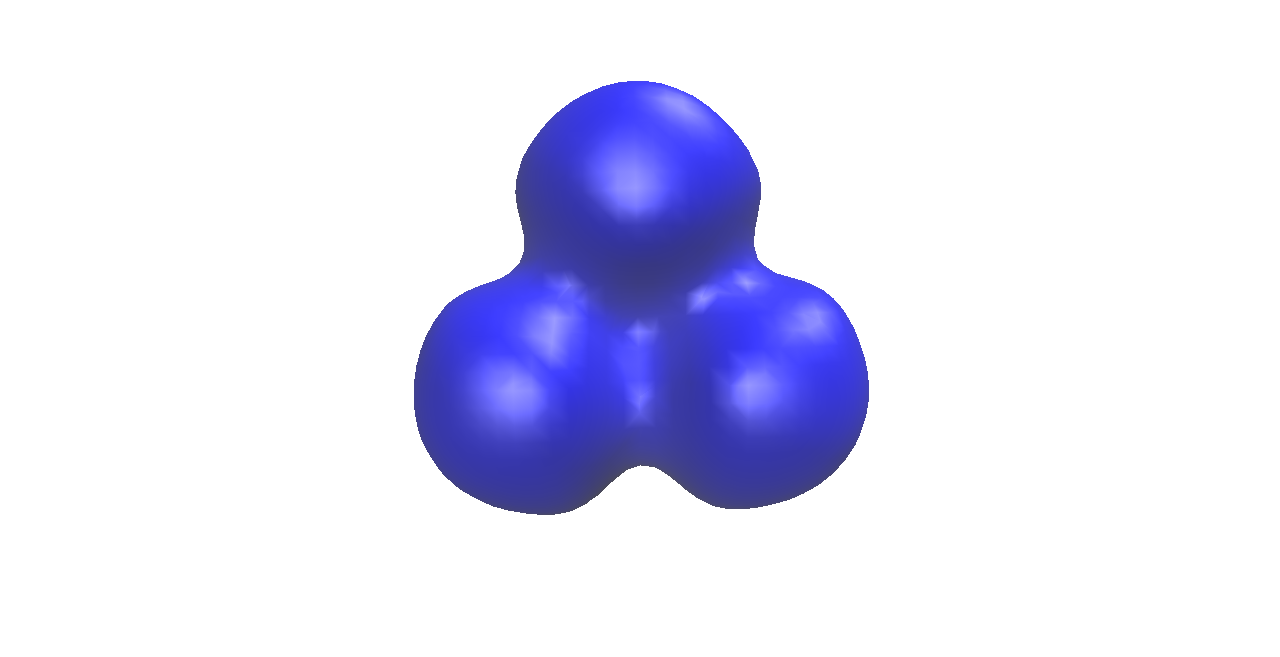}}
\subcaptionbox{}
{\includegraphics[scale=0.7,trim={2cm 0 2cm 0},clip]{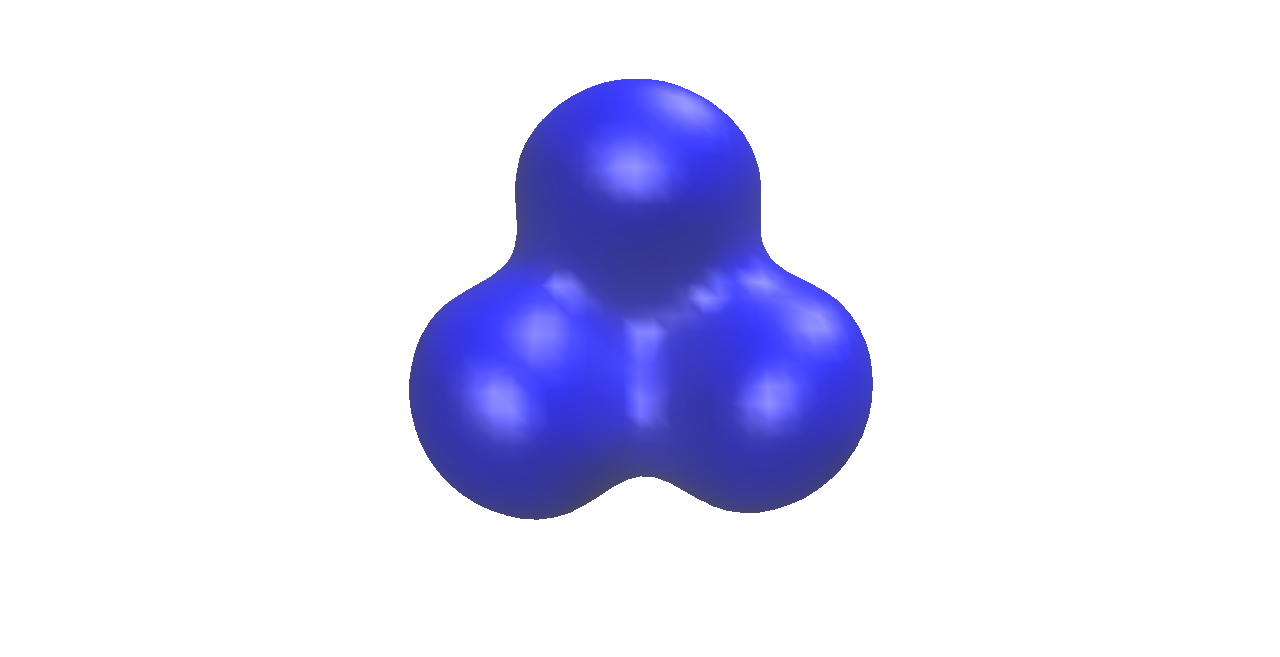}}\\
\subcaptionbox{}
{\includegraphics[scale=0.7,trim={2cm 0 2cm 0},clip]{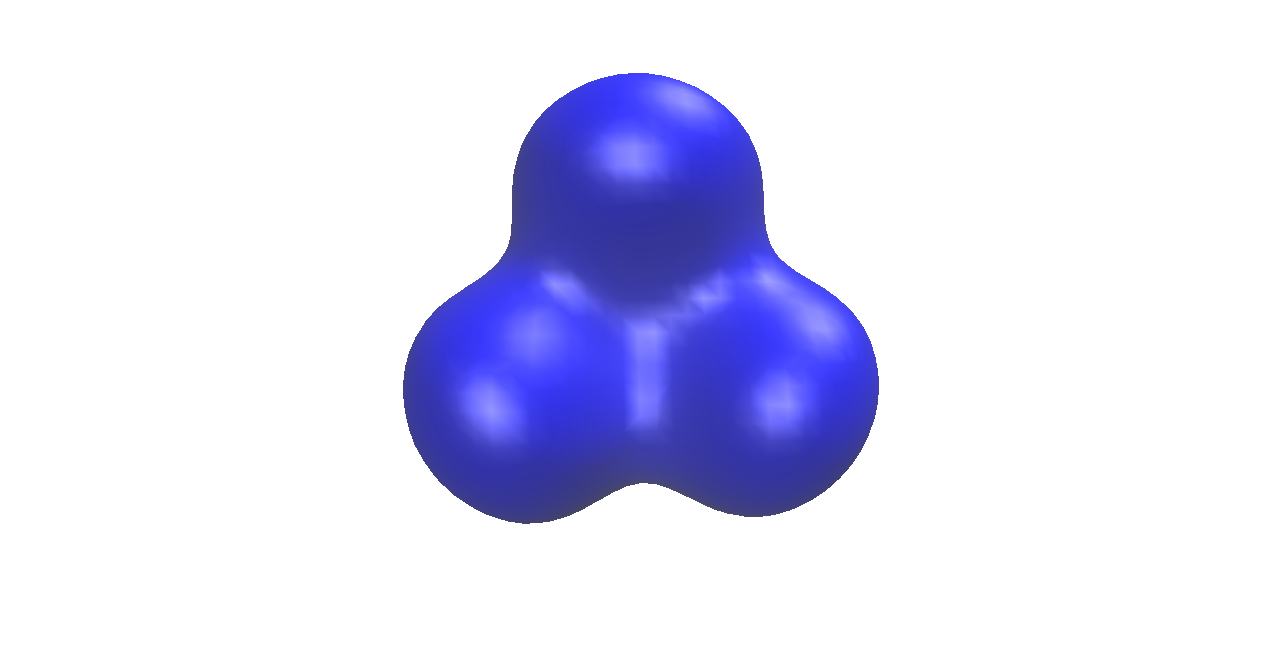}}
\subcaptionbox{}
{\includegraphics[scale=0.7,trim={2cm 0 2cm 0},clip]{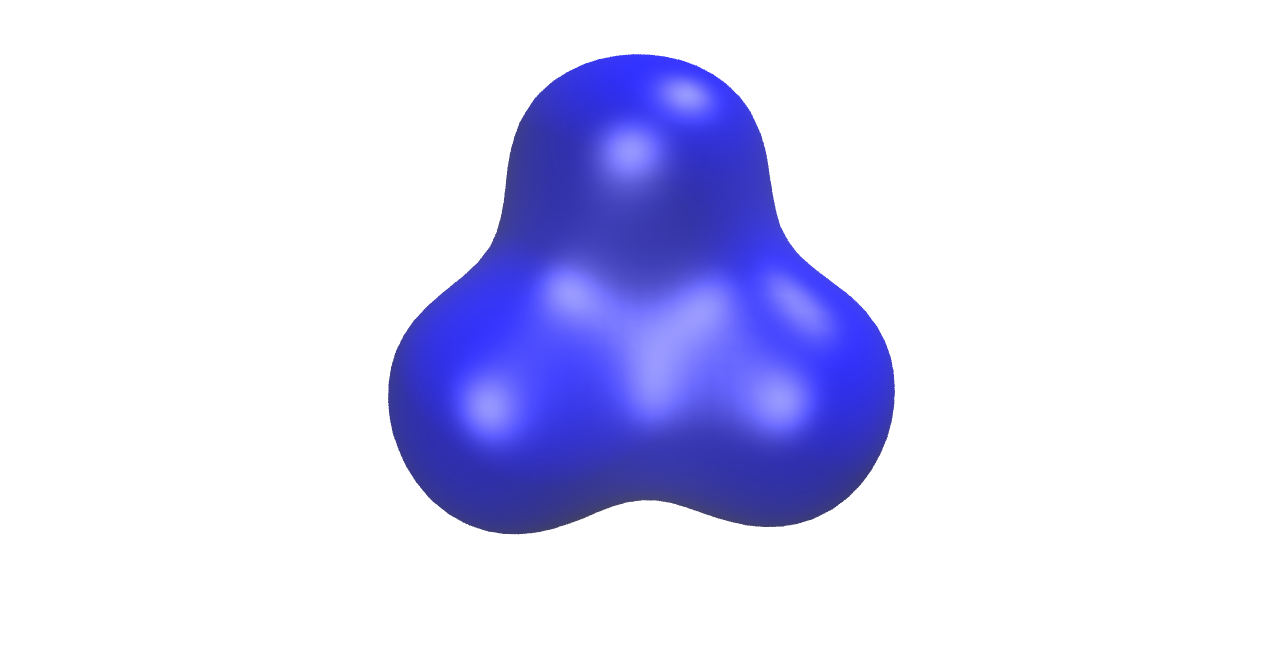}}
\subcaptionbox{}
{\includegraphics[scale=0.7,trim={2cm 0 2cm 0},clip]{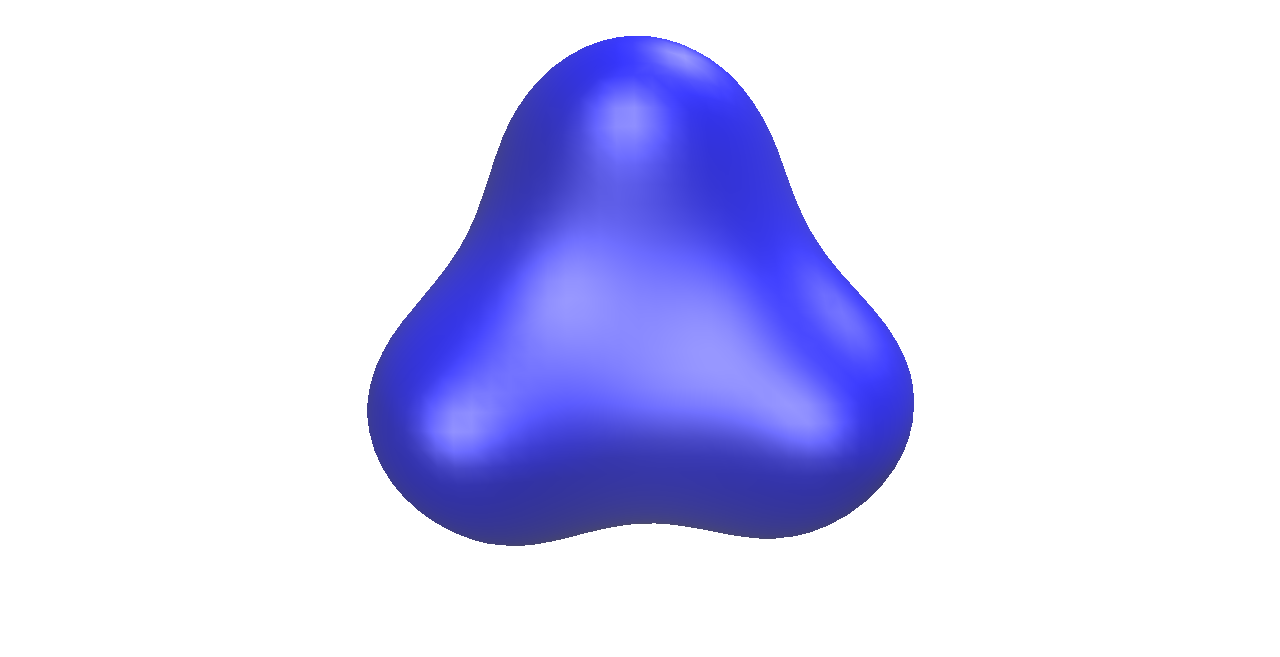}}
\end{tabular}
\caption{The isosurfaces of the three-atom-molecule generated with the Clifford-Fourier transform method using the piecewise initial data defined in Equation \ref{eq:pwinitial} and extracted at isovalue $=0.9$ and different propagation times.(a) The initial surface. (b) The isosurface with propagation time $t=10^1$. (c) $t=10^2$. (d) $t=10^3$. (e) $t=10^4$. (f) $t=10^5$.}
\label{fig:pwtimepropagation}
\end{figure}

The well-posedness of the generalized Perona-Malik equation proposed  was analyzed in terms of the existence and uniqueness of the solution \cite{bertozzi2004low,jin2010strong,xu2007existence}. The properties of Equation \ref{eq:wei1999} were shown to be different from the properties of other high order PDEs because it is not derived from a variational formulation \cite{jin2010strong}. The stability of Equation \ref{eq:wei1999} comes from appropriate choice of the coefficients $d_q (u(\mathbf{r},t),| \nabla u(\mathbf{r},t)|,t)$ \cite{wang2011partial}.

As noted in our earlier work \cite{wang2011partial}, the PDE transform can extract mode functions from some data given, say $X$, which is a very important property of the PDE transform. The solution of Equation \ref{eq:wei1999} can be found by the following equation
\begin{equation}\label{eq:modes}
\check{X}^k(\mathbf{r},t) = \mathcal{L} X^k(\mathbf{r})
\end{equation}
where $\mathcal{L}$ is a low-pass PDE tansform that satisfies
\begin{equation}
\mathcal{L}u(\mathbf{r},0) = u(\mathbf{r},t)
\end{equation}   
with $t$ being an artificial time involved in $\mathcal{L}$, $\check{X}^k(\mathbf{r},t)$ is the $k$th mode function and $X^k(\mathbf{r})$ is the $k$th residue function that is defined by
\begin{align}
X^1 &= X(\mathbf{r}) \\
X^k &= X^1 - \mathlarger{\sum}_{j=1}^{k-1}\check{X}^j, \qquad\qquad k=2,3,\cdots
\end{align} 
The original data $X$ can be reconstructed perfectly as \cite{zheng2012biomolecular}
\begin{equation}
X = \mathlarger{\sum}_{j=1}^{k-1}\check{X}^j + X^k.
\end{equation} 
Note that recursive applications of the PDE transform can generate the mode functions based on the input data, where the first mode is the trend of the data and the first residue is a general edge function. In contrast, high-pass PDE transform, proposed in our earlier work, were constructed in a way that the first mode is the edge type of information and the trend is the final residue \cite{wang2011partial}.

For the practical applications in this work, we   assume the following linearized form 
\begin{equation}\label{eq:linheat}
\partial_t u = \mathlarger{\sum}_{j=1}^{m} (-1)^{j+1} d_j \nabla^{2j}u +\epsilon (X^k - u), \qquad t\geq 0,
\end{equation}
 where $X^k \in \mathbb{R}^n$ is the $k$th residue of the data, $d_j>0$ and $\epsilon \sim 0$. This linearized equation is subject to the initial data $u(\mathbf{r},0)=X^k$. Solving this arbitrarily high-order PDE transform is computationally expensive. We use the fast Clifford-Fourier transform (FCFT) to make the computations more efficient computationally. The Clifford-Fourier transform is applied to both sides of Equation \ref{eq:modes} as follow
\begin{equation}
\hat{\check{X}}^k(\mathbf{r},t) = \hat{\mathcal{L}} \hat{X}^k(\mathbf{r})
\end{equation}
where $\hat{\mathcal{L}}$ is a frequency response function expressed as
\begin{equation}
\hat{\mathcal{L}}(\epsilon,t,d_1,d_2,\cdots,d_m) 
= e^{-\left(\mathlarger{\sum}_{j=1}^{m} d_j (w^2)^j + \epsilon\right)t}
+\dfrac{\epsilon}{\sum_{j=1}^{m} d_j (w^2)^j + \epsilon}\left(1-e^{-\left(\mathlarger{\sum}_{j=1}^{m} d_j (w^2)^j + \epsilon\right)t}\right)
\end{equation}
with $w^2 = \sum_{i=1}^{n}w_i^2$, and $\hat{\check{X}}^k$ and $\hat{X}$ are the Cilfford-Fourier transform of $\check{X}$ and $X$.

In the present work, periodic boundary condition is used whenever needed.

\section{Biomolecular surface generation}\label{sec:generation}
\begin{figure}[h!]
\center
\begin{tabular}{ c c c }
\subcaptionbox{}
{\includegraphics[scale=0.7,trim={2cm 0 2cm 0},clip]{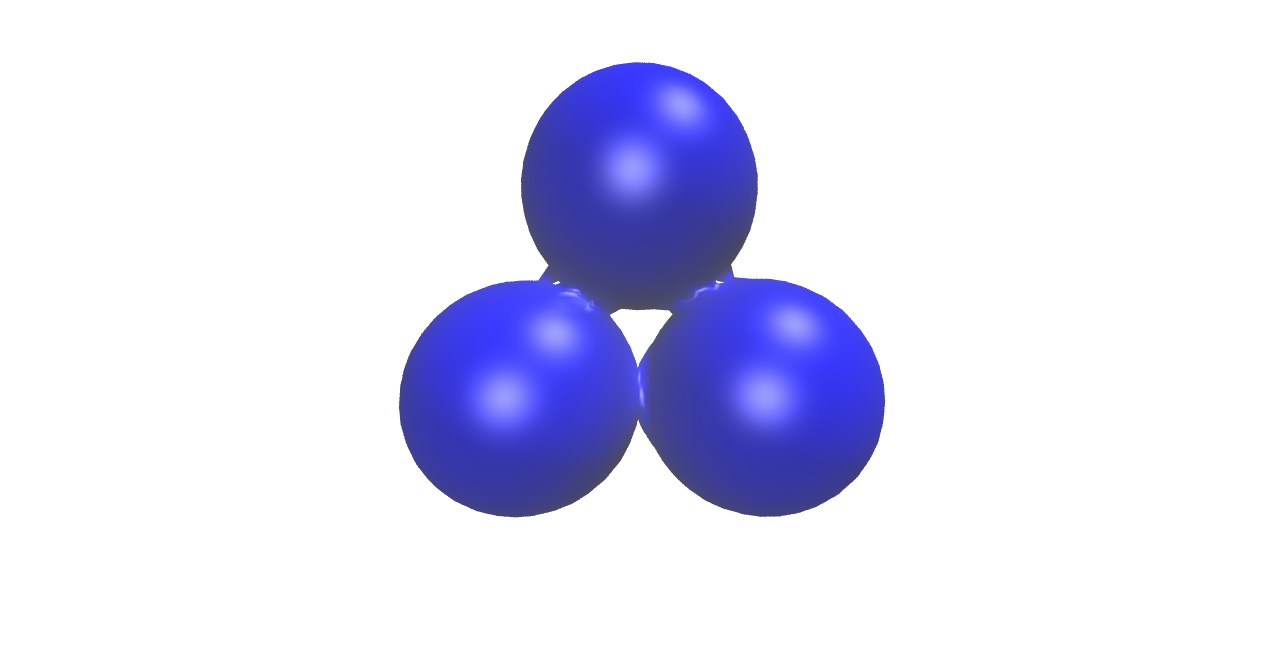}}
\subcaptionbox{}
{\includegraphics[scale=0.7,trim={2cm 0 2cm 0},clip]{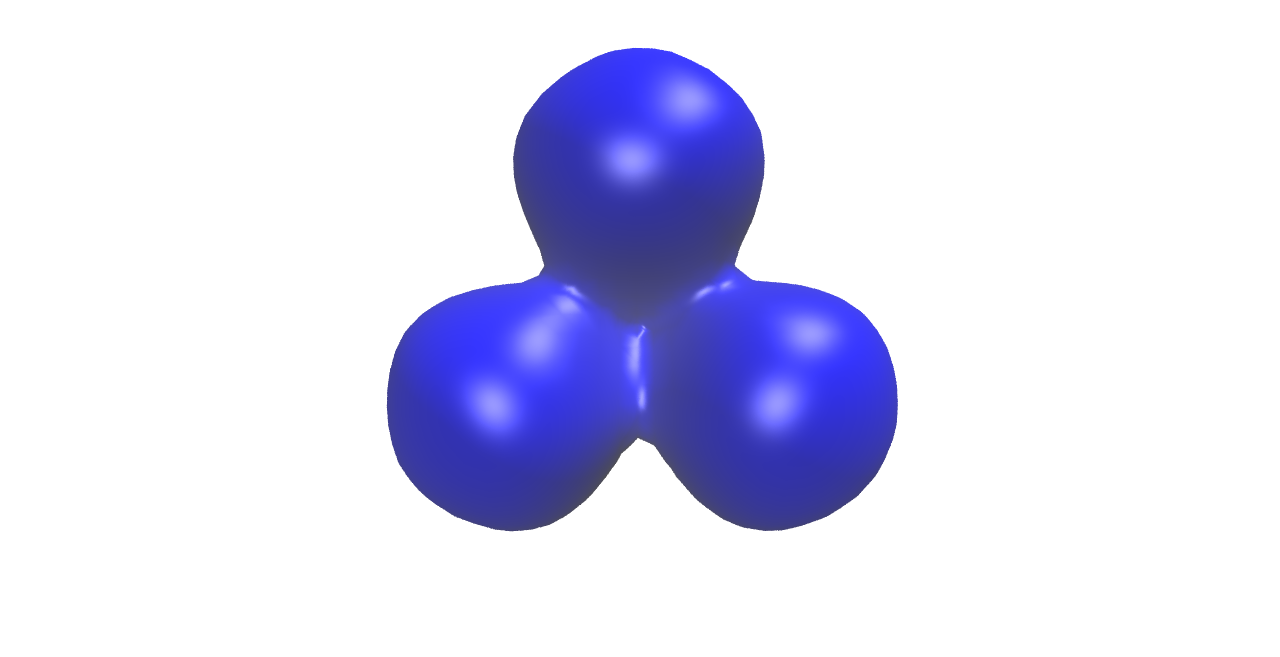}}
\subcaptionbox{}
{\includegraphics[scale=0.7,trim={2cm 0 2cm 0},clip]{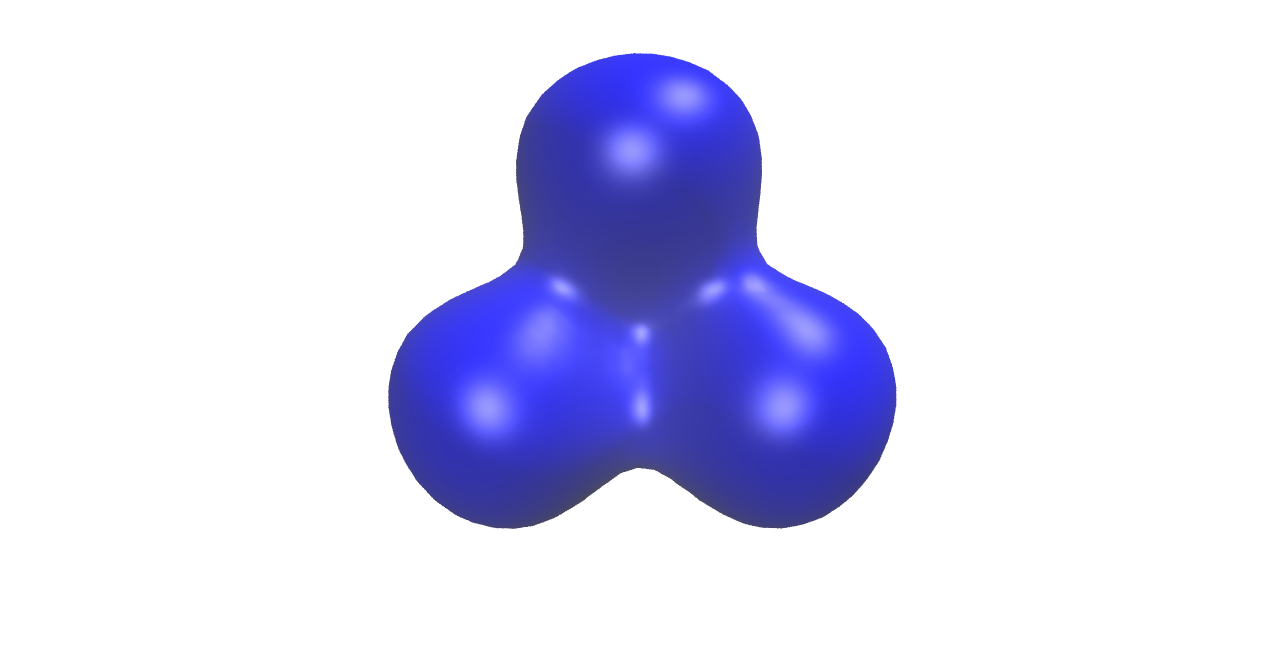}}\\
\subcaptionbox{}
{\includegraphics[scale=0.7,trim={2cm 0 2cm 0},clip]{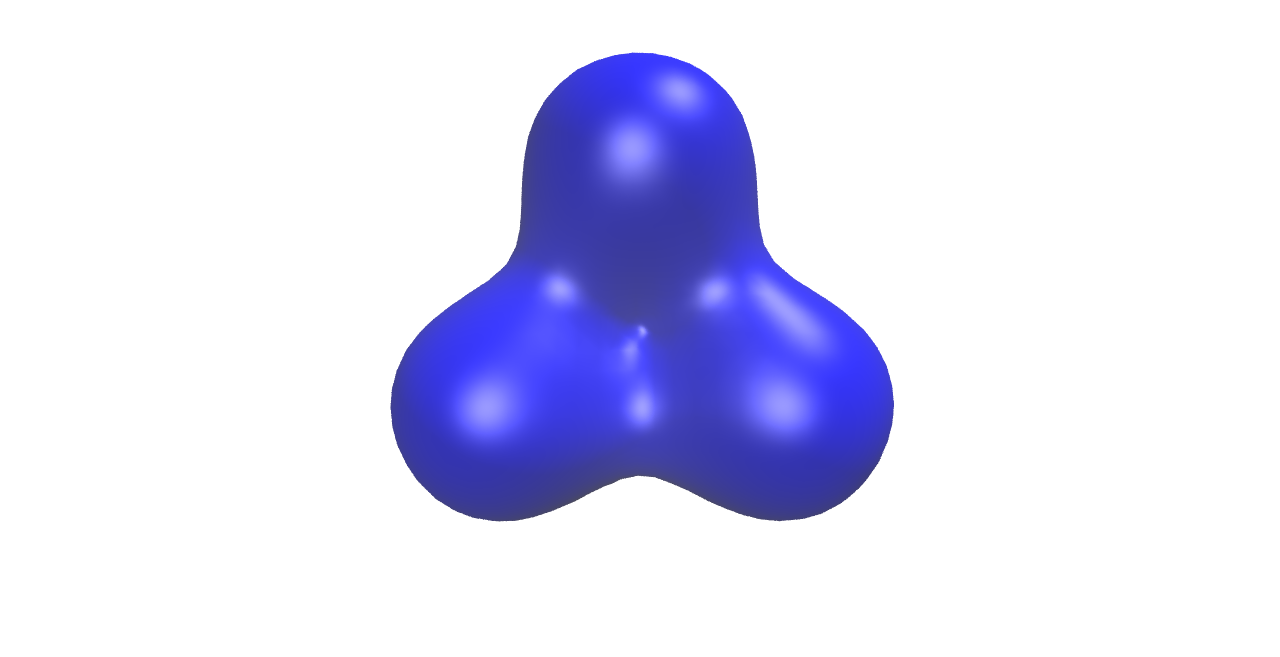}}
\subcaptionbox{}
{\includegraphics[scale=0.7,trim={2cm 0 2cm 0},clip]{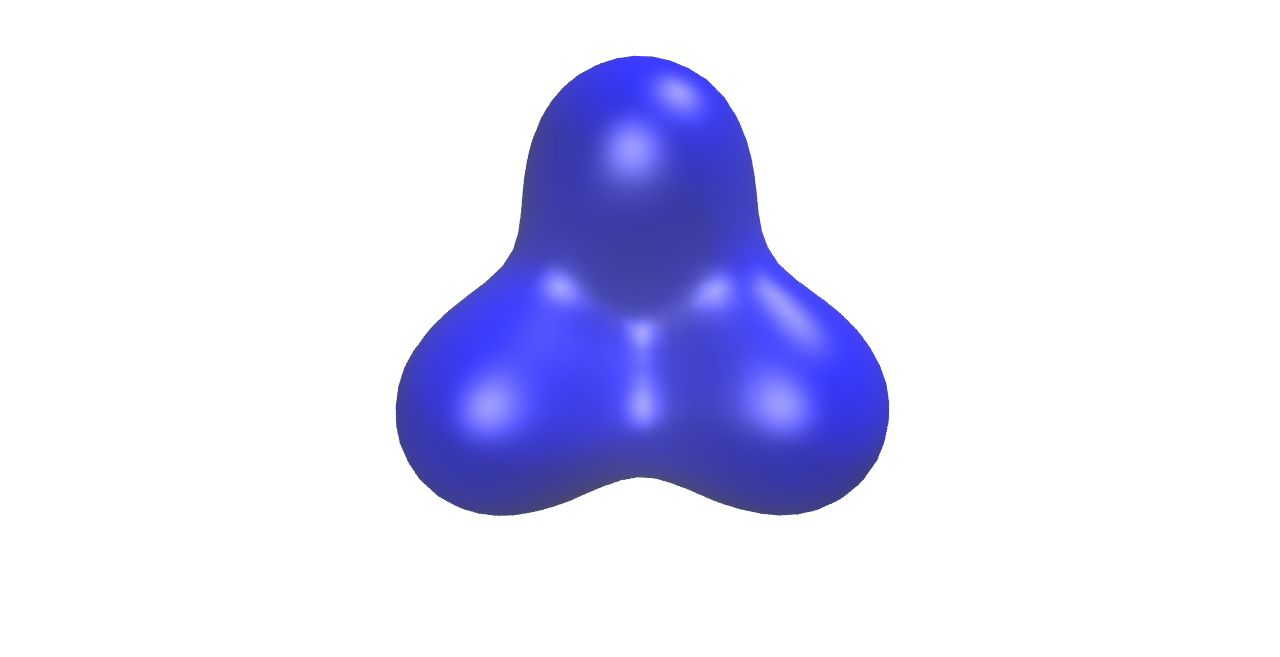}}
\subcaptionbox{}
{\includegraphics[scale=0.7,trim={2cm 0 2cm 0},clip]{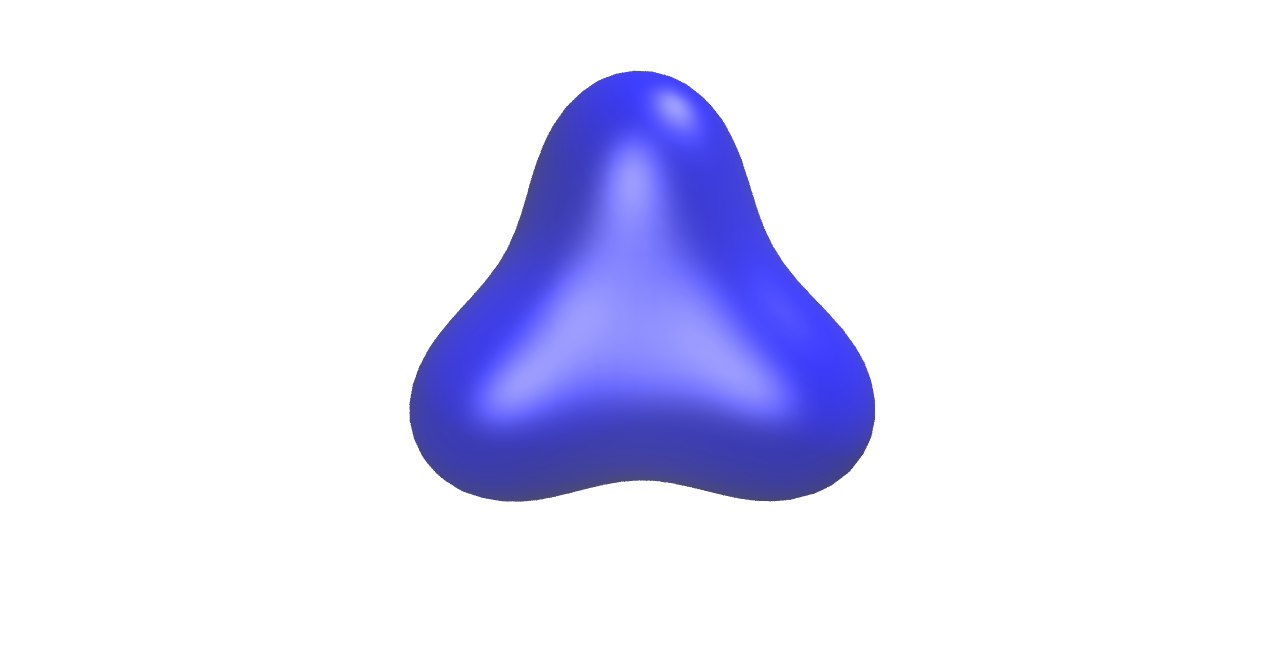}}
\end{tabular}
\caption{The isosurfaces of the imaginary three-atom-molecule generated with the Clifford-Fourier transform method using the Gaussian initial data defined in Equation \ref{eq:gaussianinitial} and extracted at isovalue $=1.0$ and different propagation times.(a) The initial surface. (b) The isosurface with propagation time $t=10^1$. (c) $t=10^2$. (d) $t=10^3$. (e) $t=10^4$. (f) $t=10^5$.}
\label{fig:gaussiantimepropagation}
\end{figure}

In this part, we propose a Clifford-Fourier transform-based surface generation method. First, we give a brief explanation of the method. Then, we show examples of surfaces generated using this method. In the next section, surface electrostatic potentials of some proteins are shown. The surfaces generated are compared with MSMS surfaces in terms of geometric singularities and electrostatic solvation free energy.  

The first step in a surface generation is to make an initial shape $u(\mathbf{r},0)$ driven by the coordinates and atomic radii of the atoms in the protein of interest. Then, we apply a PDE transform with specific parameters of time and order. After that, the Clifford-Fourier transform is applied to attain the final shape. From the final shape, we generate the molecular surface by extracting a specific isosurface. For the initial data, we have two cases for the initial data used in this method: piece-wise initial data and Gaussian initial data. 
\subsection{Initial data}
The first one of the two initial data cases is to use the piecewise initial data that we used in our earlier work \cite{zheng2012biomolecular,wei2005molecular}. The piecewise initial value $u(\mathbf{r},0)$ is defined as 

\begin{equation}\label{eq:pwinitial}
  u(\mathbf{r},0) =
  \begin{cases}
  0, & \mathbf{r} \in \bigcup\limits_{\beta = 1,\cdots,N_{\beta}} O(\mathbf{r}_{\beta},r_{\beta}),  \\
  1, & \text{otherwise,} 
  \end{cases}
\end{equation}
where $O(\mathbf{r}_{\beta},r_{\beta})$ is the sphere centered at $\mathbf{r}_{\beta}$ and has a radius $r_{\beta}$, i.e. $O(\mathbf{r}_{\beta},r_{\beta}): \lbrace \mathbf{r} \in \mathbb{R}^3 , \| \mathbf{r} - \mathbf{r}_{\beta}   \| \leq r_{\beta} \rbrace $ with $ \mathbf{r}_{\beta}$ and $r_{\beta}$ being the coordinates of a specific atom in the molecule and its atomic radius respectively and $N_{\beta}$ is the total number of atoms in that molecule. In our present work, we   take the van der Waals radius to be the atomic radius. So, this equation means that if $\mathbf{r}$ is in the sphere $O(\mathbf{r}_{\beta},r_{\beta})$ of any atom in the molecule, then $u(\mathbf{r},0)=0$, and if it is outside of any sphere, then $u(\mathbf{r},0)=1$. As noted, this initial shape represents the van der Waals surface which is non-smooth. Another non-smooth definition of a piecewise case can be achieved by switching the region with value 0 to 1 and vice versa. This latter case was used in our earlier work extensively \cite{bates2008minimal,zheng2012biomolecular}.
\begin{figure}[h!]
\center
\begin{tabular}{ c c c }
\subcaptionbox{}
{\includegraphics[scale=0.7,trim={2cm 0 2cm 0},clip]{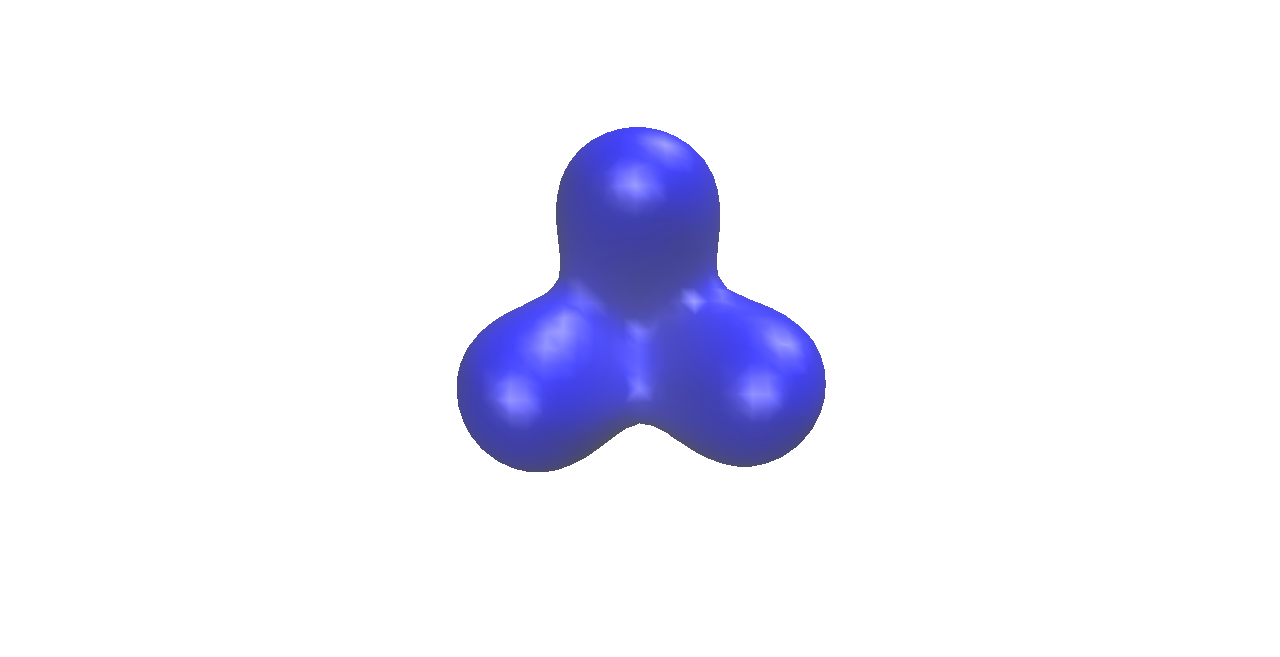}}
\subcaptionbox{}
{\includegraphics[scale=0.7,trim={2cm 0 2cm 0},clip]{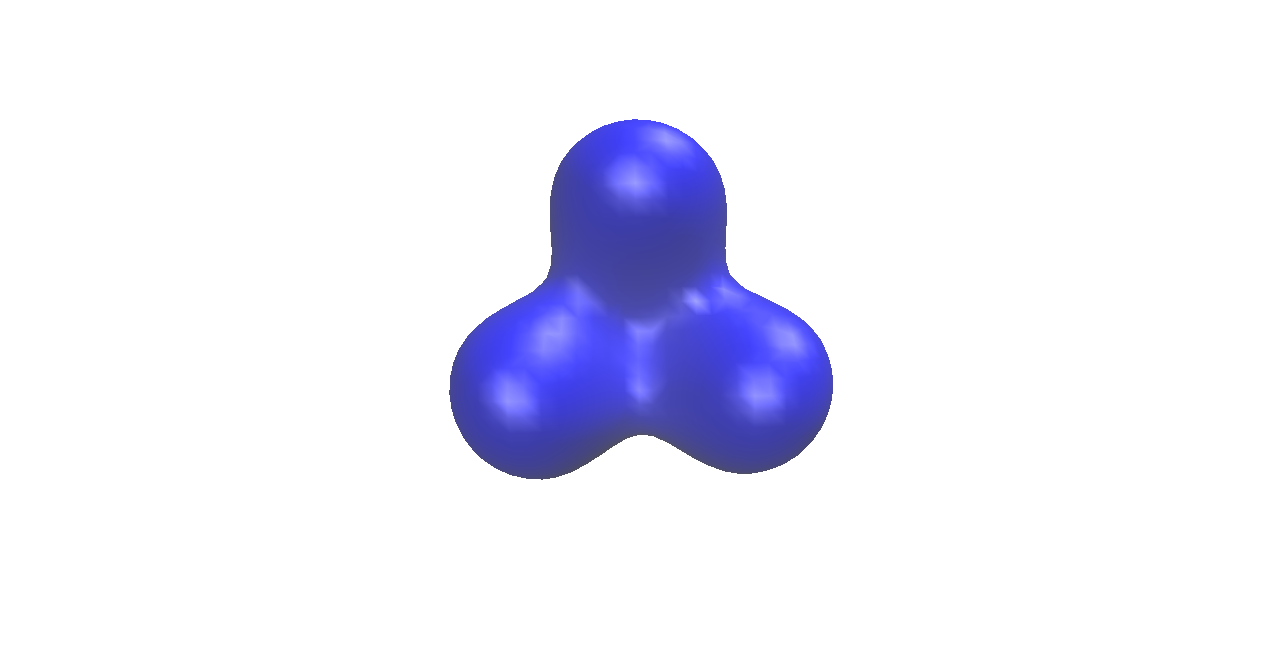}}
\subcaptionbox{}
{\includegraphics[scale=0.7,trim={2cm 0 2cm 0},clip]{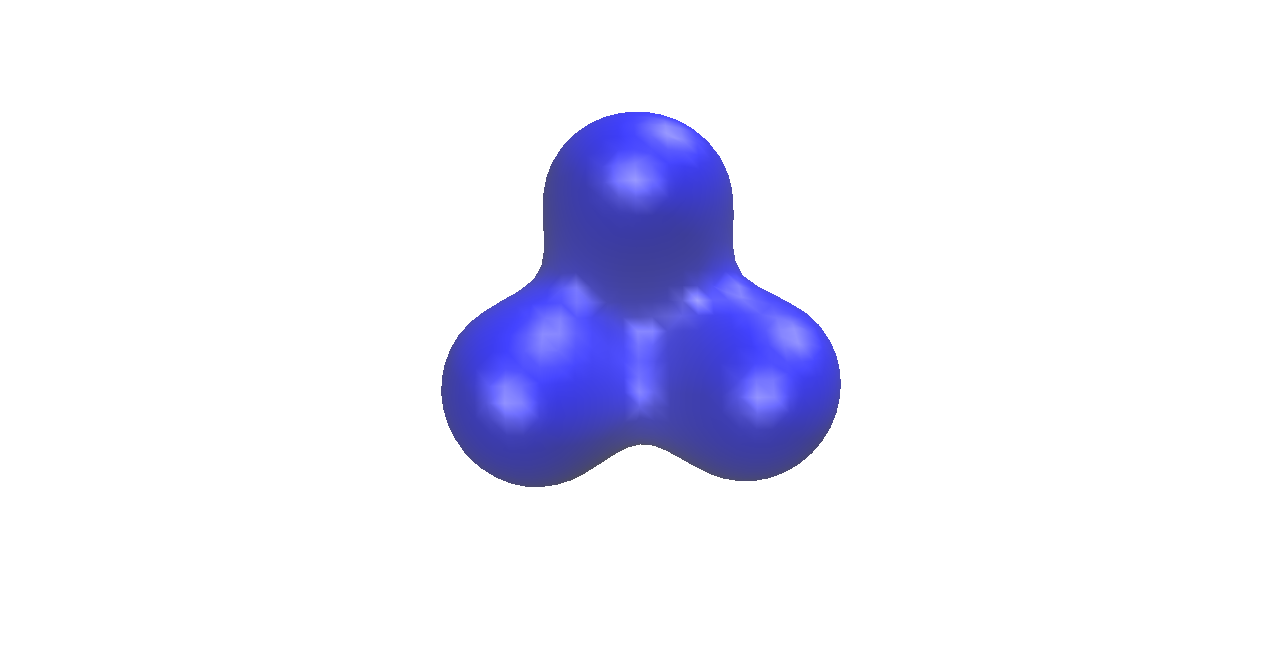}}\\
\subcaptionbox{}
{\includegraphics[scale=0.7,trim={2cm 0 2cm 0},clip]{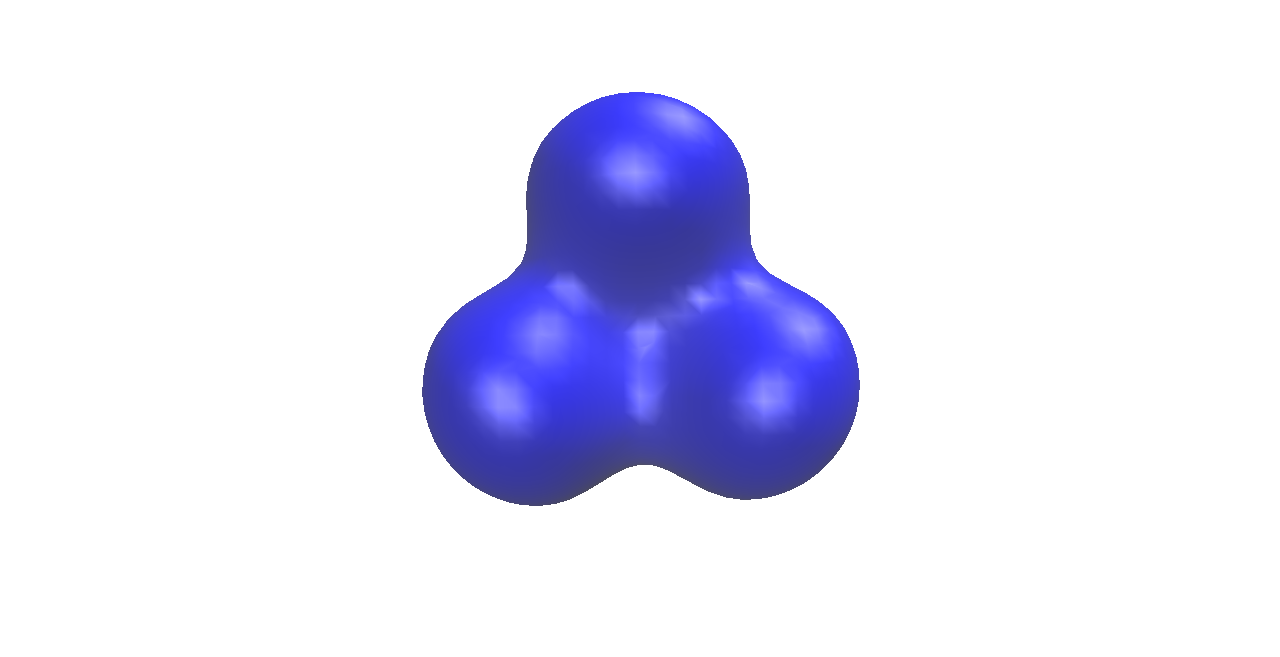}}
\subcaptionbox{}
{\includegraphics[scale=0.7,trim={2cm 0 2cm 0},clip]{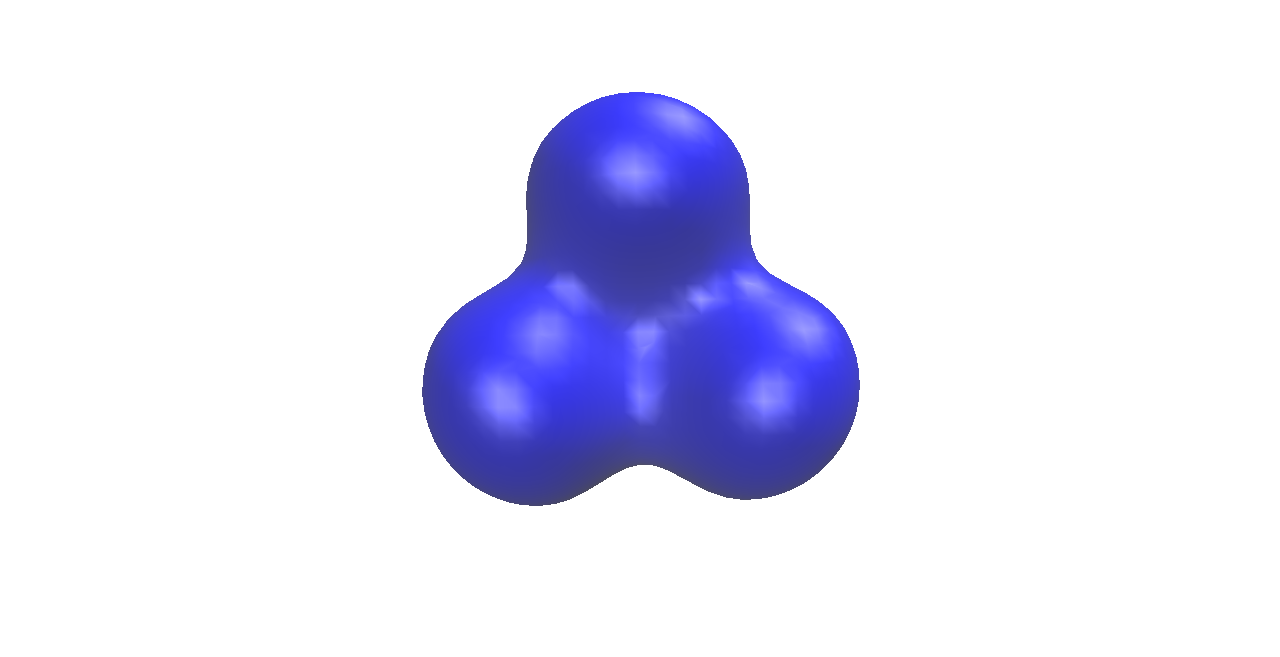}}
\subcaptionbox{}
{\includegraphics[scale=0.7,trim={2cm 0 2cm 0},clip]{pw_09_100}}
\end{tabular}
\caption{The isosurfaces of the imaginary three-atom-molecule generated with the Clifford-Fourier transform method using the piecewise initial data defined in Equation \ref{eq:pwinitial} and extracted at propagation time $t=10^2$ and different isovalues.(a) The isosurface extracted at isovalue $=0.4$. (b) isovalue $=0.5$. (c) isovalue $=0.6$.  (d) isovalue $=0.7$. (e) isovalue $=0.8$. (f) isovalue $=0.9$. }
\label{fig:pw100isovalue}
\end{figure}

The second case of initial data is achieved using Gaussian functions. Gaussian functions have been used to generate molecular surfaces in the literature \cite{yu2008feature,giard2010molecular,zhang2006origin} and we exploited them in our earlier work of surface generation using PDE transform \cite{zheng2012biomolecular}. In this work, we adopt the Gaussian function proposed in our earlier work \cite{zheng2012biomolecular} which is a modified version of the one Giard and Macq \cite{giard2010molecular} defined as
\begin{equation}\label{eq:gaussianinitial}
u(\mathbf{r},0) = \max_{\beta} \left(se^{-\dfrac{\| \mathbf{r} - \mathbf{r}_{\beta}   \|^2 - r_{\beta}^2}{r_{e}^2}}\right),
\end{equation}

where $s$ is the threshold parameter and $r_e$ is set to $3$ \AA. In this case of the initial value, the surface is not represented directly by the Gaussian function but the surface is indeed embedded within the Gaussian function. As noted, this case represents a smooth function but this does not give it any superiority over the non-smooth case in surface generation as discussed later.

\begin{figure}[h!]
\center
\begin{tabular}{ c c c }
\subcaptionbox{}
{\includegraphics[scale=0.7,trim={2cm 0 2cm 0},clip]{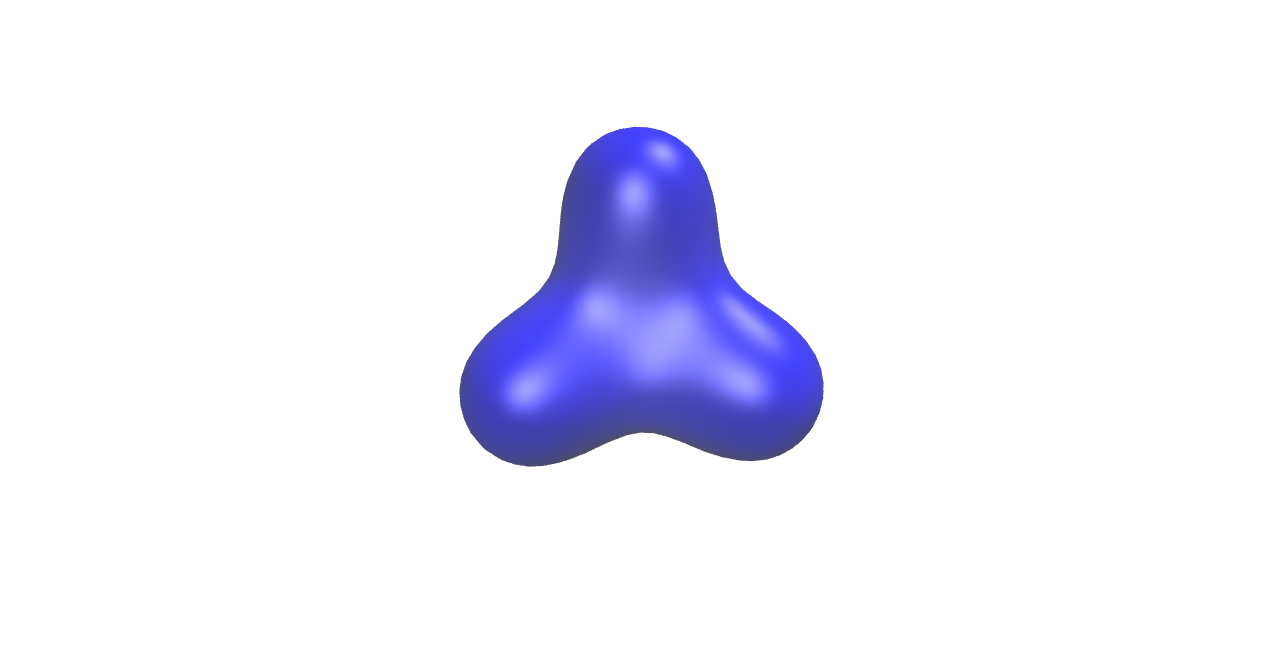}}
\subcaptionbox{}
{\includegraphics[scale=0.7,trim={2cm 0 2cm 0},clip]{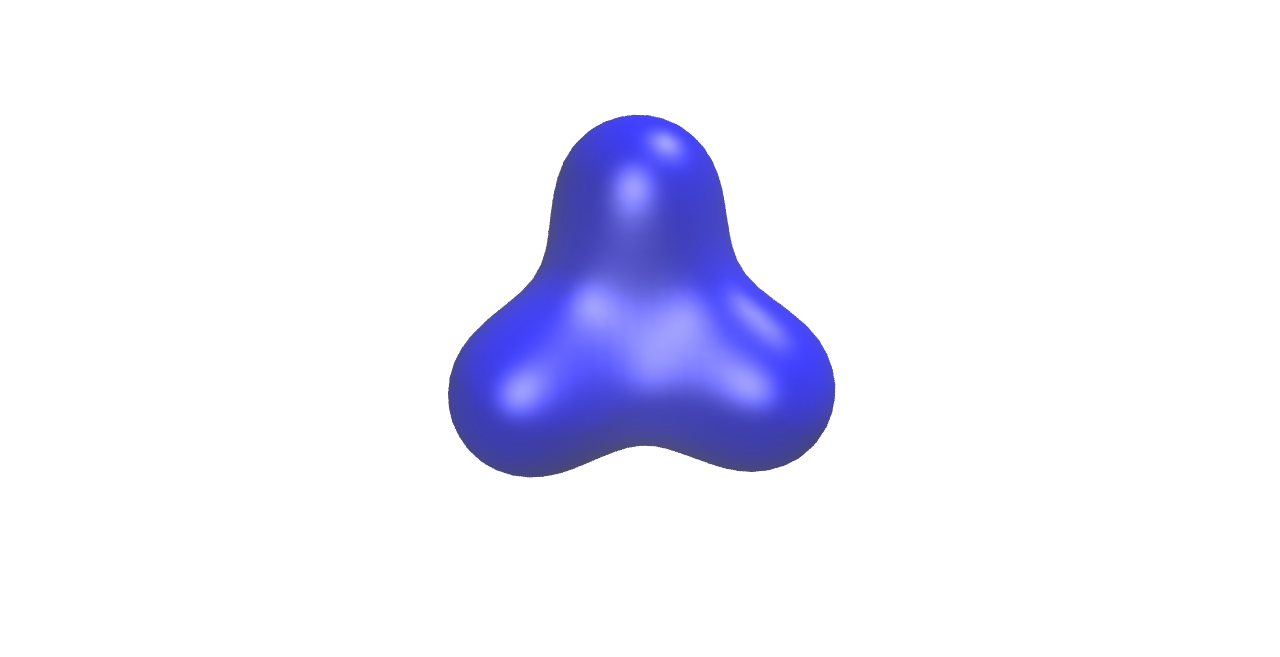}}
\subcaptionbox{}
{\includegraphics[scale=0.7,trim={2cm 0 2cm 0},clip]{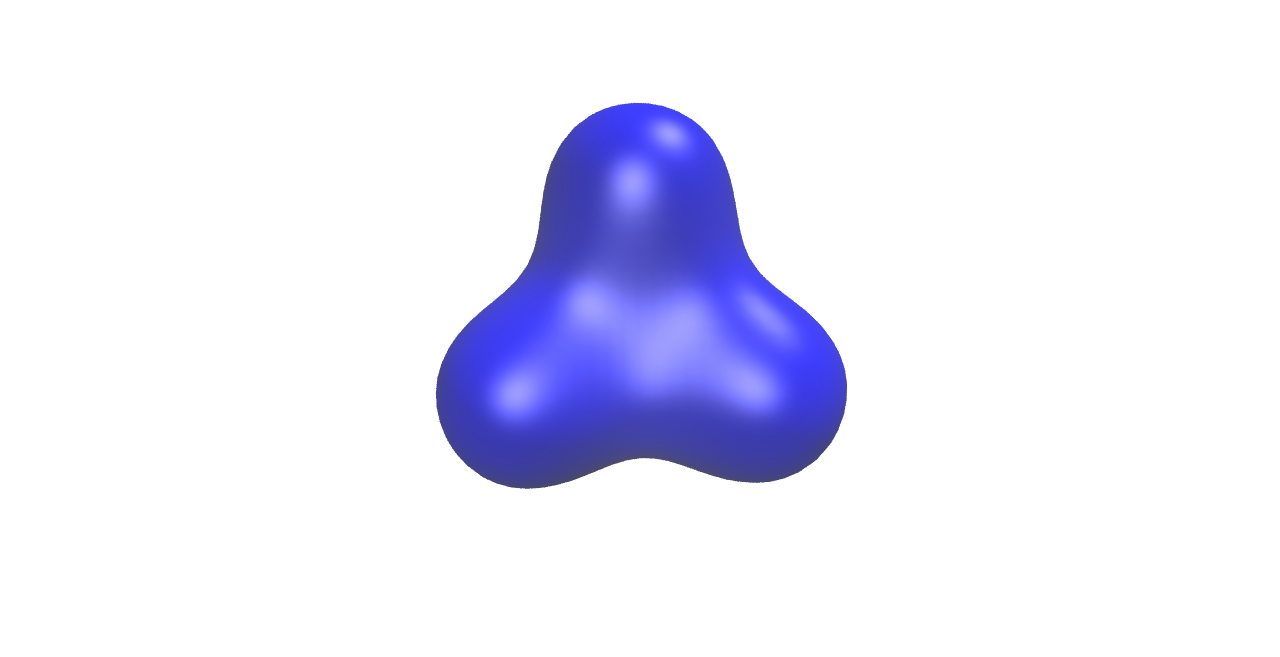}}\\
\subcaptionbox{}
{\includegraphics[scale=0.7,trim={2cm 0 2cm 0},clip]{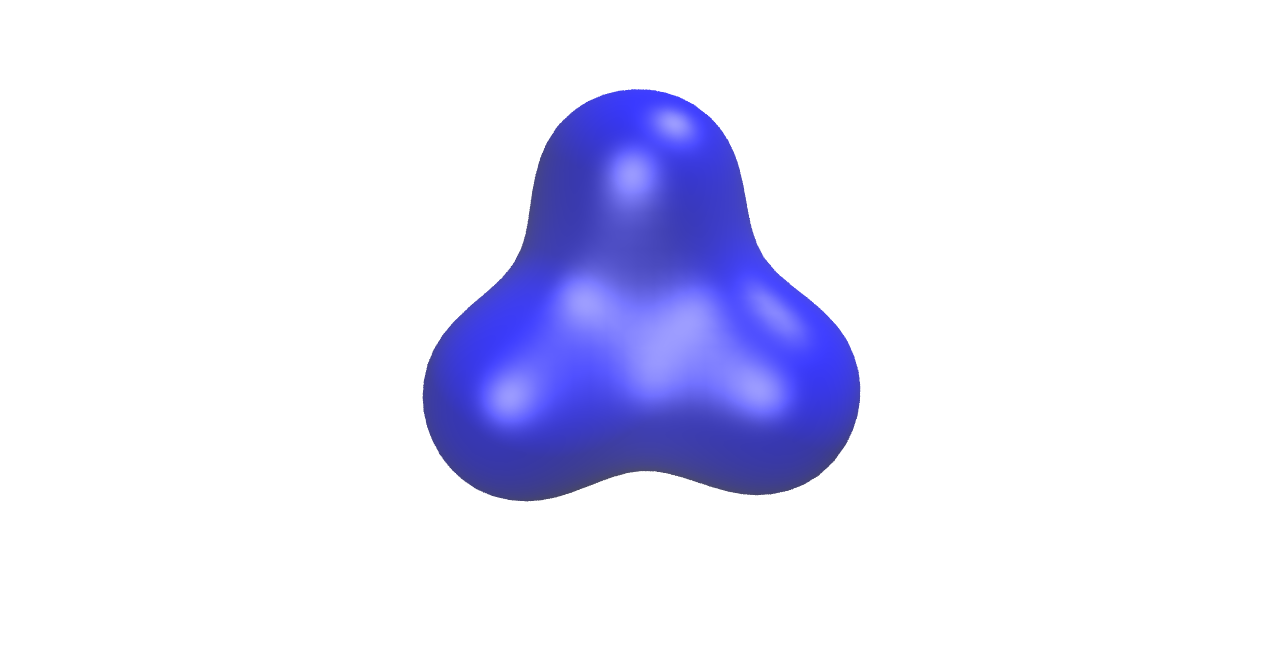}}
\subcaptionbox{}
{\includegraphics[scale=0.7,trim={2cm 0 2cm 0},clip]{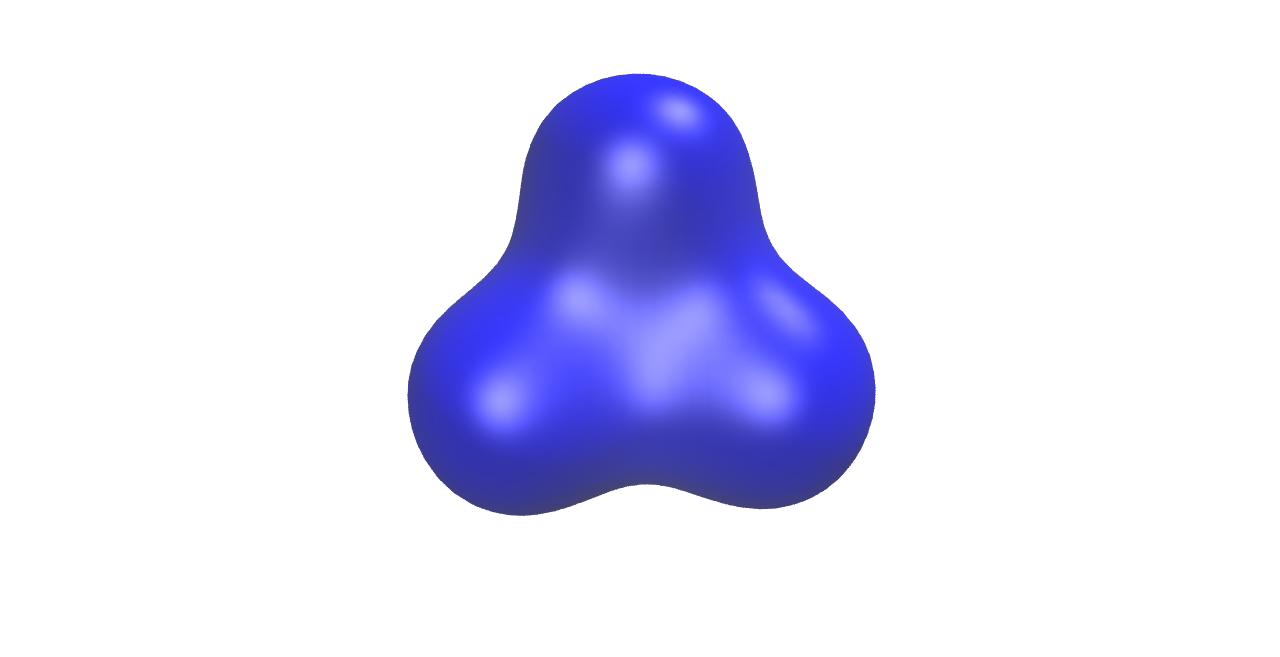}}
\subcaptionbox{}
{\includegraphics[scale=0.7,trim={2cm 0 2cm 0},clip]{pw_09_10000}}
\end{tabular}
\caption{The isosurfaces of the imaginary three-atom-molecule generated with the Clifford-Fourier transform method using the piecewise initial data defined in Equation \ref{eq:pwinitial} and extracted at propagation time $t=10^4$ and different isovalues.(a) The isosurface extracted at isovalue $=0.4$. (b) isovalue $=0.5$. (c) isovalue $=0.6$.  (d) isovalue $=0.7$. (e) isovalue $=0.8$. (f) isovalue $=0.9$. }
\label{fig:pw10000isovalue}
\end{figure}
\subsection{Test cases}
Now, we conduct different experiments on a test case of an imaginary molecule that is composed of three atoms. The coordinates of the centers of atoms are $(0, 0, 1.8), (0, 0,-1.8), \text{ and } (0, 3.12, 0)$ and each of them has an atomic radius of $1.8$ \AA. First, we explore the effect of the propagation time on the extracted isosurfaces. The experiments are done with the following propagation times $t=10^1,t=10^2,t=10^3,t=10^4,$ and $t=10^5$ to both initial values, piecewise and Gaussian. After that, the effect of changing the isovalue on the extracted isosurfaces is investigated. We carry out this first to the piecewise case with the following isovalues $0.4,0.5,0.6,0.7,0.8$, and $0.9$. This experiment is done twice to see the effect in two different times $t=10^2$ and $t=10^4$. Then, this investigation is carried out on the Gaussian initial shape with the following isovalues $0.5,0.6,0.7,0.8,0.9$, and $1.0$ and it is done twice to see the effect in two different times $t=10^2$ and $t=10^4$. 

\begin{figure}[h!]
\center
\begin{tabular}{ c c c }
\subcaptionbox{}
{\includegraphics[scale=0.7,trim={2cm 0 2cm 0},clip]{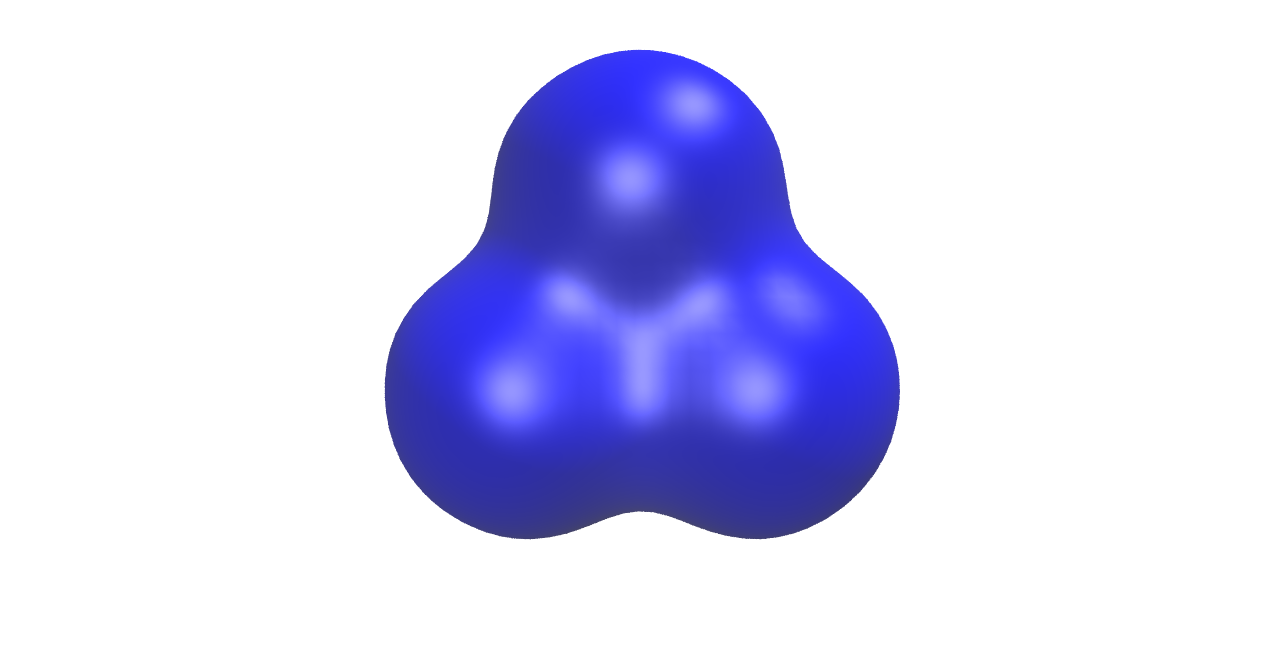}}
\subcaptionbox{}
{\includegraphics[scale=0.7,trim={2cm 0 2cm 0},clip]{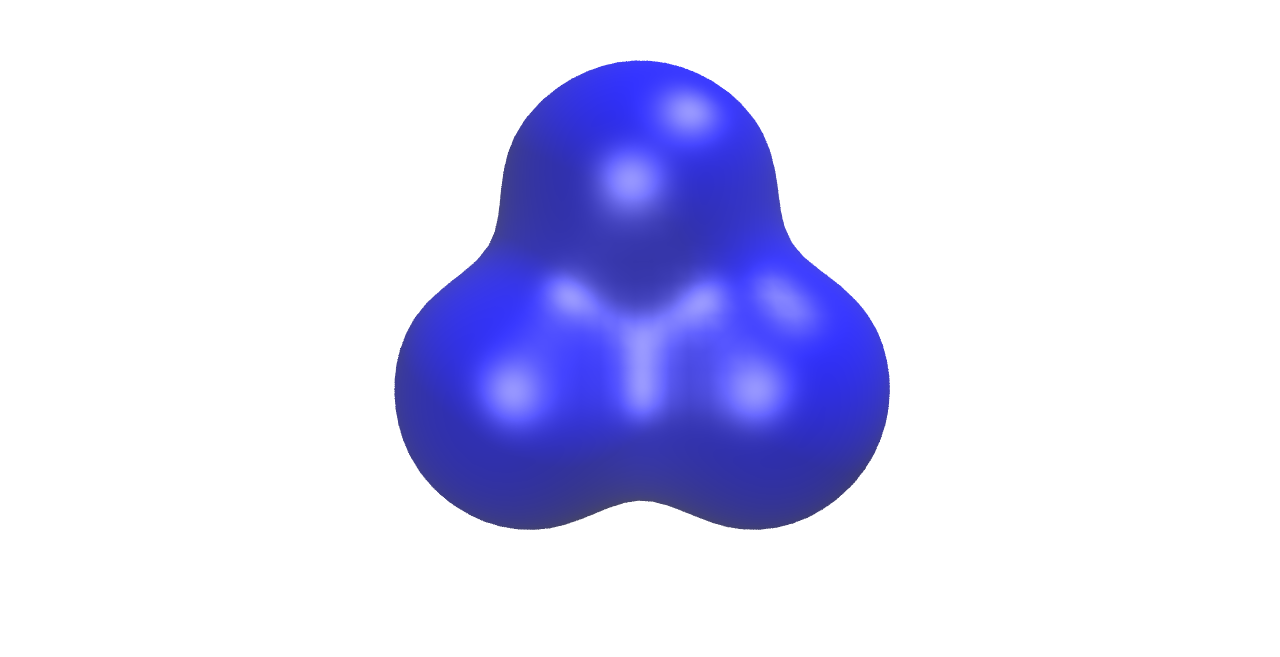}}
\subcaptionbox{}
{\includegraphics[scale=0.7,trim={2cm 0 2cm 0},clip]{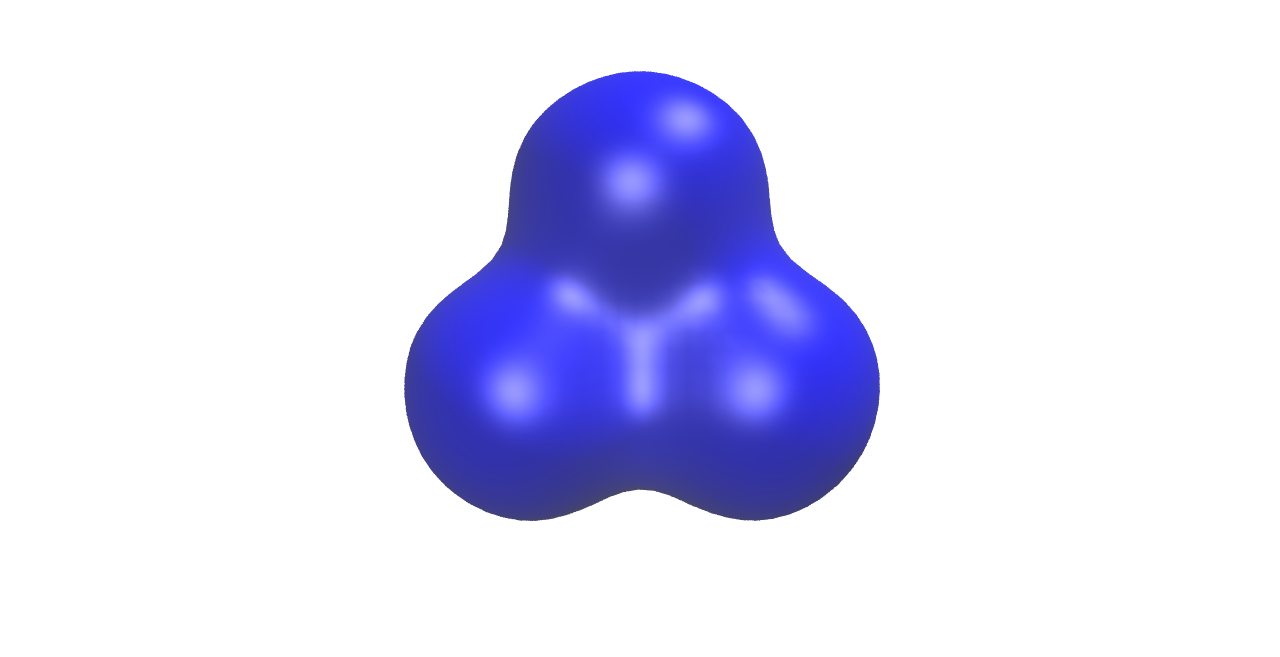}}\\
\subcaptionbox{}
{\includegraphics[scale=0.7,trim={2cm 0 2cm 0},clip]{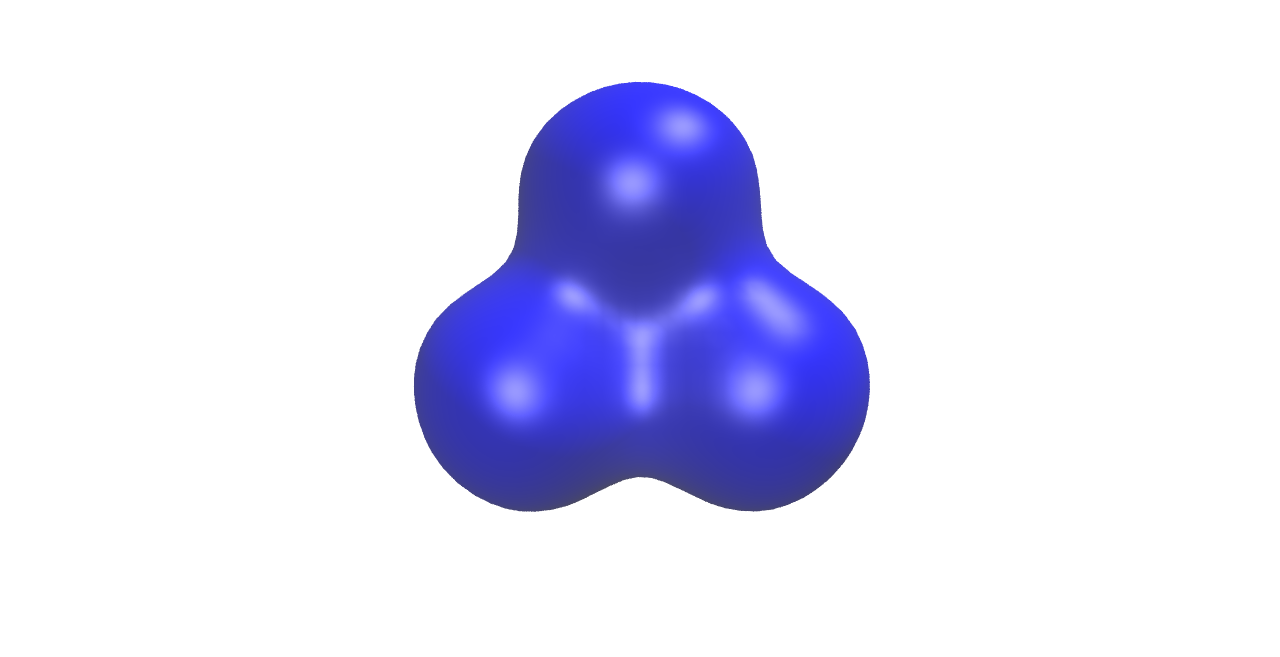}}
\subcaptionbox{}
{\includegraphics[scale=0.7,trim={2cm 0 2cm 0},clip]{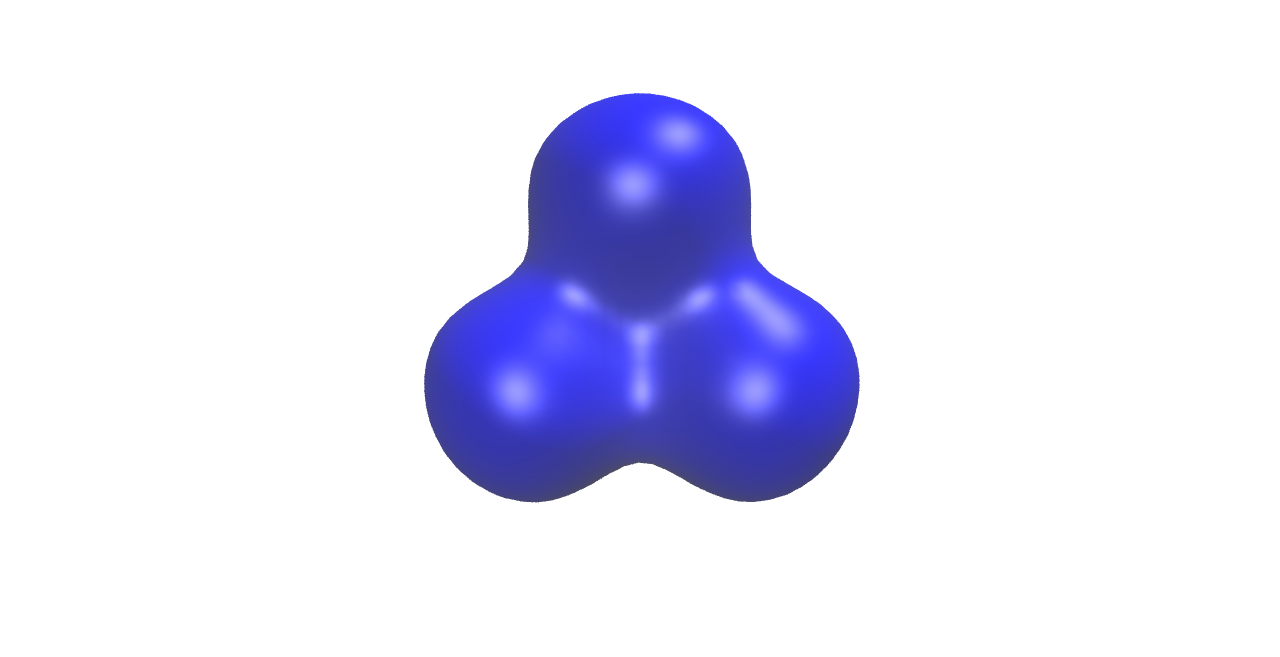}}
\subcaptionbox{}
{\includegraphics[scale=0.7,trim={2cm 0 2cm 0},clip]{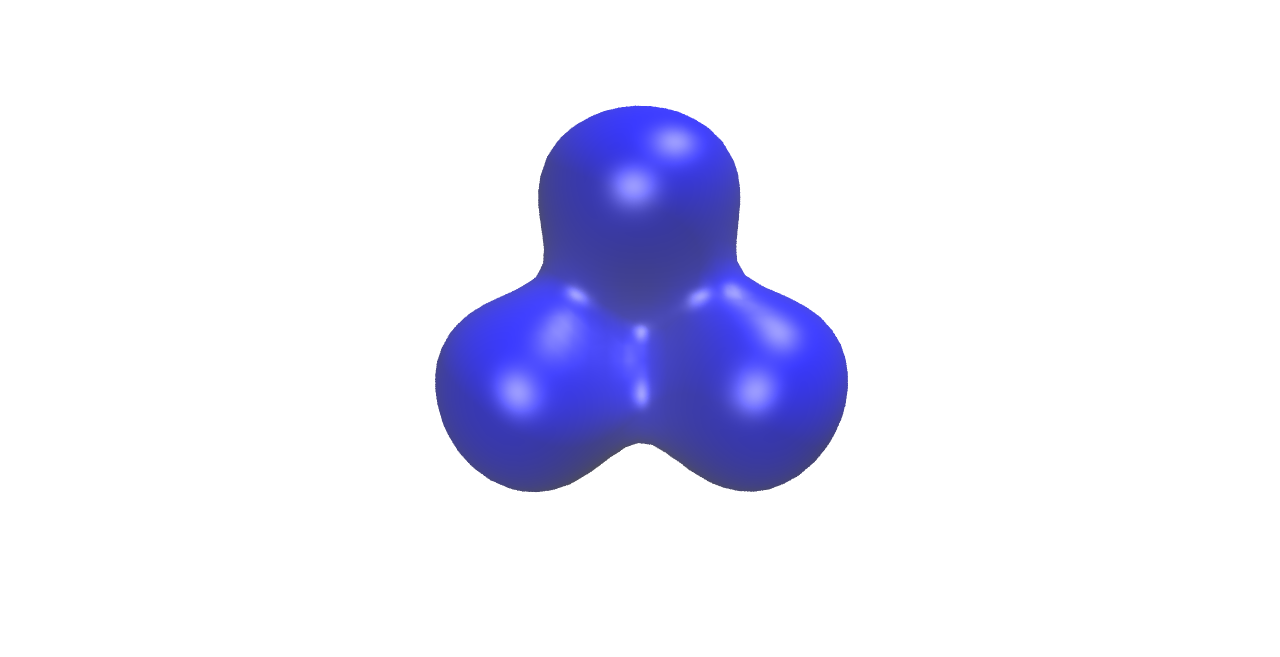}}
\end{tabular}
\caption{The isosurfaces of the imaginary three-atom-molecule generated with the Clifford-Fourier transform method using the Gaussian initial data defined in Equation \ref{eq:gaussianinitial} and extracted at propagation time $t=10^2$ and different isovalues.(a) The isosurface extracted at isovalue $=0.5$. (b) isovalue $=0.6$. (c) isovalue $=0.7$.  (d) isovalue $=0.8$. (e) isovalue $=0.9$. (f) isovalue $=1.0$. }
\label{fig:gaussian100isovalue}
\end{figure}

\begin{figure}[h!]
\center
\begin{tabular}{ c c c }
\subcaptionbox{}
{\includegraphics[scale=0.7,trim={2cm 0 2cm 0},clip]{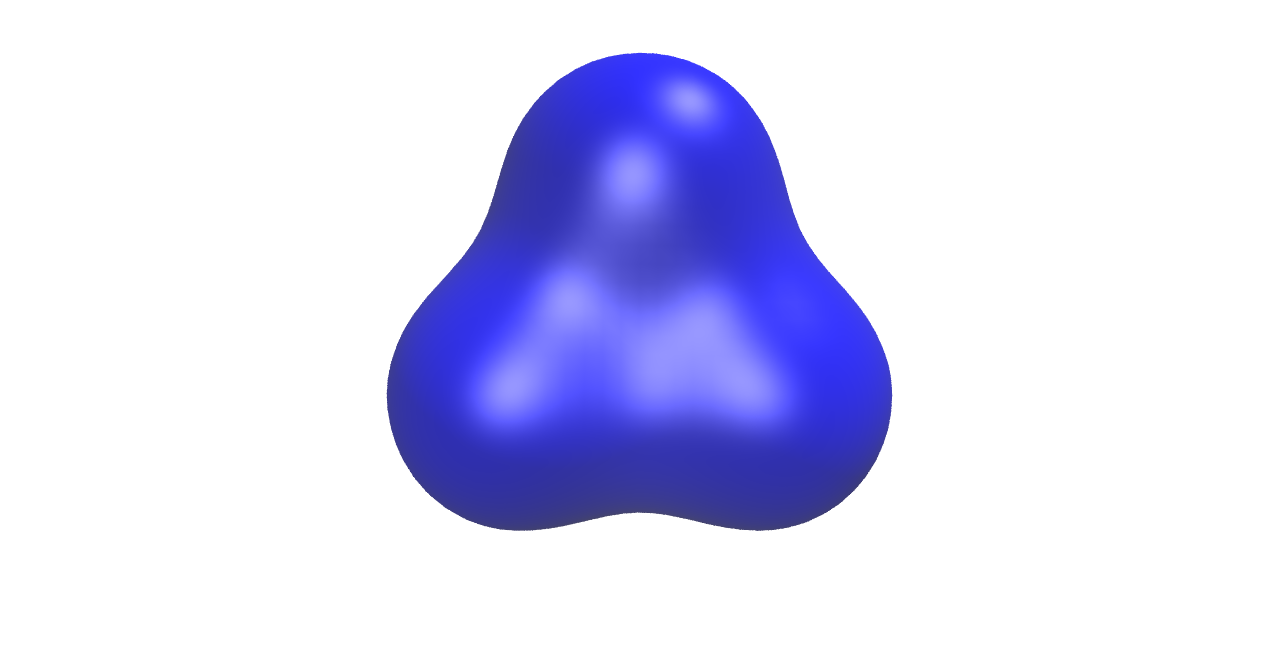}}
\subcaptionbox{}
{\includegraphics[scale=0.7,trim={2cm 0 2cm 0},clip]{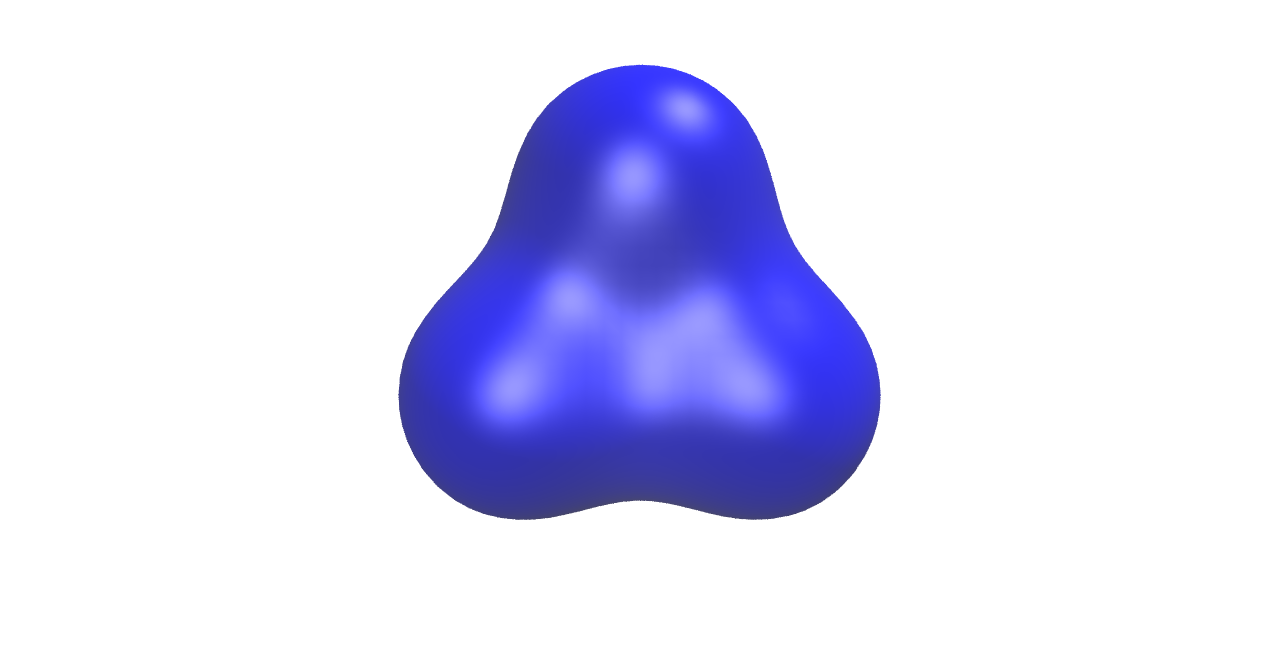}}
\subcaptionbox{}
{\includegraphics[scale=0.7,trim={2cm 0 2cm 0},clip]{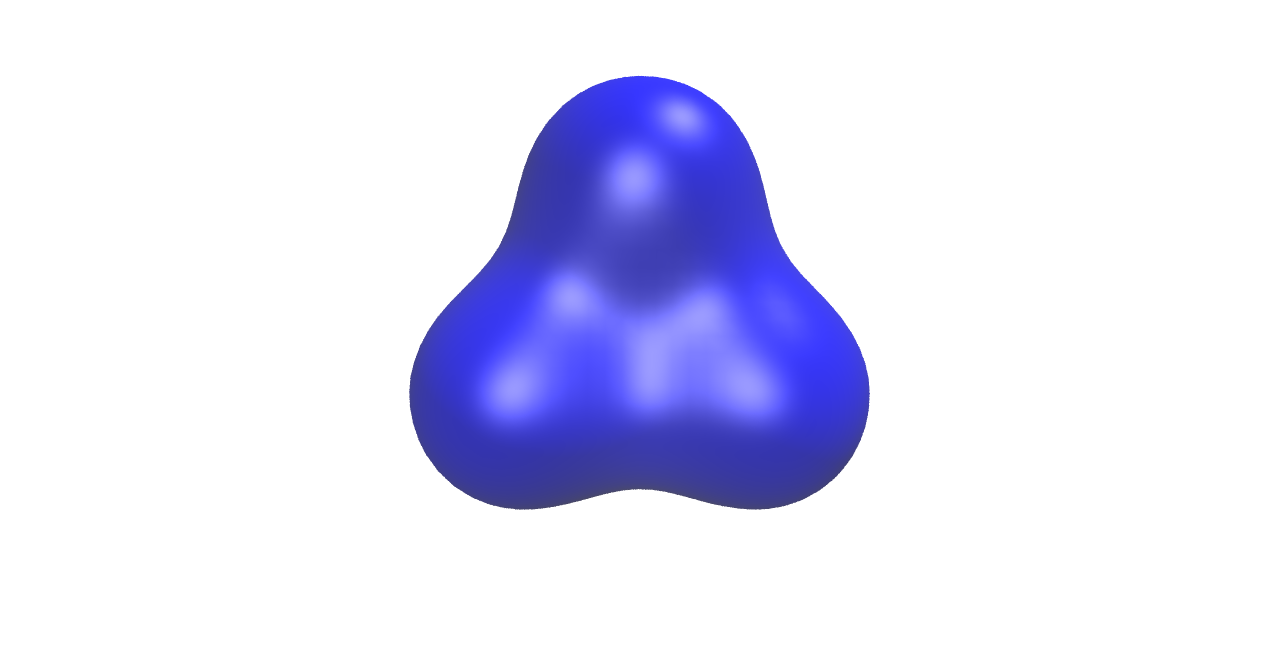}}\\
\subcaptionbox{}
{\includegraphics[scale=0.7,trim={2cm 0 2cm 0},clip]{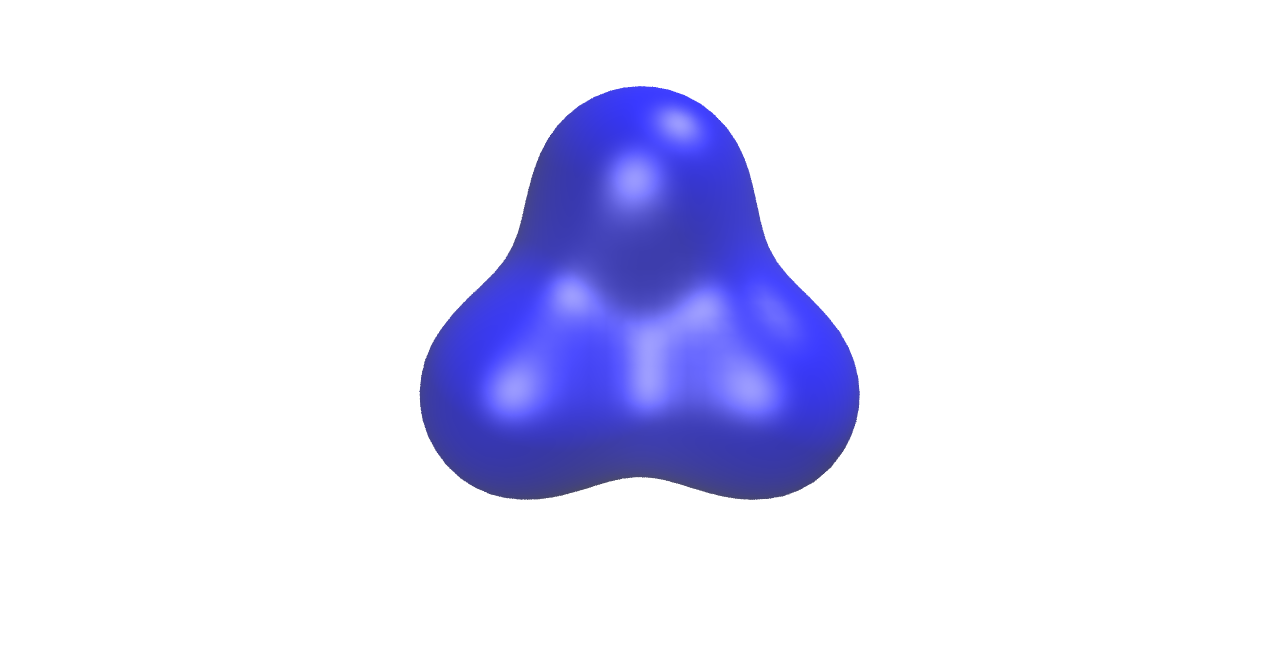}}
\subcaptionbox{}
{\includegraphics[scale=0.7,trim={2cm 0 2cm 0},clip]{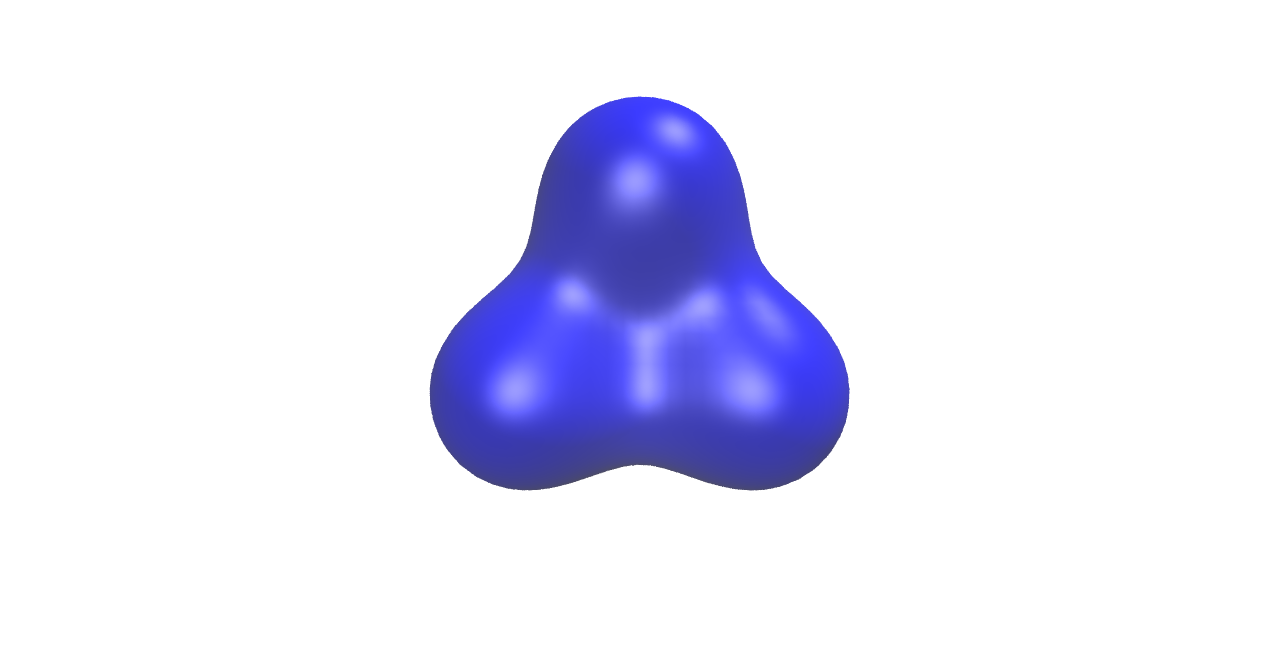}}
\subcaptionbox{}
{\includegraphics[scale=0.7,trim={2cm 0 2cm 0},clip]{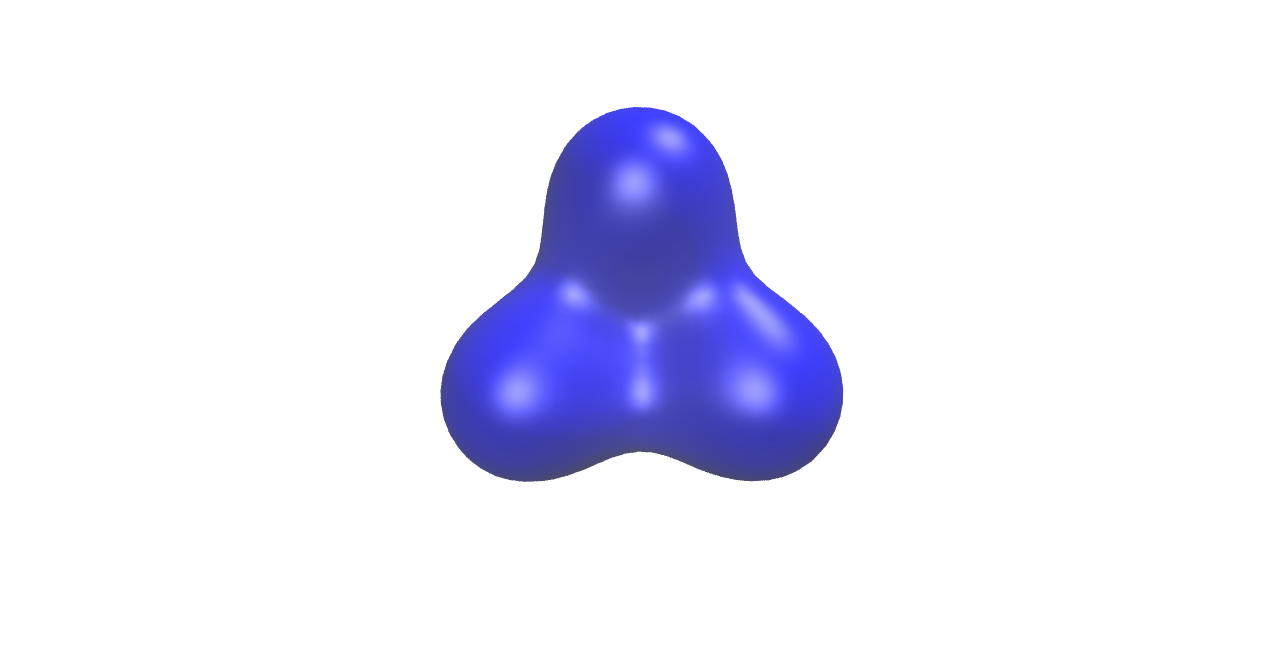}}
\end{tabular}
\caption{The isosurfaces of the imaginary three-atom-molecule generated with the Clifford-Fourier transform method using the Gaussian initial data defined in Equation \ref{eq:gaussianinitial} and extracted at propagation time $t=10^4$ and different isovalues.(a) The isosurface extracted at isovalue $=0.5$. (b) isovalue $=0.6$. (c) isovalue $=0.7$.  (d) isovalue $=0.8$. (e) isovalue $=0.9$. (f) isovalue $=1.0$. }
\label{fig:gaussian10000isovalue}
\end{figure}

\subsubsection{The effect of propagation time}
 We start our experiments by investigating the effect of propagation time on the extracted isosurfaces. We apply our algorithm on a piecewise initial shape of the imaginary three-atom-molecule. All extracted isosurfaces has isovalue $=0.9$ while time propagates in powers of $10$ as $t=10^1,t=10^2,t=10^3,t=10^4,t=10^5$. As seen in Figure \ref{fig:pwtimepropagation}, as time propagates the isosurfaces become more smooth, and geometric singularities disappear. With this being said, it is not necessarily true that the higher the propagation time the better the surface is. It is clear that in Figure \ref{fig:pwtimepropagation}(f) the isosurface is very smooth and not very helpful.

Then, the same experiment is done to the Gaussian initial shape with propagation times $t=10^1,t=10^2,t=10^3,t=10^4,$ and $ t=10^5$ and all the extracted isosurfaces has isovalue $1.0$. Figure \ref{fig:gaussiantimepropagation} shows that as time propagates the isosurfaces get more smooth and the geometric singularities disappear. However, this does not mean increasing propagation time is always better. As you can see that Figure \ref{fig:gaussiantimepropagation}(f) shows a surface that is over smoothed and hence not very useful.


\subsubsection{The effect of the isovalue}
 
After seeing the effect of the propagation time, we now show the effect of changing the isovalue. As seen in Figure \ref{fig:pw100isovalue}, we have six different isosurfaces extracted at the following isovalues $0.4,0.5,0.6,0.7,0.8,$ and $0.9$ respectively at propagation time $t=10^2$ using the piecewise initial data defined in Equation \ref{eq:pwinitial}. Due to the definition of the initial data, the surfaces get inflated as the isovalue increases. In terms of geometric singularities, there is no significant changes as the isovalue changes. Likewise, Figure \ref{fig:pw10000isovalue} shows isosurfaces extracted at propagation time $t=10^4$ and the same isovalues $0.4,0.5,0.6,0.7,0.8,$ and $0.9$ using the piecewise initial data. The same observations can be said about this figure as well, where isosurfaces get larger as isovalue increases and no noticeable changes happen to geometric characteristics. In both figures, we do not have geometric singularities.  

\begin{figure}[h!]
\center
\begin{tabular}{ c c c }
\subcaptionbox{}
{\includegraphics[scale=0.17,trim={6cm 0 6cm 0},clip]{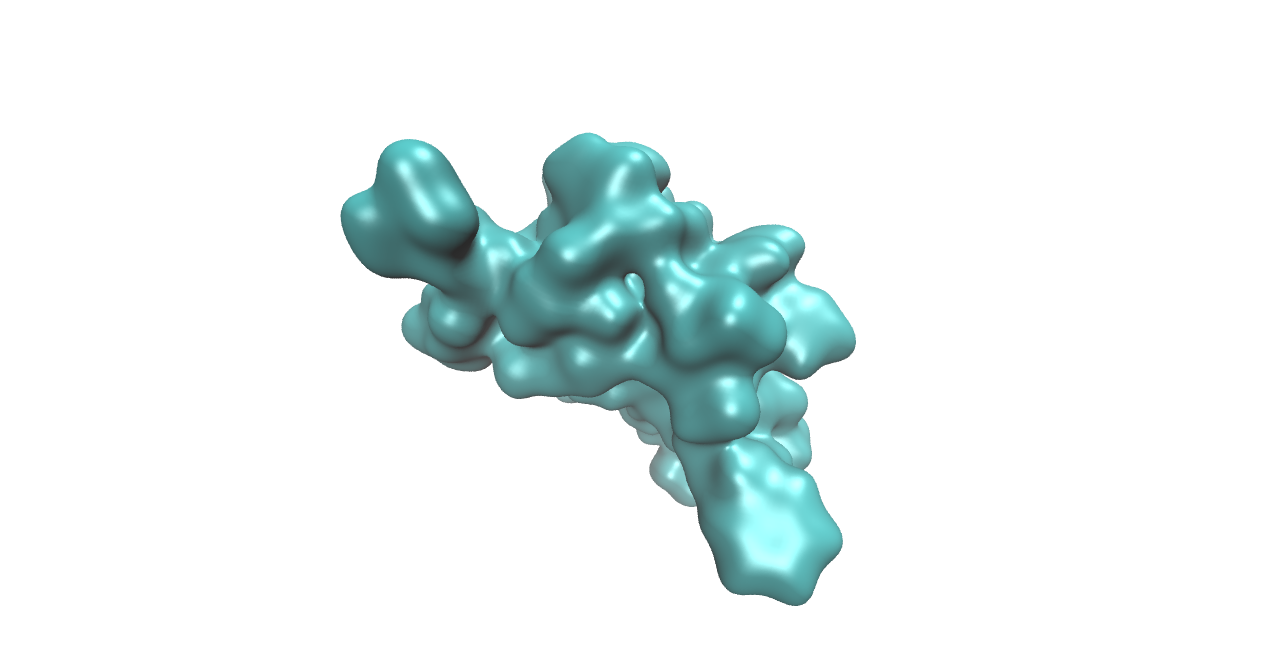}}
\subcaptionbox{}
{\includegraphics[scale=0.17,trim={6cm 0 6cm 0},clip]{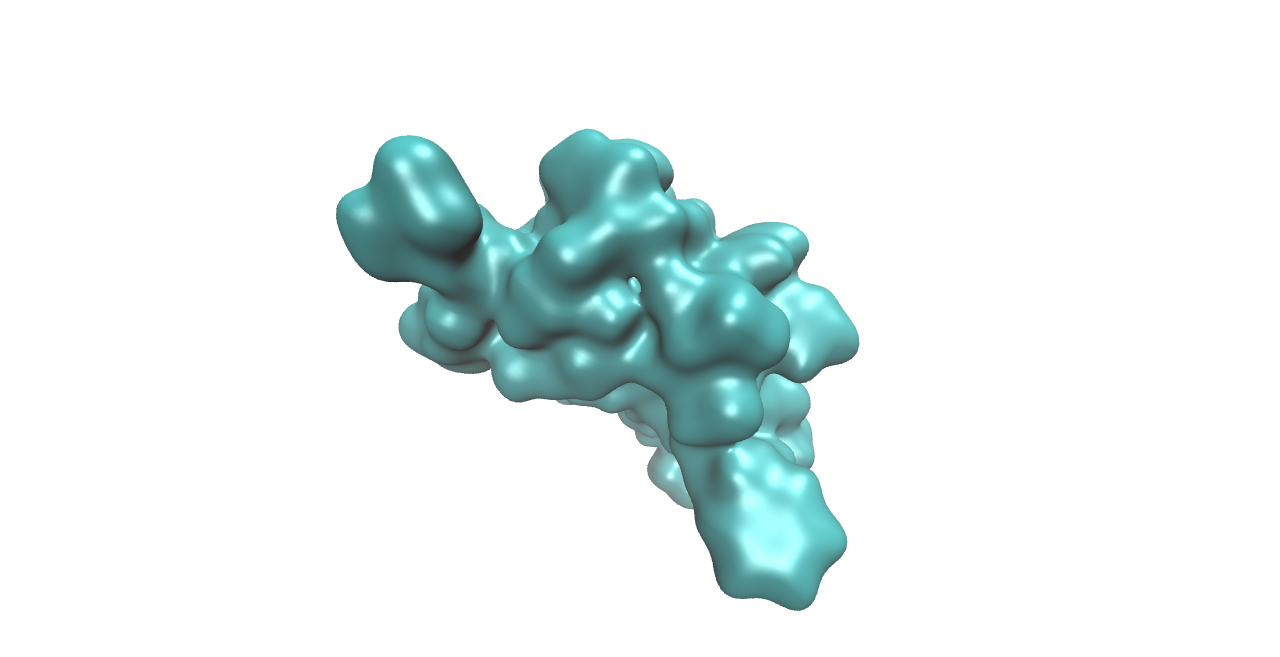}}
\subcaptionbox{}
{\includegraphics[scale=0.17,trim={6cm 0 6cm 0},clip]{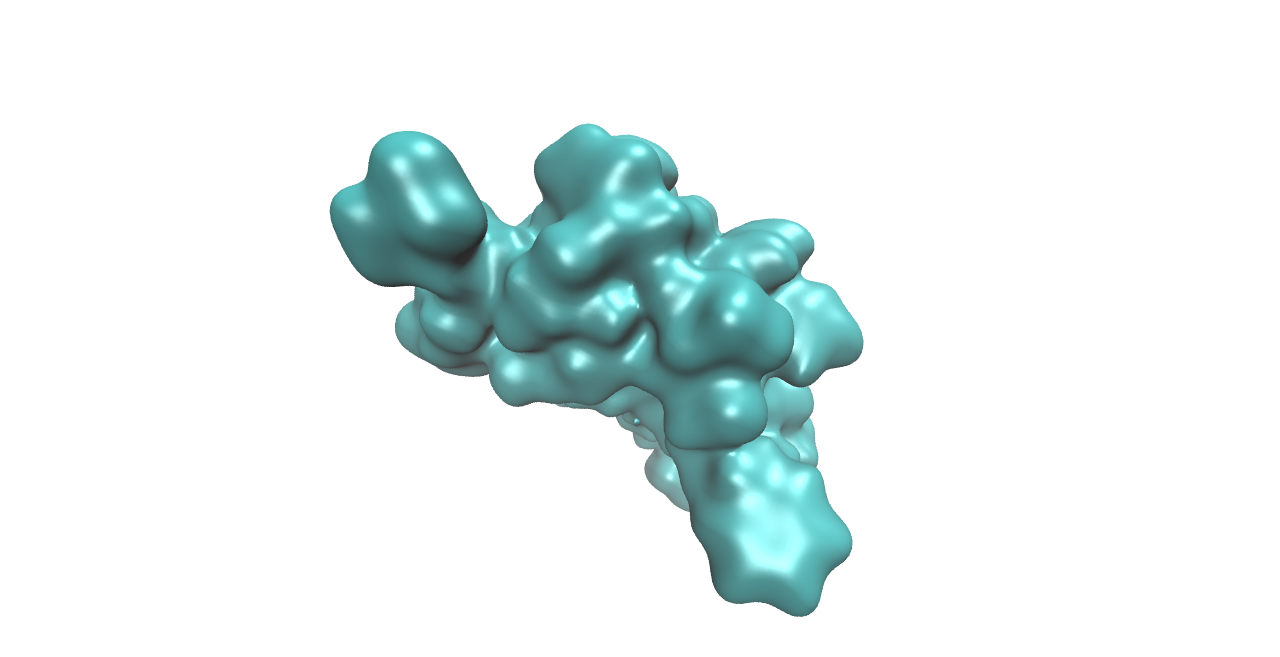}}\\
\subcaptionbox{}
{\includegraphics[scale=0.17,trim={6cm 0 6cm 0},clip]{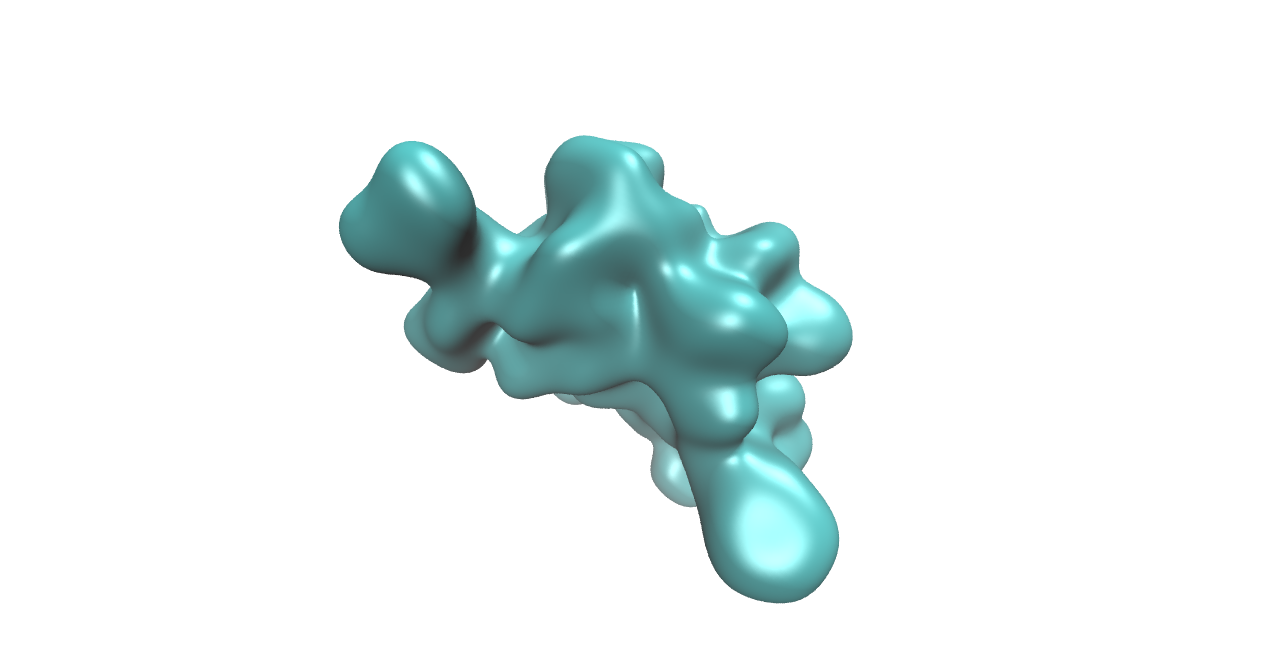}}
\subcaptionbox{}
{\includegraphics[scale=0.17,trim={6cm 0 6cm 0},clip]{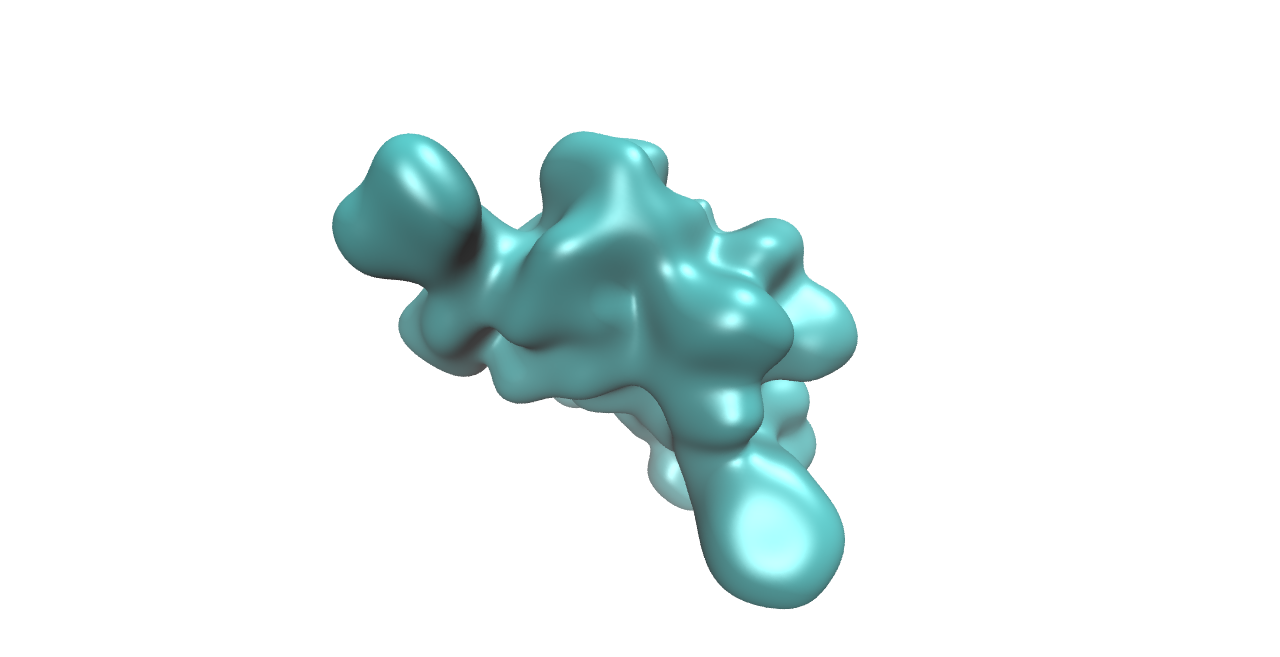}}
\subcaptionbox{}
{\includegraphics[scale=0.17,trim={6cm 0 6cm 0},clip]{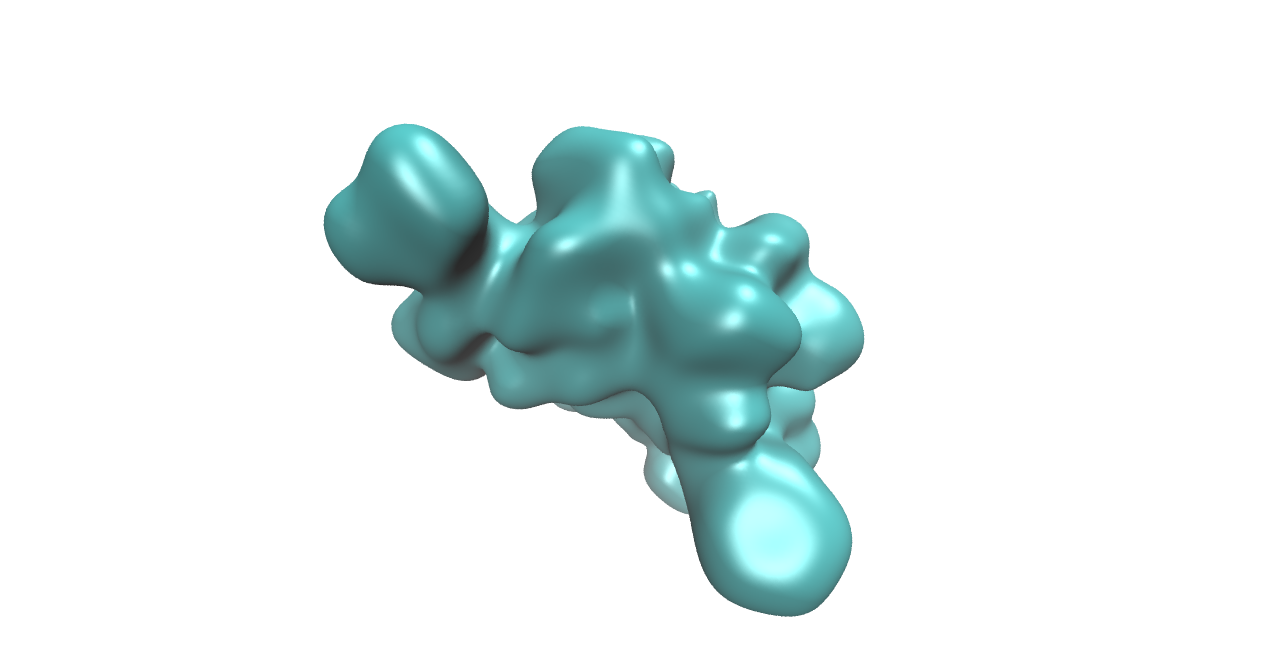}}
\end{tabular}
\caption{The molecular surfaces of protein 1ajj generated by Clifford-Fourier transform method using piecewise initial data. (a) isovalue $=0.7$, $t=10^2$, (b) isovalue $=0.8$, $t=10^2$, (c) isovalue $=0.9$, $t=10^2$, (d) isovalue $=0.7$, $t=10^4$, (e) isovalue $=0.8$, $t=10^4$, (f) isovalue $=0.9$, $t=10^4$.}
\label{fig:1ajj_pw}
\end{figure}

Now, we conduct the same experiment on the Gaussian initial shape. Figure \ref{fig:gaussian100isovalue} shows six isosurfces extracted at the following isovlaues $0.5,0.6,0.7,0.8,0.9$ and $1.0$ respectively at propagation time $t=10^2$ using the Gaussian initial data defined in Equation \ref{eq:gaussianinitial}. In contrary to the piecewise case, the isosurfaces of the Gaussian shape get deflated as the isovalues get larger. That is also due to the definition of Gaussian initial equation which is monotonically decreasing. Regarding the geometric characteristics, it is clear that the isosurfaces tend to be more meaningful as the isovalues get larger. However, very large isovalues may cause some geometric singularities. Likewise, the same experiment has been done for the same settings but with propagation time $t=10^4$ and is demonstrated in Figure \ref{fig:gaussian10000isovalue}. Clearly, the same conclusions can be made about the isosurfacs in terms of surface size as well as geometric characteristics and singularities.

\subsection{Biomolecular surfaces}
In this part, we show the biomolecular surfaces of real proteins generated using our algorithm exploiting the Clifford-Fourier transform to validate our proposed method. To gain acceptance within the molecular visualization community, our method is compared to the well-established method SES using MSMS package \cite{sanner1996reduced} available in the software visual molecular dynamics (VMD) \cite{HUMP96}. All protein structures and atomic coordinates used in our computations were obtained from the Protein Data Bank (PDB) website (https://www.rcsb.org). We then used the package PDB2PQR\cite{dolinsky2004pdb2pqr} to add the missing hydrogen atoms and assign point charges at atomic centers based on the CHARMM force field \cite{mackerell1998all}.

\begin{figure}[h!]
\center
\begin{tabular}{ c c c }
\subcaptionbox{}
{\includegraphics[scale=0.17,trim={6cm 0 6cm 0},clip]{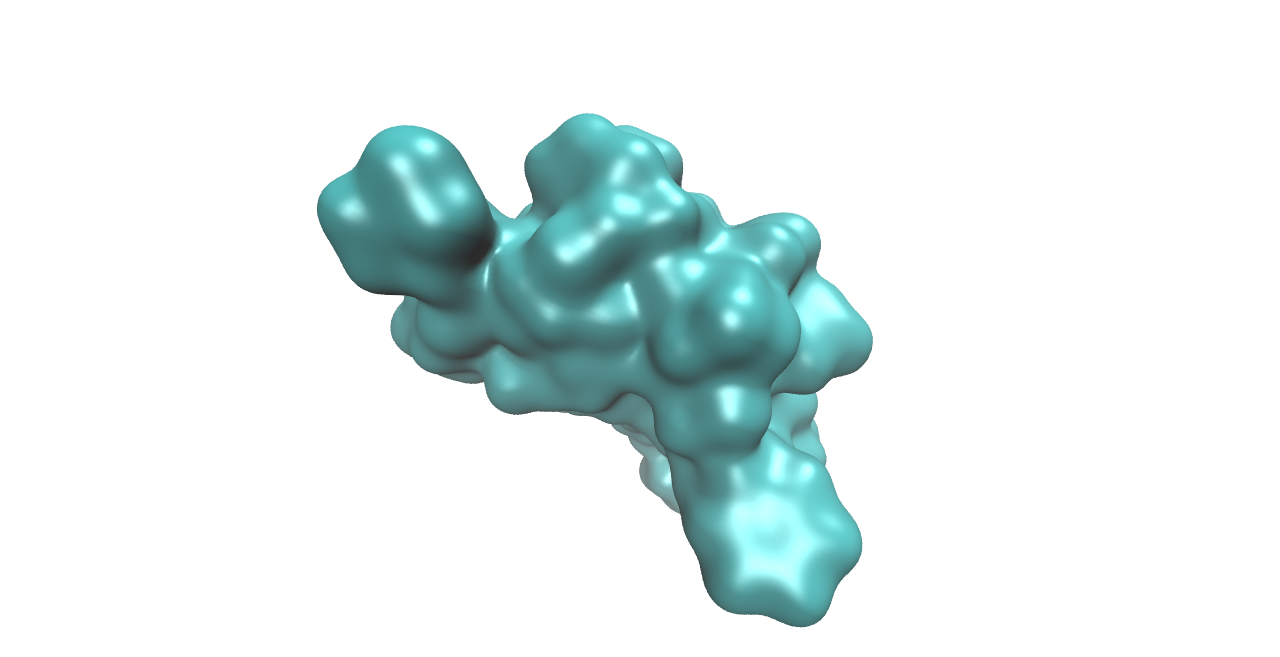}}
\subcaptionbox{}
{\includegraphics[scale=0.17,trim={6cm 0 6cm 0},clip]{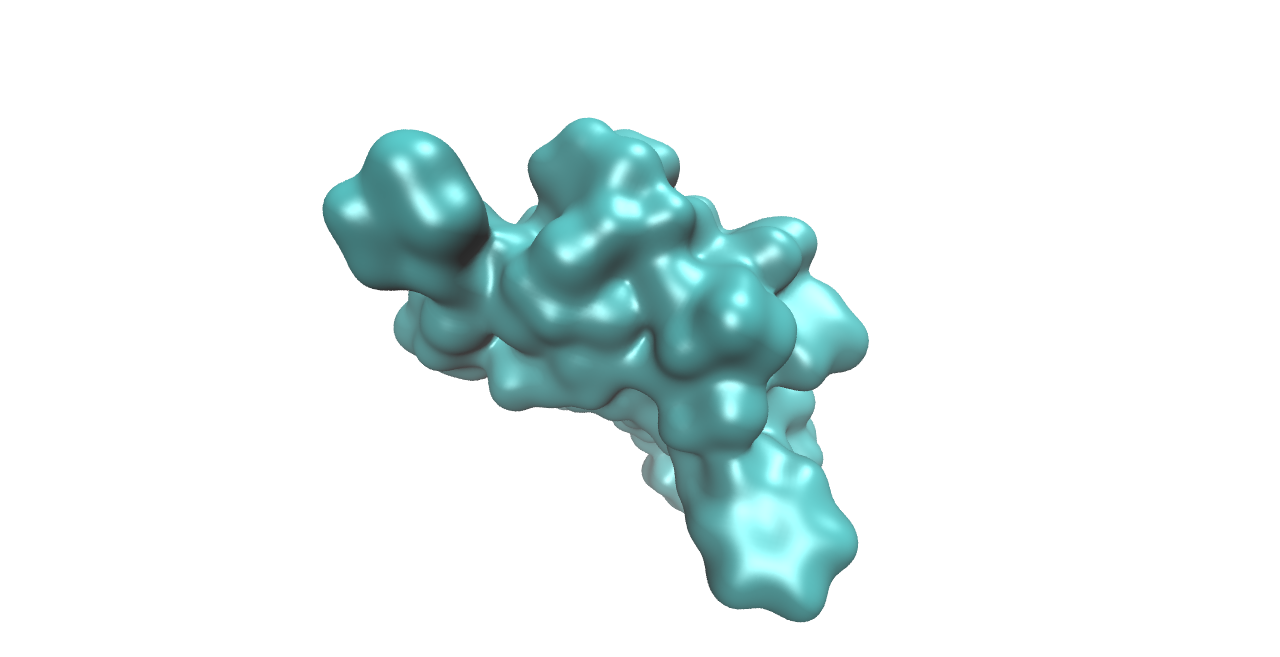}}
\subcaptionbox{}
{\includegraphics[scale=0.17,trim={6cm 0 6cm 0},clip]{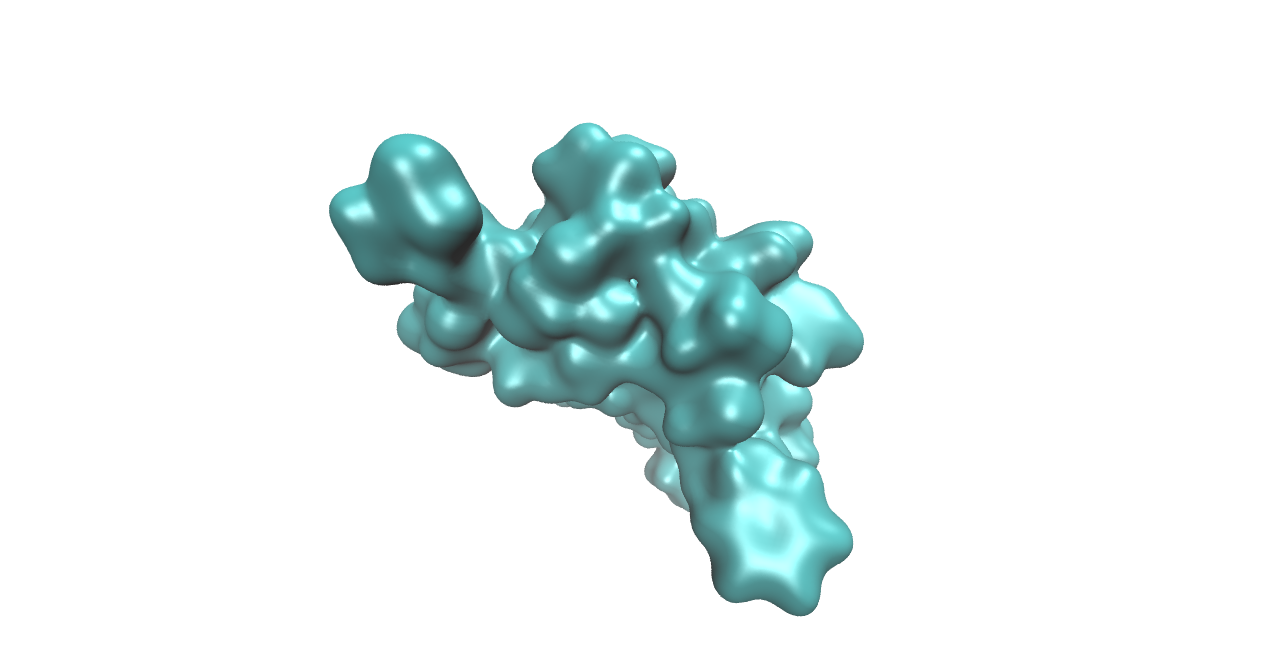}}\\
\subcaptionbox{}
{\includegraphics[scale=0.17,trim={6cm 0 6cm 0},clip]{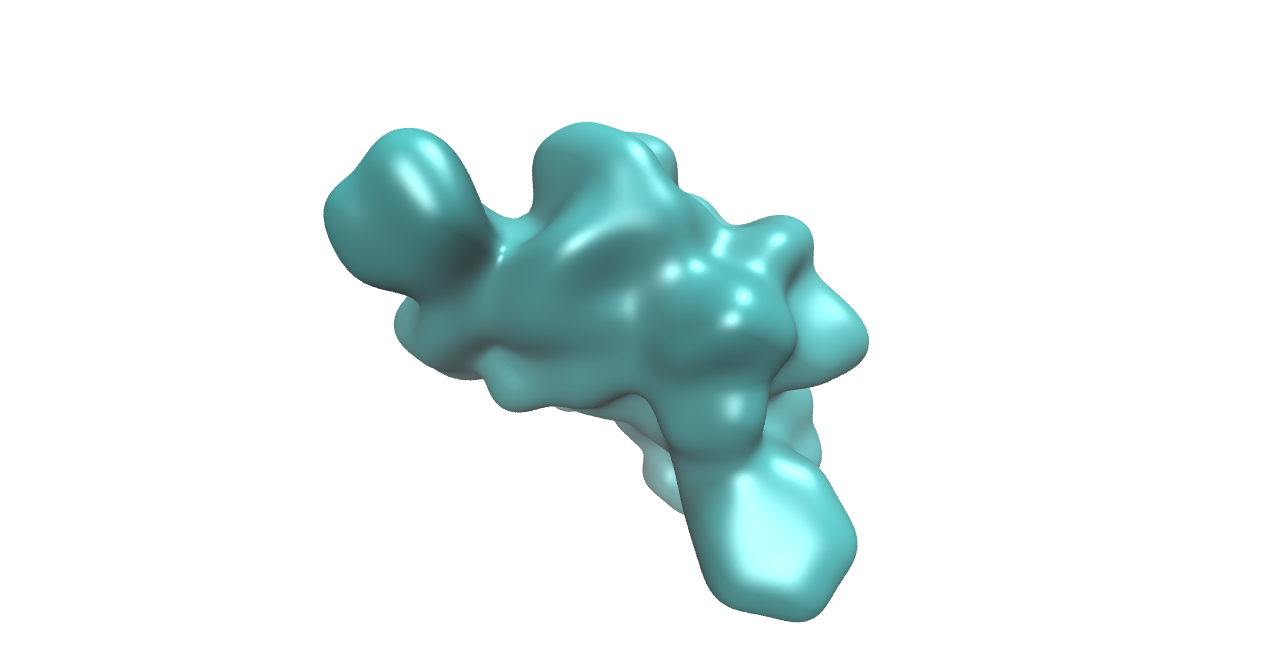}}
\subcaptionbox{}
{\includegraphics[scale=0.17,trim={6cm 0 6cm 0},clip]{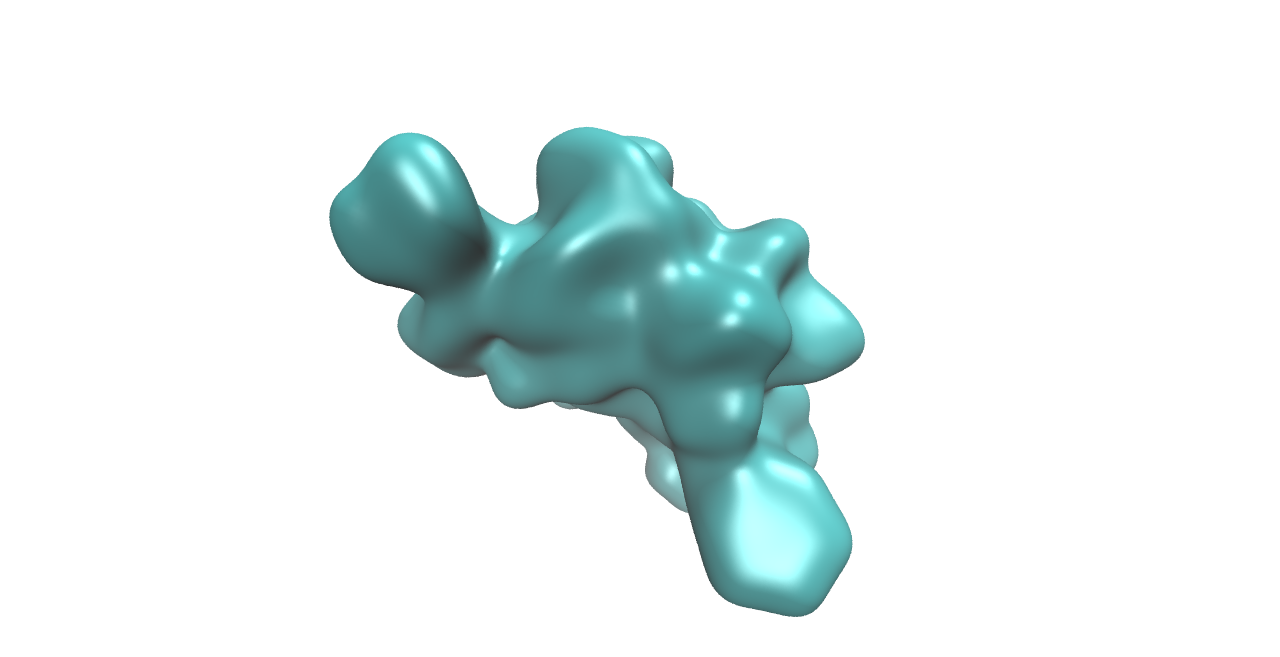}}
\subcaptionbox{}
{\includegraphics[scale=0.17,trim={6cm 0 6cm 0},clip]{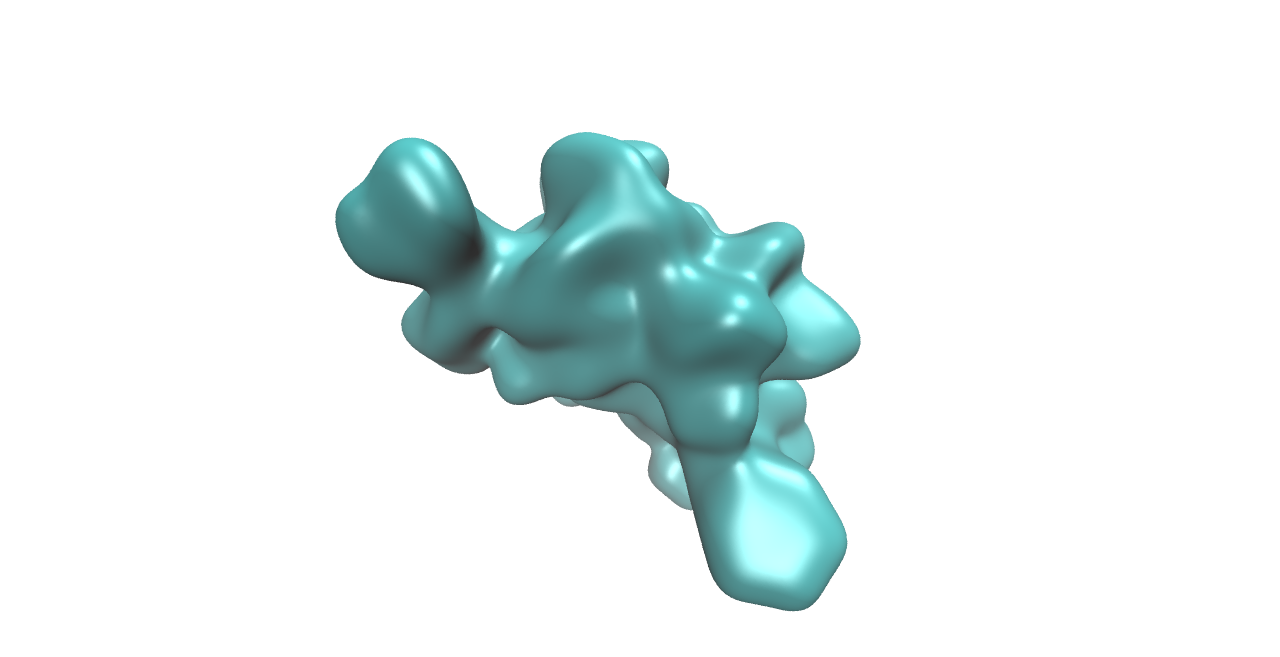}}
\end{tabular}
\caption{The molecular surfaces of protein 1ajj generated by Clifford-Fourier transform method using Gaussian initial data. (a) isovalue $=0.7$, $t=10^2$, (b) isovalue $=0.8$, $t=10^2$, (c) isovalue $=0.9$, $t=10^2$, (d) isovalue $=0.7$, $t=10^4$, (e) isovalue $=0.8$, $t=10^4$, (f) isovalue $=0.9$, $t=10^4$.}
\label{fig:1ajj_gaussian}
\end{figure}

Now, we start with the protein 1ajj by applying our method with both cases of initial shapes: piecewise and Gaussian. Figure \ref{fig:1ajj_pw} shows the piecewise initial shape where the first row has propagation time $t=10^2$ and the second row has propagation time $t=10^4$. In each row, we have the following isovalues $0.7,0.8,$ and $0.9$ respectively. As expected from the above test cases, the isosurfaces get inflated slightly as the isovalues increase. Also as expected, the geometric singularities do not disappear as we change the isovalues. The second row in Figure \ref{fig:1ajj_pw}, which is for $t=10^4$, goes in agreement with our predictions in the above test cases experiments as well. The isosurfaces for propagation time $t=10^4$ are very smooth and hence not preferred for biomolecular surfaces.

Now, we carry out the same experiments on protein 1ajj with Gaussian initial data. This is to investigate the combinations of propagation times and isovalues. Likewise, Figure \ref{fig:1ajj_gaussian} shows isosurfaces extracted at isovalues $0.7,0.8,0.9$ respectively where the first row show the results at propagation time $t=10^2$ and the second row show the results at $t=10^4$. Isosurfaces demonstrated affirm our predictions made earlier where the geometric singularities tend to appear as the isovalues increase. Moreover, the surfaces get deflated as the isovalues increase due to the definition of the Gaussian function. Picking the appropriate isovalue has two competing factors which are the surface volume and the presence of geometric singularities.   
\begin{figure}[h!]
\center
\begin{tabular}{M{0.5cm} | M{4.5cm} | M{4.5cm} | M{4.5cm} }
 & CFT (piecewise) & CFT (Gaussian) & MSMS \\
 & isovalue$=0.9$, $t=10^2$ & isovalue$=0.8$, $t=10^2$& density$=9.5$, probe$=1.5$ \\
\hline
\begin{turn}{-90}1ajj\end{turn} &
\subcaptionbox{}
{\includegraphics[scale=0.15,trim={6cm 0 6cm 0},clip]{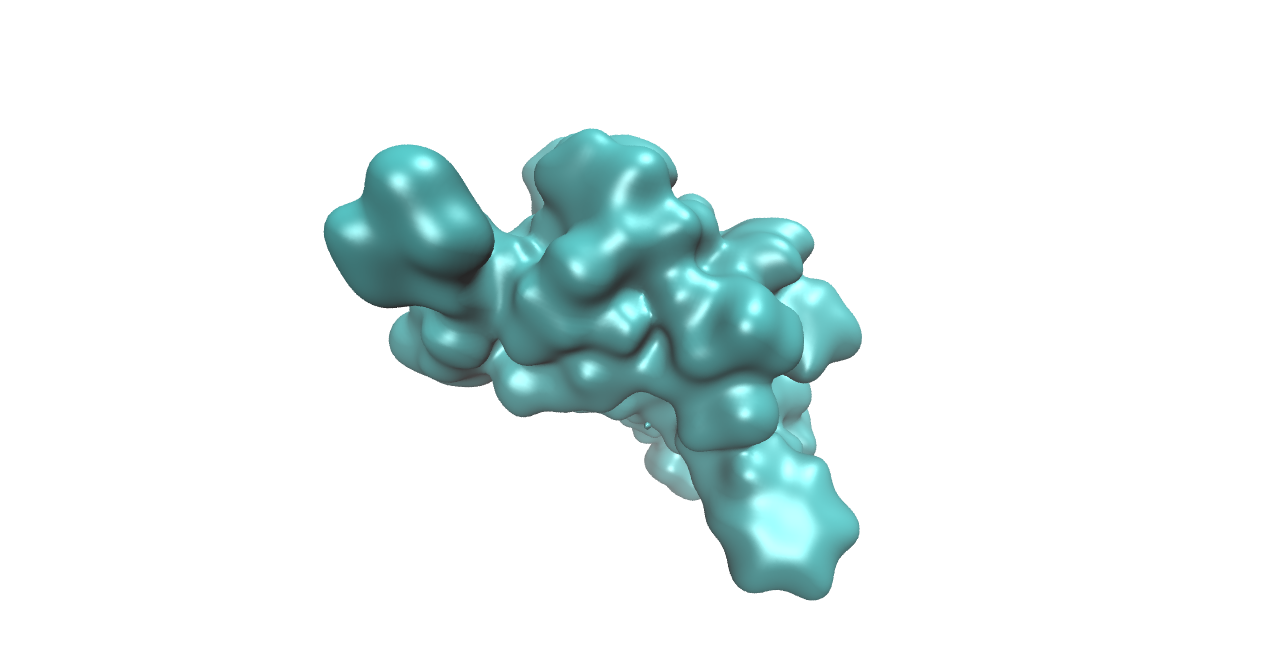}} &
\subcaptionbox{}
{\includegraphics[scale=0.15,trim={6cm 0 6cm 0},clip]{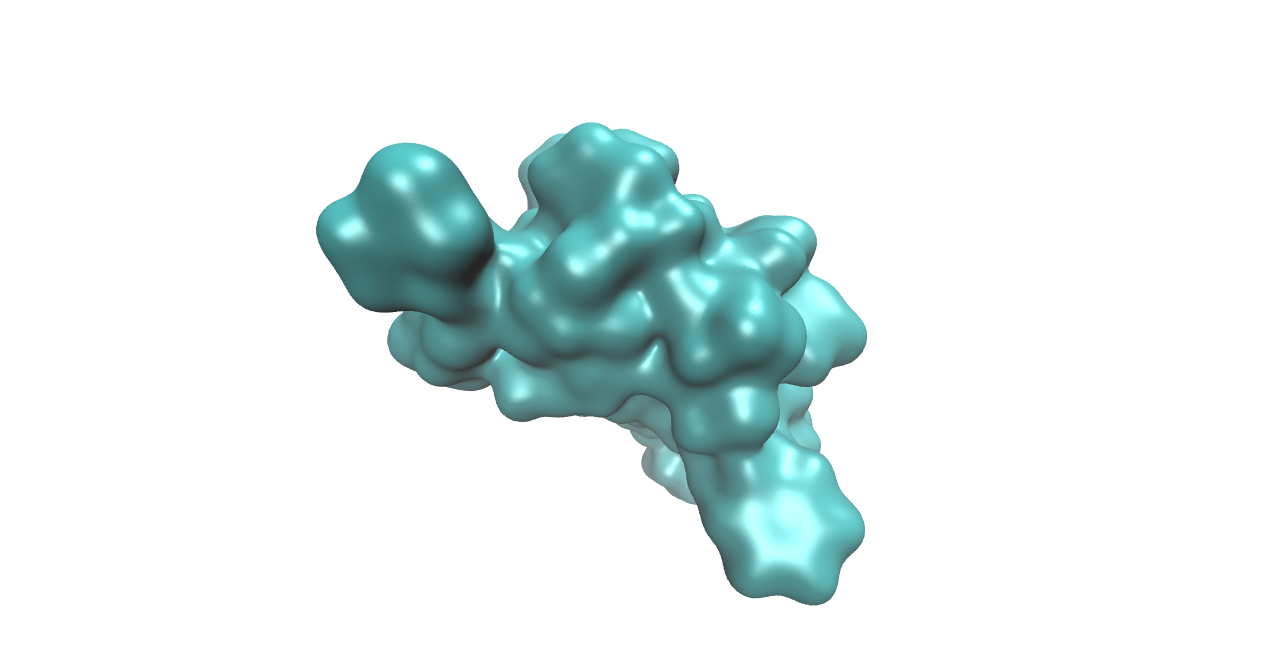}}&
\subcaptionbox{}
{\includegraphics[scale=0.15,trim={6cm 0 6cm 0},clip]{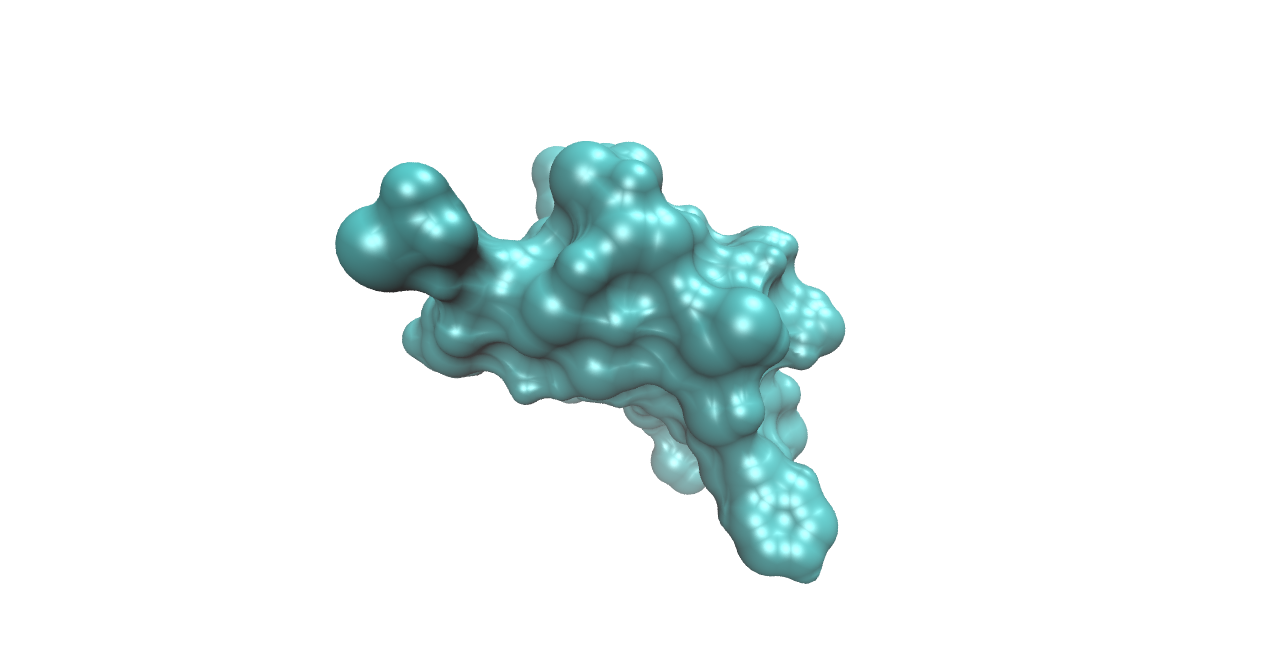}}\\
\hline
\begin{turn}{-90}1bor \end{turn} &
\subcaptionbox{}
{\includegraphics[scale=0.15,trim={6cm 0 6cm 0},clip]{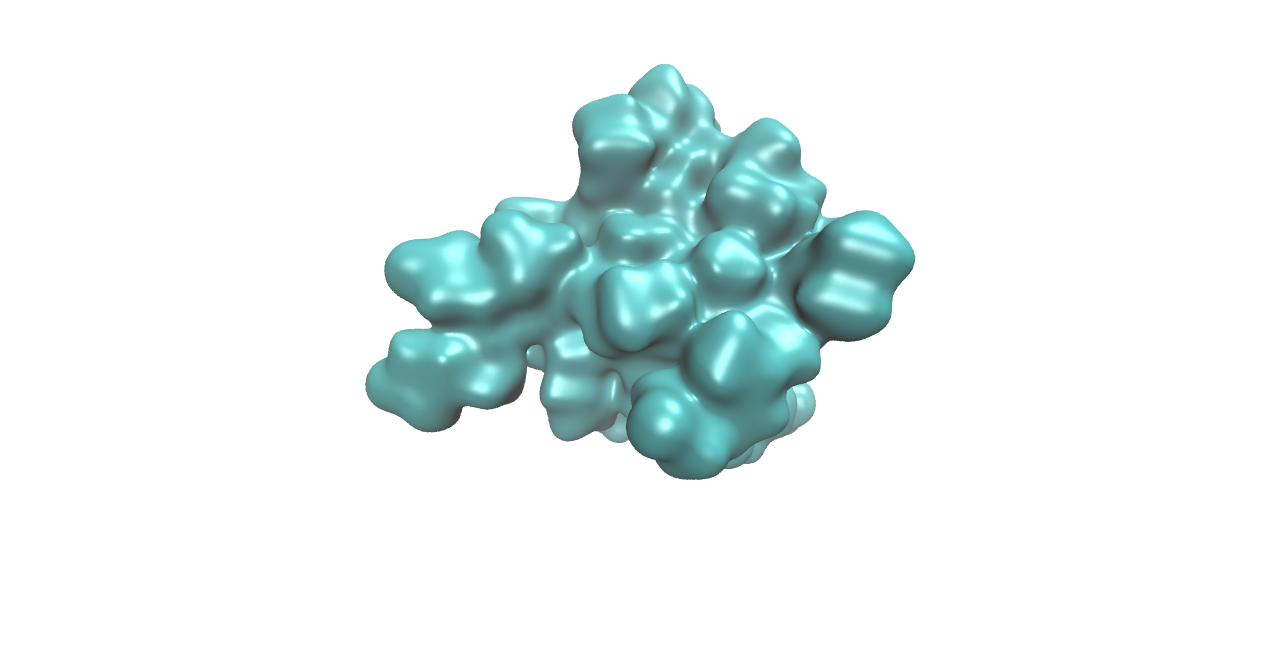}} &
\subcaptionbox{}
{\includegraphics[scale=0.15,trim={6cm 0 6cm 0},clip]{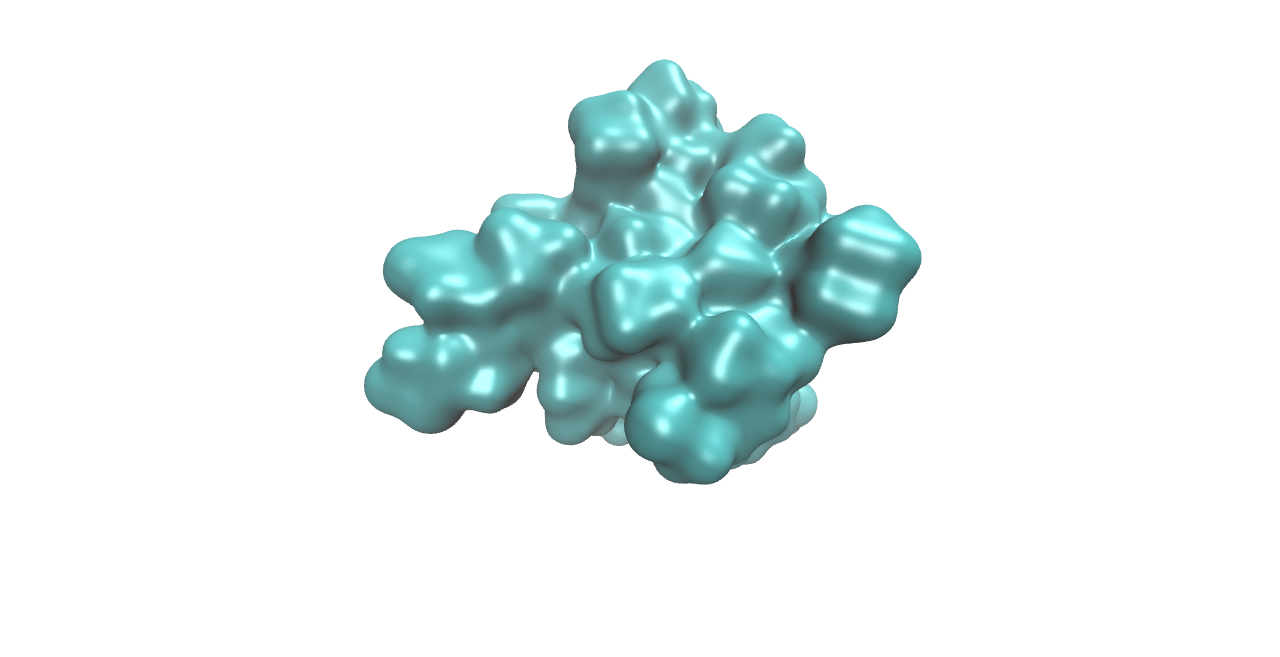}} &
\subcaptionbox{}
{\includegraphics[scale=0.15,trim={6cm 0 6cm 0},clip]{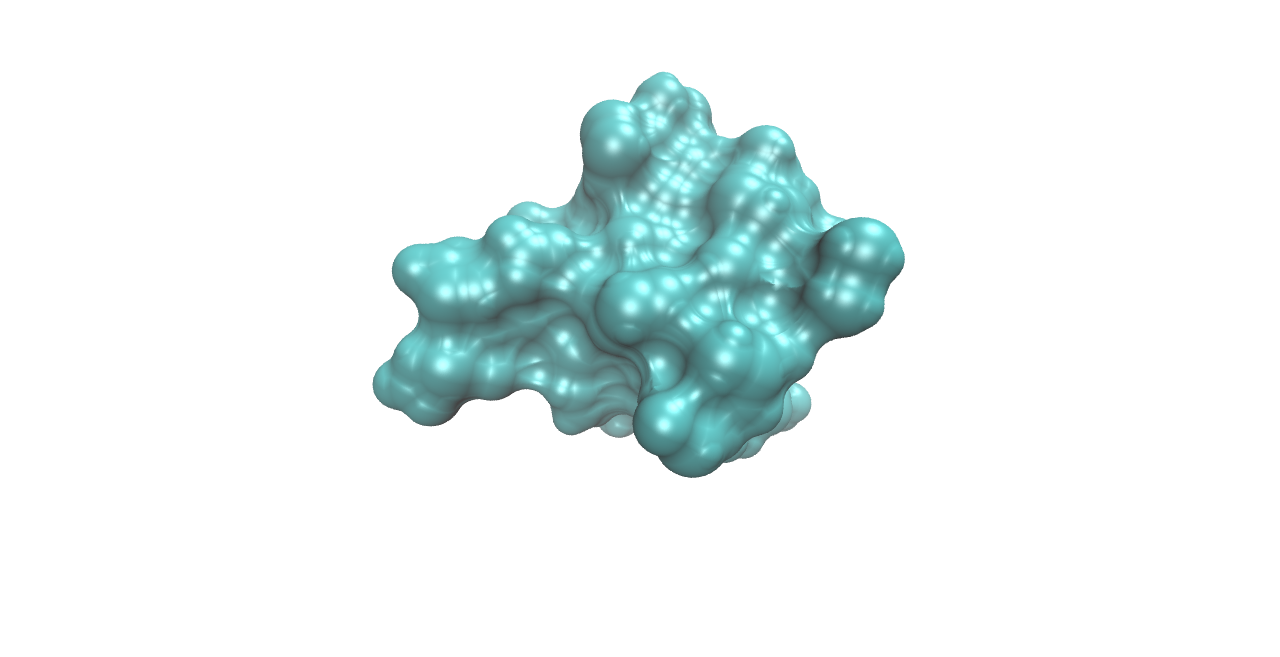}}\\
\hline
\begin{turn}{-90}1mbg \end{turn}&
\subcaptionbox{}
{\includegraphics[scale=0.15,trim={6cm 0 6cm 0},clip]{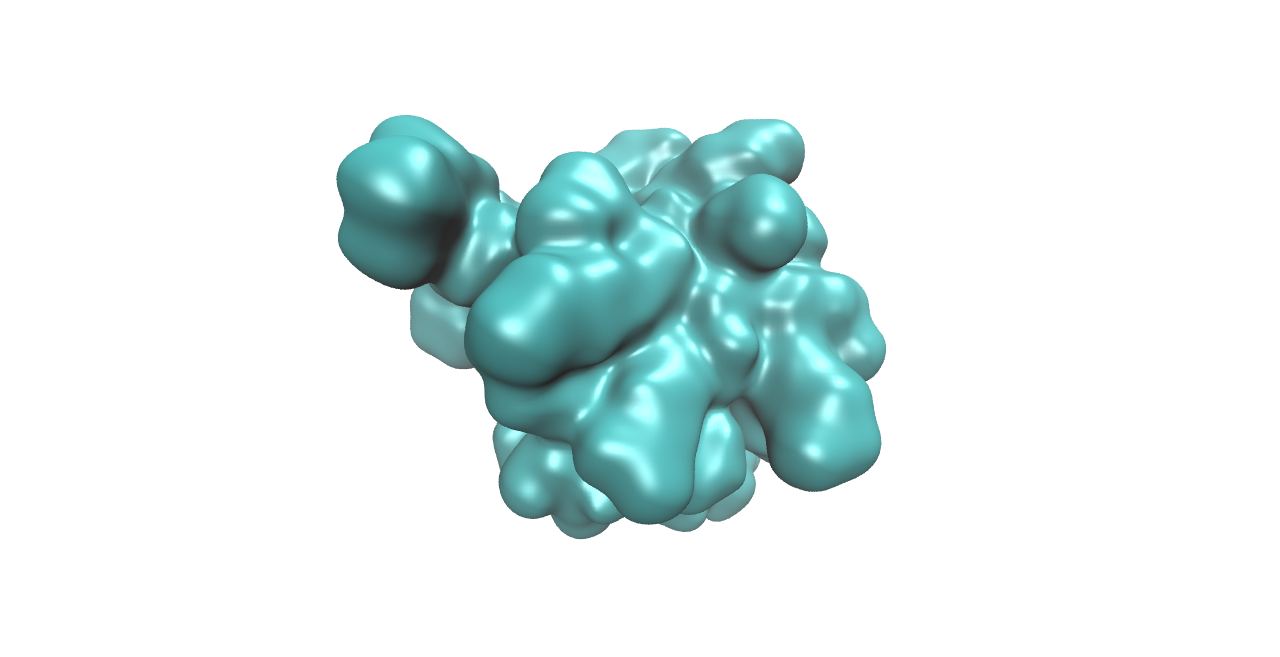}} &
\subcaptionbox{}
{\includegraphics[scale=0.15,trim={6cm 0 6cm 0},clip]{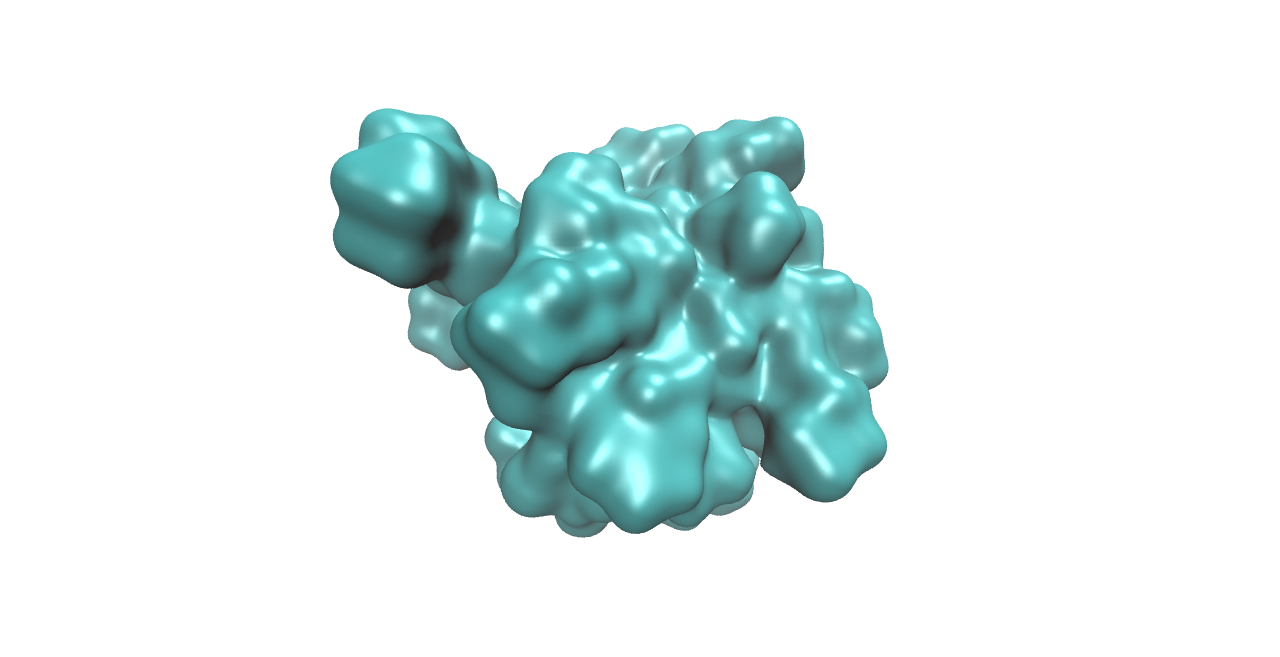}} &
\subcaptionbox{}
{\includegraphics[scale=0.15,trim={6cm 0 6cm 0},clip]{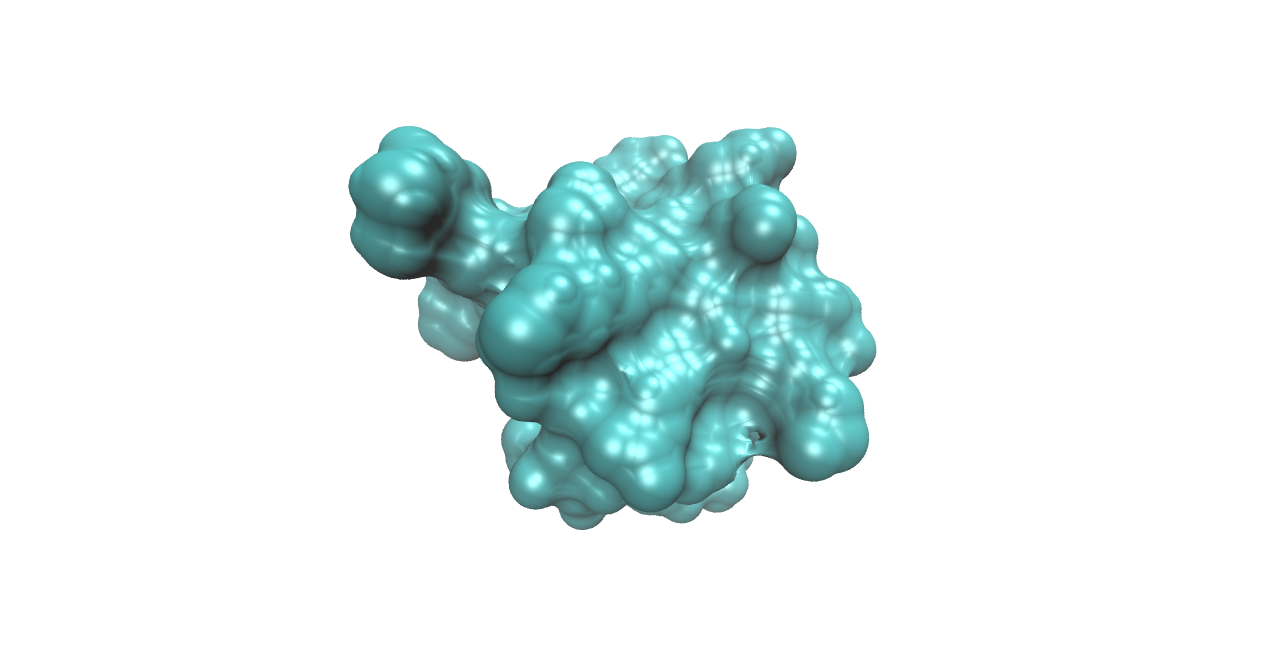}}\\
\hline
\begin{turn}{-90}1sh1 \end{turn}&
\subcaptionbox{}
{\includegraphics[scale=0.15,trim={6cm 0 6cm 0},clip]{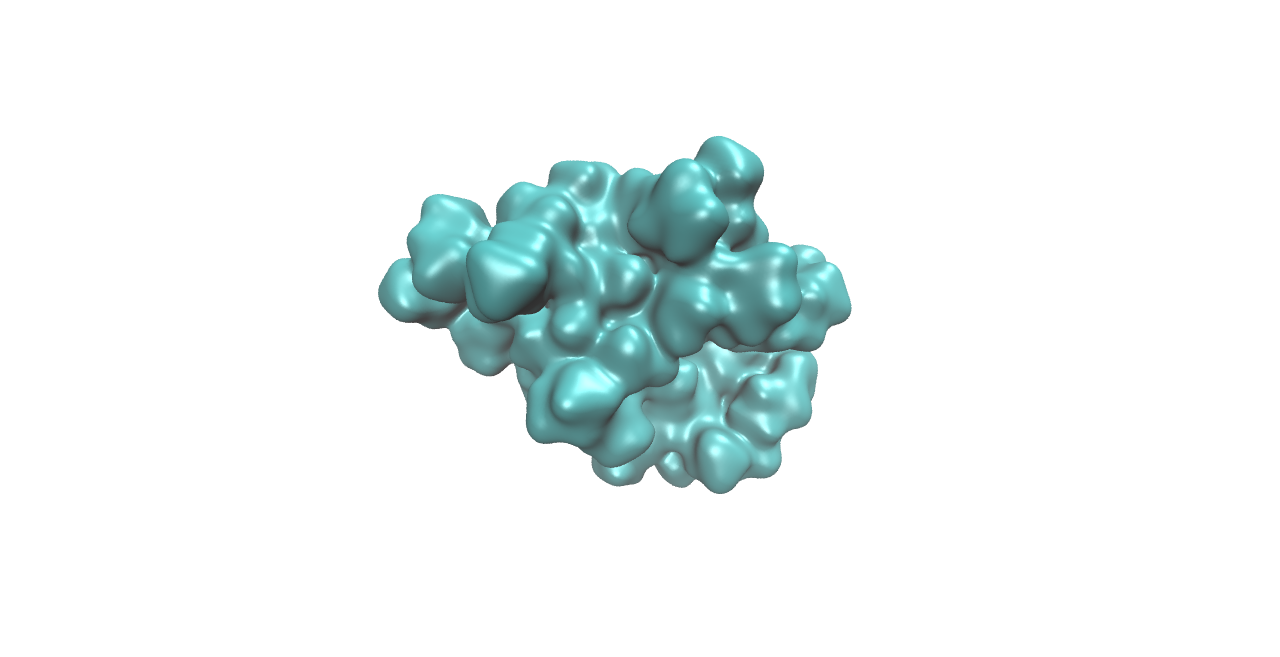}} &
\subcaptionbox{}
{\includegraphics[scale=0.15,trim={6cm 0 6cm 0},clip]{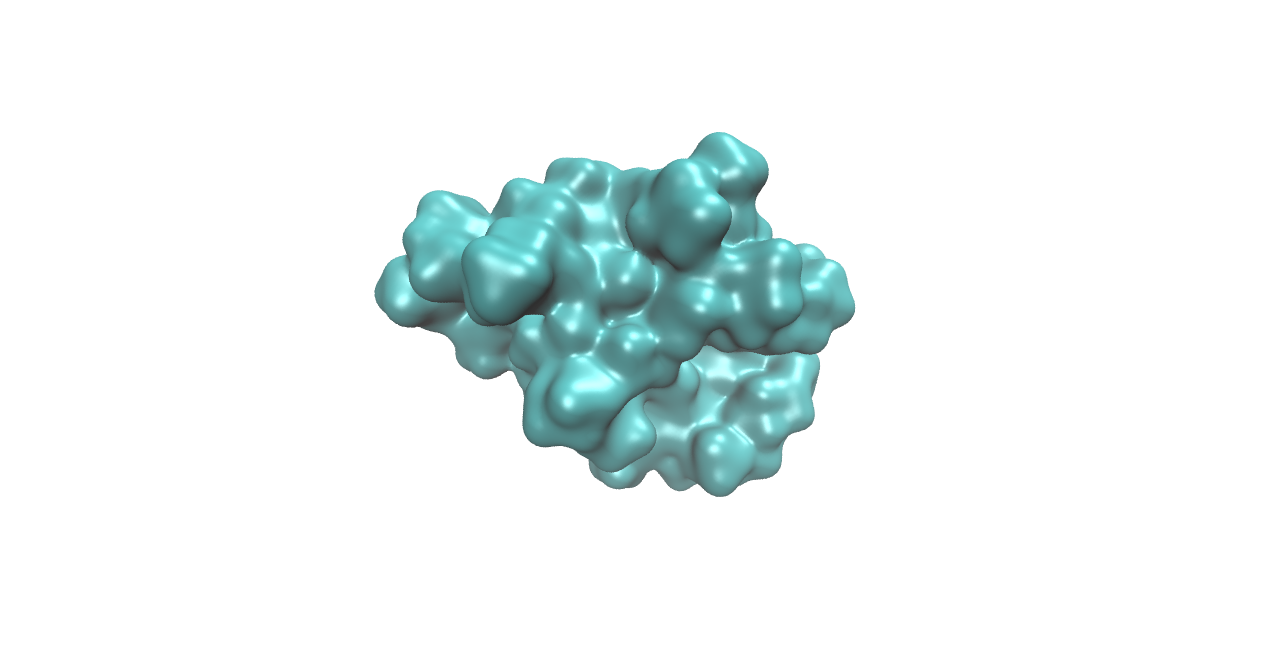}} &
\subcaptionbox{}
{\includegraphics[scale=0.15,trim={6cm 0 6cm 0},clip]{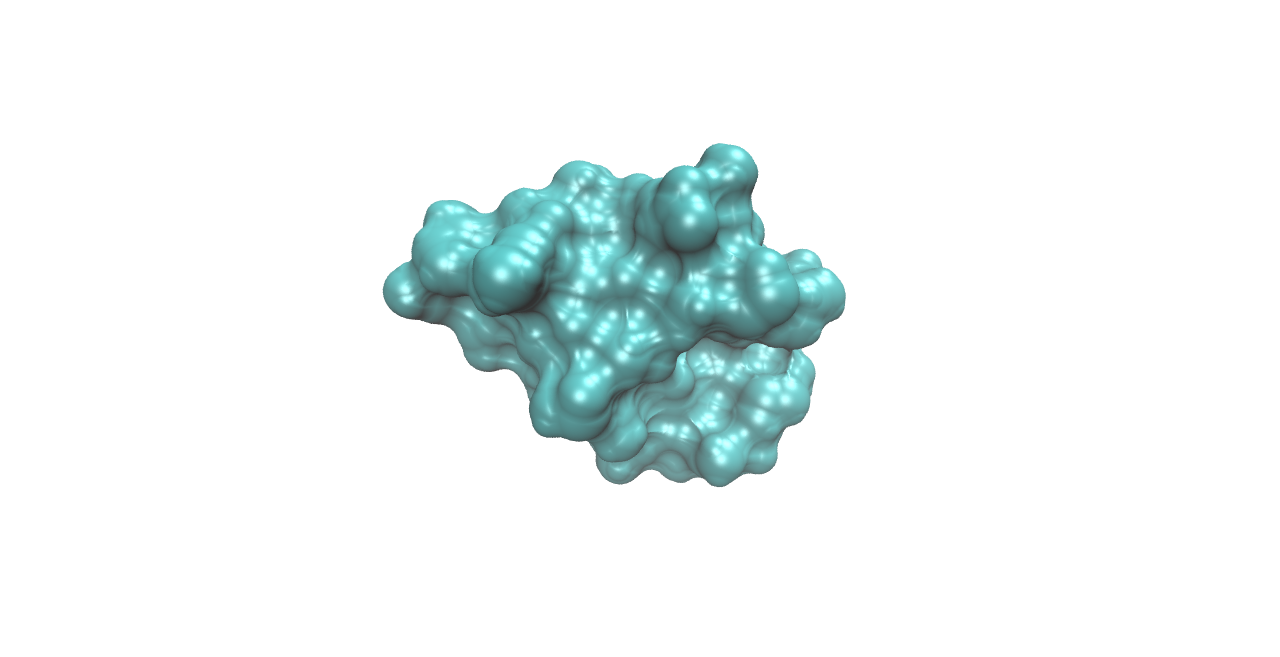}}
\end{tabular}
\caption{Comparison of molecular surfaces generated by Clifford-Fourier transform (CFT)  using piecewise initial data (first column), Clifford-Fourier transform using Gaussian initial data (second column) and MSMS (third column). The first row is for protein 1ajj. The second row is for protein 1bor. The third row is for protein 1mbg. The fourth row is for protein 1sh1. All Clifford-Fourier transform surfaces generated $t=10^2$. }
\label{fig:many_proteins}
\end{figure}
After applying our method to the protein 1ajj and experimenting with different combinations of propagation times, isovalues and initial values, we demonstrate more biomolecular surfaces. We show the biomolecular surfaces of the proteins: 1ajj, 1bor, 1mbg, and 1sh1 using our method with two different sets of parameters and initial values. For the piecewise initial shape, we  choose isovalue $0.9$, the propagation time $t=10^2$ and order $2m=12$. For the Gaussian initial shape, we  choose isovalue $0.8$, the propagation time $t=10^2$ and order $2m=12$. These two surfaces of each protein are compared with SES surface generated using MSMS package in VMD with probe radius set to $1.5$ and density set to $9.5$. Figure \ref{fig:many_proteins} illustrates the biomolecular surfaces of the proteins mentioned above where each row corresponds to a specific protein. The first row corresponds to 1ajj, the second to 1bor, the third to 1mbg, and the fourth to 1sh1. In each row, the first surface is corresponding to the piecewise initial shape, the second surface is corresponding to the Gaussian initial shape and the third is corresponding to the MSMS surface.

\begin{figure}[h!]
\center
\begin{tabular}{ c c }
\subcaptionbox{}
{\includegraphics[scale=1.4,trim={3cm 0 3cm 0},clip]{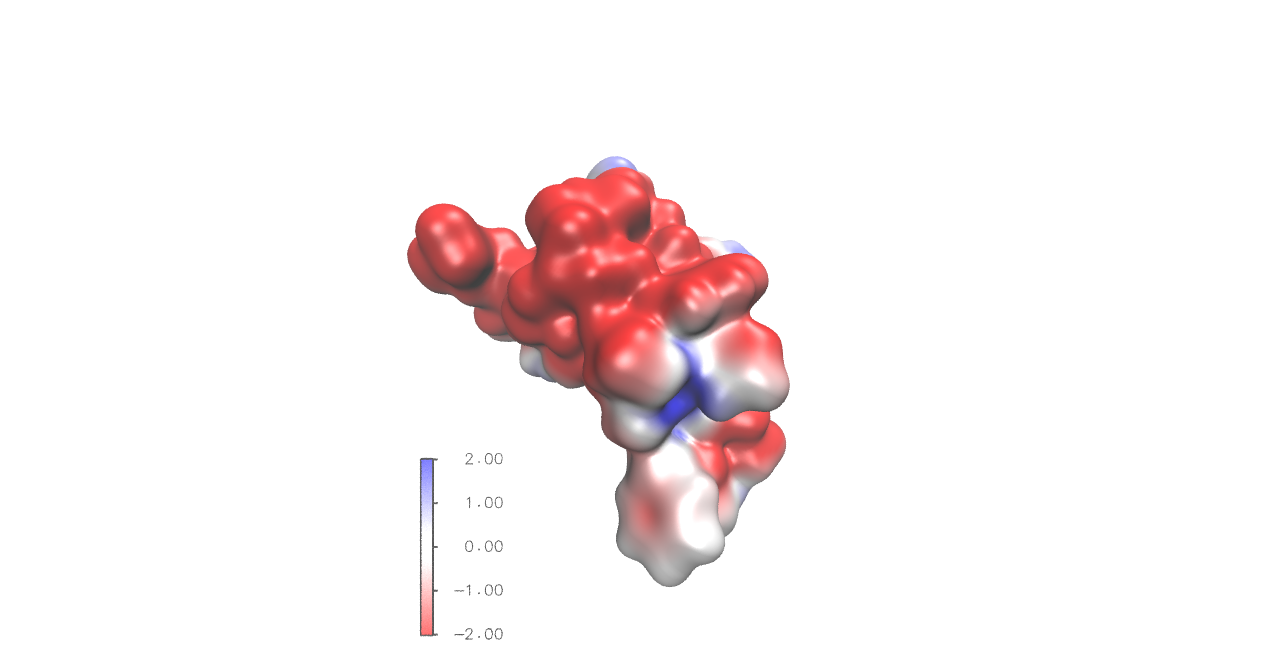}}
\subcaptionbox{}
{\includegraphics[scale=1.4,trim={3cm 0 3cm 0},clip]{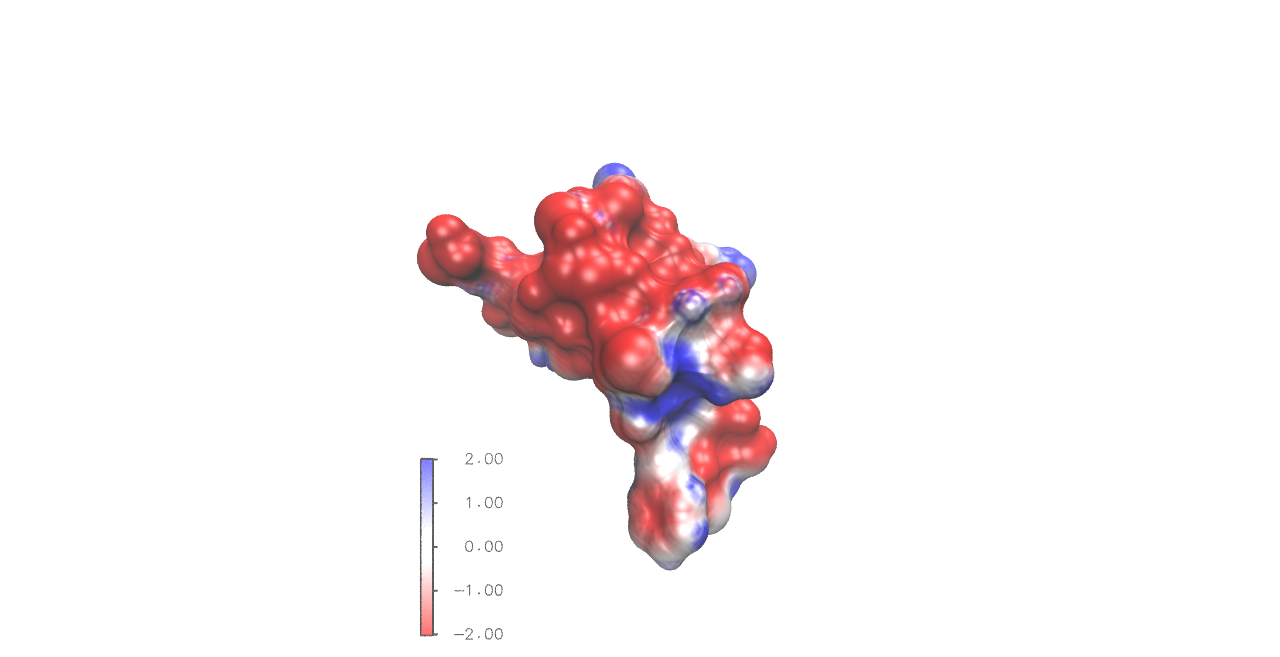}}
\end{tabular}
\caption{The electrostatic surface potentials of protein 1ajj mapped on two surfaces. (a) Surface generated by Clifford-Fourier transform method. (b) The MSMS surface. }
\label{fig:apbs_1ajj}
\end{figure}

\section{Applications}\label{sec:application}
The surfaces generated for a specific biomolecule represents its boundary region, which is regarded as the interface between the biomolecule region and the solvent region. The solvent-solute interface is essential in many models and applications such as electrostatic calculations \cite{jurrus2018improvements},  diffusion analysis \cite{zheng2011poisson}, and differential geometry-based solvation models \cite{wei2010differential}.

In this section, we show some applications on the biomolecular surfaces generated with our Clifford-Fourier transform method. We use the Poisson-Boltzmann model for electrostatic calculations. First, electrostatic surface potentials of some proteins are mapped to the surfaces generated with our Clifford-Fourier transform method and then they are compared with the mapping on surfaces generated by MSMS package. The electrostatic surface potentials are calculated using APBS method \cite{jurrus2018improvements} available in VMD. Then, the electrostatic solvation free energy is calculated and compared with three other methods of biomolecular surface generation. We calculate the energies using the match interface and boundary (MIB) method \cite{chen2011mibpb} and compared our results with MSMS, two different methods that are developed using flexibility and rigidity index (FRI) \cite{mu2017geometric}.  

\begin{figure}[h!]
\center
\begin{tabular}{ c c }
\subcaptionbox{}
{\includegraphics[scale=1.2,trim={2cm 0 2cm 0},clip]{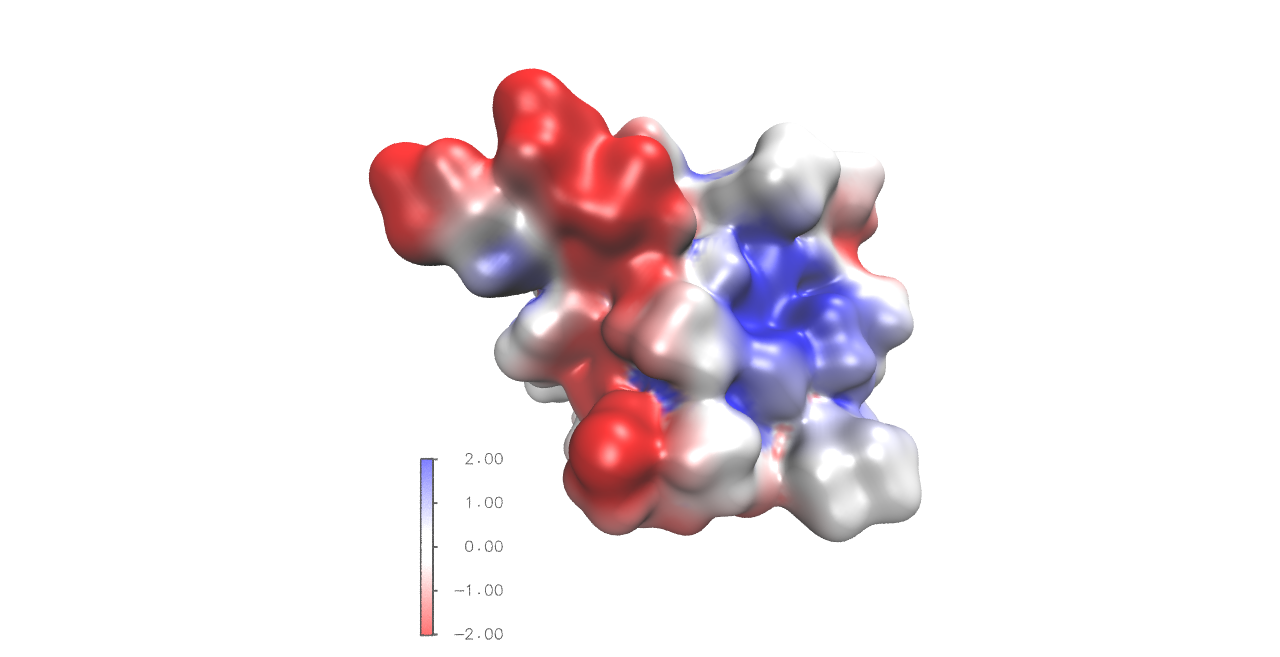}}
\subcaptionbox{}
{\includegraphics[scale=1.2,trim={2cm 0 2cm 0},clip]{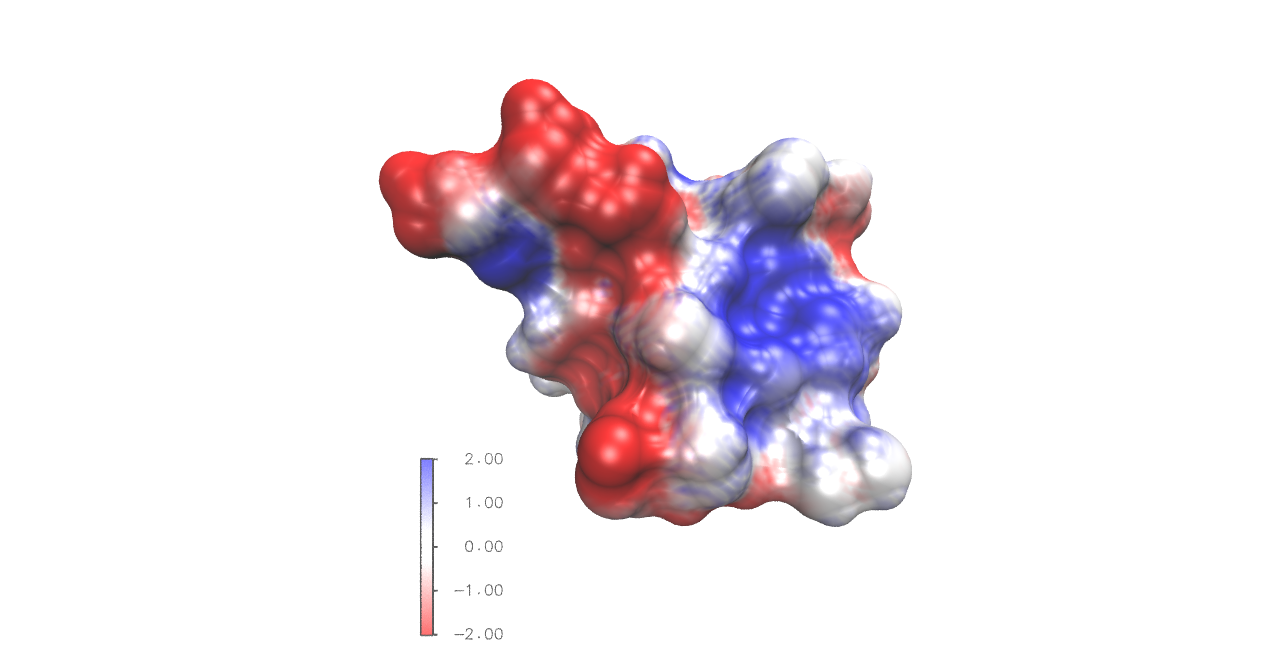}}
\end{tabular}
\caption{The electrostatic surface potentials of protein 1bor mapped on two surfaces. (a) Surface generated by Clifford-Fourier transform method. (b) The MSMS surface. }
\label{fig:apbs_1bor}
\end{figure}
\subsection{Electrostatic surface potentials}
The electrostatic surface potential is an important property for many applications in biology and biophysics \cite{zheng2012biomolecular}. It is essential in studies of drug design, protein-protein interactions, and other applications. The electrostatic potentials are calculated by solving the Poisson-Boltzmann equations. The calculation is done using a package available in VMD. After that, the potentials are mapped to the surface generated using our Clifford-Fourier transform method and to the SES surface generated by MSMS package in VMD. We validate our method by demonstrating the potentials on two proteins: 1ajj and 1bor. First, Figure \ref{fig:apbs_1ajj} shows the Clifford-Fourier transform surface on the left and the MSMS surface to the right where both figures show a very good match in surface potentials. Then, Figure \ref{fig:apbs_1bor} shows another example of surface potentials mapped to the Clifford-Fourier transform surface on the left and mapped to MSMS surface on the right. This latter figure illustrates also that the two surfaces matched very well in their surface potentials. On top of that, the Clifford-Fourier transform surfaces shown below do not have geometric singularities as can be verified from the figures. Moreover, Clifford-Fourier transform surfaces are more smooth in both proteins.

\begin{figure}[h!]
\centering
\begin{tabular}{ c|c|c|c|c } 
\hline
\multicolumn{5}{c}{Electrostatic solvation free energies (kcal/mol)} \\
\hline
Protein ID &  MSMS surface & FRI surface 1 &  FRI surface 2 & CFT surface \\
\hline
1ajj & -1100.754 &  -1155.158 &  -1258.784 &  -1154.039\\ 
1vii &  -862.865 &  -761.195 &  -846.414 &   -793.326\\ 
1bor &  -927.310 &  -1021.579 &  -1140.463 &  -833.647 \\ 
451c &  -1003.17  & -971.629 &  -1152.786  &  -946.969\\ 
1svr &  -1582.131 &  -1530.159  & -1700.400 &  -1518.964 \\ 
1uxc &  -1097.189 &  -975.583 &  -1075.070 &   -982.813 \\ 
1mbg &  -1340.086 &  -1329.440 &  -1412.048 &  -1226.894 \\ 
1ptq &  -800.130 &  -765.680 &  -866.256 &  -726.299 \\ 
1sh1 &  -729.626  & -648.315  & -785.820  &  -724.799 \\ 
2pde &  -1234.229 &  -863.702  & -990.056  &  -1192.055\\ 
1hpt &  -788.626 &  -768.424  & -873.410 &   -685.036\\ 
1a7m &  -2173.814 &  -2159.535  & -2492.651 &  -2196.988 \\  
1neq &  -1683.679  & -1661.077 &  -1811.667  &  -1552.394 \\ 
1r69 &  -1115.733 &  -983.129 &  -1087.688  &   -983.534 \\ 
1a2s &  -1868.827 &  -2216.464  & -2356.554  &  -1819.385\\  
2erl &  -894.960 &  -1127.617  & -1189.821 &  -856.024\\
1bbl &  -970.053  & -993.386  & -1052.617  &  -876.956 \\
1fca &  -1148.672  & -1427.109 &  -1534.746  & -1166.321 \\ 
1frd &  -2691.339 &  -2935.189 &  -3106.112  & -2438.468 \\ 
1bpi &  -1267.063  & -1170.959  & -1269.508  & -1128.051 \\ 
1a63 &  -2291.449 &  -2233.139 &  -2495.664  & -2213.839\\ 
\hline
\end{tabular}
\captionof{table}{Comparison of electrostatic solvation free energies for surfaces generated by MSMS, exponential FRI method with $\eta = 1.85$ Å, $\nu = 2$, and $\mu_0 = 1.5$, Lorentz FRI method with $\eta = 1.86$ Å, $\nu = 8$, and $\mu_0 = 1.5$, and Clifford-Fourier transform surface with grid size =$0.25$, isovalue= $0.9$, propagation time=$10^2$. }
\label{tab:table-electrostatic}
\end{figure}


\subsection{Electrostatic solvation free energy}
Now, we calculate the electrostatic solvation free energies of 21 proteins to validate our method of surface generation. These calculation are performed using the match interface and boundary method. To show the validity, the results are compared with three other methods discussed in the work of Mu \textit{et al}\cite{mu2017geometric}. Table \ref{tab:table-electrostatic} shows the electrostatic solvation free energies of molecular surfaces generated by MSMS package, exponential kernel based rigidity (FRI surface 1) with $\eta = 1.85$ Å, $\nu = 2$, and $\mu_0 = 1.5$, Lorentz kernel based rigidity (FRI surface 2) with $\eta = 1.86$ Å, $\nu = 8$, and $\mu_0 = 1.5$ and Clifford-Fourier transform with grid step = $0.25$, isovalue = $0.9$ and propagation time = $10^2$. The results of the three methods mentioned above were

\section{Concluding remarks}

Molecular surface generation is an important topic in computational biophysics and is crucial to the understanding of biological processes. A variety of computational methods have been developed for molecular surface generation, including those based on geometry, differential geometry, partial differential equation   (PDE) transform, level sets, etc. Therefore, molecular surface generation has been a research topic where biology meets physics, mathematics, and computer science. 

Geometric algebra has been widely applied to physics, computer vision, image analysis, and molecular biophysics, etc. However, it has not been used for molecular surface generation. This work introduces geometric algebra for molecular surface generation. More specifically, we utilize Clifford-Fourier transform (CFT), an important technique in geometric algebra, to define biomolecular surfaces.  
  
  We presented geometric algebra definitions of main calculus concepts such as integration and derivative. We also discussed in detail the $2$-dimensional CFT and $3$-dimensional CFT. We pointed to the fact that pseudoscalars in $\mathbb{R}_2$ are not commutative with multivectors in contrast to $\mathbb{R}_3$ which maintains the commutativity. This impacts the derivative property leading the $2$-dimensional CFT to not have a general rule. On the other hand, the $3$-dimensional CFT goes in parallel to the classical Fourier transform in this regard. Hence, we preferred the $3$-dimensional CFT. This choice is due to the importance of the derivative property in solving PDEs using CFT.   The PDE transform would then be implemented using CFT to generate molecular surfaces.

After introducing the Clifford-Fourier transform and PDE transform, we proposed our geometric algebra method of surface generation. We started by providing two cases of initial data: piecewise and Gaussian. Then, we conducted many experiments on an imaginary three-atom-molecule the real protein 1ajj to see the effect of changing the propagation times and the effect of changing the isovalues. These experiments showed that both initial cases were valid and able to generate good molecular surfaces. After that, molecular surfaces of real proteins were generated using both initial cases and compared to MSMS surfaces. Our surfaces showed superiority over MSMS surfaces due to the free of geometric singularities. To validate our method, we calculated electrostatic surface potentials and mapped them to our surfaces and to MSMS surfaces to visualize the electrostatic consistency and singularity free. Furthermore, we computed the electrostatic solvation free energies of 21 proteins on our surfaces and compared them to those from MSMS and FRI surfaces. Our surfaces showed very good results in terms of energies as well. One more feature of our method is that it gets the terminal state in a single-step approach which makes it time-efficient.

Finally, the proposed geometric algebra surface generation method opens the door for new multiscale methods for surface generation where different propagation times might be applied at once, and then the resulting surfaces can be combined to get the final surface. Also, the presented Clifford-Fourier transform has great potential to be exploited in biophysical problems like protein-protein docking and protein-ligand binding.

\section*{Acknowledgment}
This work was supported in part by NIH grant  GM126189, NSF grants DMS-2052983,  DMS-1761320, and IIS-1900473,  NASA grant 80NSSC21M0023,  Michigan Economic Development Corporation, MSU Foundation,  Bristol-Myers Squibb 65109, and Pfizer.
AA thanks Dr. Jiahui Chen for technical assistance.

\end{document}